\documentclass[10pt,reqno,final]{amsart}
\usepackage{amsfonts}
\usepackage[notcite,notref]{showkeys}
\usepackage{url,hyperref,multirow}
\usepackage{color}
\usepackage[dvipsnames]{xcolor}
\usepackage{cases}
\usepackage{subfigure}
\usepackage{float}
\usepackage{graphicx}
\usepackage{subcaption} % Include this package
\usepackage{epsfig,amsmath,bm,epstopdf}
\usepackage{amssymb,version,graphicx,fancybox,mathrsfs,pifont,booktabs,mathtools}
\usepackage{esint}
\usepackage{algorithm}
\usepackage{algpseudocode}
\usepackage[marginparwidth=2.5cm, marginparsep=3mm]{geometry}
\catcode`\@=11 \theoremstyle{plain}
\@addtoreset{equation}{section}   % Makes \section reset 'equation' counter.

\@addtoreset{figure}{section}
\renewcommand\thefigure{\thesection.\@arabic\c@figure}
\@addtoreset{table}{section}
\renewcommand\thetable{\thesection.\@arabic\c@table}
\newtheorem{thm}{Theorem}[section]  % Master counter
\newtheorem{assu}{Assumption}[section]
\newtheorem{cor}{Corollary}[section]
\newtheorem{prop}{Proposition}[section]
\newtheorem{lmm}{Lemma}[section]
\newtheorem{defn}{Definition}[section]

\newtheorem{exm}{Example}[section]
\newtheorem{remark}{Remark}[section]

\newenvironment{lemma}{\begin{lmm}}{\end{lmm}}

\newcommand \R {\mathbb{R}}

\newcommand \nmesh {N_\mathrm{mesh}}

\DeclareMathOperator{\sech}{sech}

\title[PODNO: Proper Orthogonal Decomposition Neural Operators] {PODNO: Proper Orthogonal Decomposition Neural Operators}
\author{Zilan Cheng$^*$}
\thanks{$^*$Division of Mathematical Sciences, School of Physical and Mathematical Sciences, Nanyang Technological University, Singapore. \texttt{zilan001@e.ntu.edu.sg}}

\author{Zhongjian Wang$^\dagger$}
\thanks{$^\dagger$Corresponding author. Division of Mathematical Sciences, School of Physical and Mathematical Sciences, Nanyang Technological University, Singapore. \texttt{zhongjian.wang@ntu.edu.sg}}

\author{Li-Lian Wang$^\ddagger$}
\thanks{$^\ddagger$Division of Mathematical Sciences, School of Physical and Mathematical Sciences, Nanyang Technological University, Singapore. \texttt{lilian@ntu.edu.sg}}

\author{Mejdi Azaiez$^\S$}
\thanks{$^\S$Bordeaux University, Bordeaux INP and I2M (UMR CNRS 5295), 33400 Talence, France. \texttt{azaiez@u-bordeaux.fr}}

\begin{document}
\begin{abstract}
In this paper, we introduce Proper Orthogonal Decomposition Neural Operators (PODNO) for solving partial differential equations (PDEs) dominated by high-frequency components. Building on the structure of Fourier Neural Operators (FNO), PODNO replaces the Fourier transform with (inverse) orthonormal transforms derived from the Proper Orthogonal Decomposition (POD) method to construct the integral kernel. Due to the optimality of POD basis, the PODNO has potential to outperform FNO in both accuracy and computational efficiency for high-frequency problems. From analysis point of view, we established the universality of a generalization of PODNO, termed as Generalized Spectral Operator (GSO). In addition, we evaluate PODNO's performance  numerically on dispersive equations such as the Nonlinear Schr\"odinger (NLS) equation and the Kadomtsev–Petviashvili (KP) equation.

\end{abstract}
\maketitle
% REQUIRED
% \begin{keywords}
% Neural operators, POD, dispersive PDEs, oscillatory solutions,  universal approximation theory 
% \end{keywords}

\section{Introduction}\label{intro}
Conventional methods, such as the popular Finite Difference \cite{quarteroni2008numerical}, Finite Element \cite{brenner2002mathematical}, Finite Volume \cite{leveque1992numerical}, Spectral \cite{gottlieb1977numerical,shen2011spectral} and Particle methods \cite{cottet2000vortex} 
are widely used to solve partial differential equations (PDEs) and provide high accuracy, particularly with refined meshes for well-defined problems. To deal with the high-dimensional and mesh-free problems, machine learning methods are designed for solving PDEs. By embedding physics into the loss function, PINNs \cite{raissi2019physics} can effectively learn even with limited data. The DeepRitz method \cite{yu2018deep} reformulates PDEs as variational problems by minimizing an energy functional derived from the weak form. Leveraging a neural network as the trial function, DeepRitz is particularly effective for solving elliptic PDEs and problems with variational structures. Inspired by the DeepRitz method, the works \cite{wang2020mesh,cui2022variational} explore the expressive power of deep neural networks in representing solutions to interface problems and problems with complex geometry, where efficient and accurate conventional solvers are not available. 

However, the above mechanism-driven methods typically focus on a single instance with the explicit form of the problem, which motivates the exploration of data-driven approaches. Neural operators have been proposed to solve a family of partial differential equations (PDEs) without requiring the explicit form of the problem, making them mesh-invariant and significantly faster during evaluation, though training may take more time. The Fourier Neural Operators (FNO) \cite{lifourier} is one such approach, designed to extend neural networks by learning mappings between infinite-dimensional spaces. FNO operates in the frequency domain using Fourier transforms, capturing global interactions efficiently. Variations of FNO, such as Graph Neural Operator \cite{li2020neural} and Multipole Graph Neural Operator(MGNO) \cite{li2020multipole}, extend FNO for complex geometry by incorporating graph networks and multipole expansions to capture long-range dependencies. The Orthogonal Polynomial Neural Operator (OPNO) \cite{liu2024render} framework extends FNO to non-periodic boundary conditions, offering competitive accuracy. The Semiperiodic Fourier Neural Operator 
(SPFNO)~\cite{liu2023spfno} further improves boundary condition handling by using trigonometric bases, ensuring errors reach machine precision. Additionally, the Laplace Neural Operator (LNO) 
 \cite{cao2024laplace} outperforms FNO by effectively capturing both transient and steady-state responses, making it more suitable for non-periodic signals and complex dynamics in time-dependent systems. Another data-driven method for solving PDEs is the Random Feature Model \cite{nelsen2021random}, which learns the mapping between the input and output spaces by projecting input functions into a high-dimensional space using randomly generated features. While more efficient and flexible than FNO, it may struggle with capturing global dependencies in complex spatial problems. Deep Operator Network (DeepONet) \cite{lu2021learning} is also a powerful framework for learning operators that map functions to other functions. Physics-Informed 
DeepONets \cite{wang2021learning} extend DeepONets by embedding physical laws. Unlike FNO, which relies on the Fourier basis, DeepONet uses an encoder-decoder structure to map between function spaces, offering greater flexibility. While FNO emphasizes spectral methods, DeepONet's adaptable architecture suits diverse problems.

The aforementioned operator learning methods are primarily tested on diffusion-dominated problems, such as the Darcy problem and Navier-Stokes equations. While high-frequency components may initially be present, these problems are typically characterized by a dominance of low-frequency patterns, as the energy tends to dissipate and spread out over time. Solving the complex dispersive equations, for example, the Nonlinear Schr\"odinger equation (NLS) \cite{kato1987nonlinear} and the Kadomtsev-Petviashvili (KP) equations \cite{kadomtsev1970stability} are also of significant importance. The NLS equation is a classical field equation primarily applied to the propagation of light in nonlinear optical fibers \cite{fibich2015nonlinear}, as well as to Bose-Einstein condensates in highly anisotropic traps \cite{griffin1996bose}. The KP equation, a two-dimensional generalization of the Korteweg–de Vries (KdV) equation, is crucial for studying the propagation of nonlinear waves in shallow water \cite{johnson2003classical}, as well as sound waves in
ferromagnetic media \cite{klein2007numerical} and nonlinear matter wave
pulses in Bose-Einstein condensates \cite{huang2003two, jones1982motions}. These equations often exhibit complex, nonlinear behaviors that are challenging to solve analytically, especially in high-dimensional or turbulent regimes. Efficient and accurate numerical methods are essential for understanding the dynamics of such systems, which makes the development of advanced solvers like neural operators particularly valuable for tackling these nonlinear problems with high precision and generalizability.

Due to the high demand for solving these equations that are rich in dynamics and high-frequency-dominated, we propose Proper Orthogonal Decomposition Neural Operators (PODNO) in this paper. While the standard FNO effectively captures global interactions in the frequency domain, it struggles with the complexity and nonlinearity inherent in these equations, particularly at high frequencies. By replacing the Fourier basis with some sparse basis, POD basis \cite{berkooz1993proper}, in the essential kernel integration part, PODNO captures the highest energy of the data space with both high frequency and low frequency, and thus enhances its ability to manage the complex wave dynamics typical of NLS and KP equations. Through extensive experiments, we demonstrate that our approach outperforms the standard FNO in terms of both accuracy and computational efficiency for high-frequency problems, offering a more robust and reliable tool for solving these complex PDEs. 

It is noteworthy that Reduced Order Modeling (ROM) has been applied to operator learning in PCA-Net \cite{bhattacharya2021model}, which projects input and output functions onto low-dimensional subspaces spanned by dominant PCA modes extracted from the training data. However, PCA-Net employs a relatively simple deep neural network (DNN) architecture to map the coefficients, which limits its applicability to problems where a small number of modes- typically corresponding to smooth or low-rank solutions- are sufficient to capture the dominant behavior of the solution operator. After that, the POD-DeepONet framework combines POD with DeepONet \cite{demo2023deeponet}. While POD captures the dominant modes of a dataset, DeepONet learns the residual between low-fidelity POD approximations and high-fidelity solutions. However, this model largely relies on the initial low-fidelity POD model, which may restrict its performance on highly complex or nonlinear problems. In contrast, PODNO can directly operate in function space, effectively capturing more intricate patterns.

The remainder of this paper is organized as follows: Section \ref{neural operators} introduces neural operators (NO). In Subsection \ref{Problem setting}, we outline the problem setting for neural operators, while Subsection \ref{FNO structure} discusses the structure of FNO, with particular emphasis on its essential kernel. In Section \ref{podno}, we propose PODNO through four key steps. Subsection \ref{GSO} introduces the concept of Generalized Spectral Operators (GSO), whose kernel is defined based on an orthonormal basis. Subsection \ref{POD method} explains the generation of the POD basis in our problem setting and describes operations on the integration kernel layers. Subsection \ref{podno_alg} then integrates these components and presents the PODNO algorithm. Finally, Subsection \ref{universality} establishes the universal approximation theorem for GSO. Section \ref{numerical results} presents numerical results. Subsection \ref{setting} defines the model setup, and Subsection \ref{darcy results} validates PODNO for steady-state problems. Subsections \ref{nls results} and \ref{kp results} explore applications to the NLS equation and the KP equation, respectively. In Subsection \ref{pod v.s. podno}, we compare the numerical results of PODNO with those of the POD-accelerated splitting method for the NLS equation. Finally, Subsection \ref{ablation} presents an ablation study for both PODNO and the ground truth solvers, assessing parameter dependencies in implementation and evaluating algorithmic robustness. The appendices provide additional theoretical insights. Appendix \ref{proof} contains the proof for GSO's universal approximation theorem. In Appendix \ref{app:pod}, we present the POD method in discrete setting, and in  Appendix \ref{dispersive} we discuss the high-frequency properties of the NLS and KP equations.  
\section{Neural Operators}\label{neural operators}
In this section, we first outline the problem setting for neural operators in Subsection \ref{Problem setting}, which naturally leads to the discussion of the FNO structure in Subsection \ref{FNO structure}.
\subsection{General setup in operator learning}\label{Problem setting}
 A neural operator is a neural network model of a mapping between two infinite-dimensional spaces using a finite collection of observed input-output pairs. Let $\mathcal{A}$ and $\mathcal{U}$ denote the input and output function spaces defined on the bounded domains $\Omega \subset \mathbb{R}^d$ and $\Omega' \subset \mathbb{R}^{d'}$, respectively. For simplicity, we assume $\Omega=\Omega'$ in the following discussion. The functions in $\mathcal{A}$ are $\mathbb{R}^{d_a}$-valued, while the functions in $\mathcal{U}$ are $\mathbb{R}^{d_u}$-valued. The mapping $\mathcal{G}$ between these function spaces can be expressed as follows \cite{lifourier}:  
\begin{equation*}
    \mathcal{G}: \mathcal{A}\rightarrow \mathcal{U}.
\end{equation*}
Suppose we have observations $\{a_j,u_j\}^M_{j=1}$ where $a_j\sim\mu$ is an i.i.d. sequence from the probability measure $\mu$ supported on $\mathcal{A}$ and $u_j=\mathcal{G}(a_j)$ is possibly corrupted with noise. Neural operators aims to build $\mathcal{N}$, an approximation of $\mathcal{G}$ by constructing a parametric map 
\begin{equation*}
    \mathcal{N}: \mathcal{A}\rightarrow\mathcal{U}.
\end{equation*}
We denote the network with latent parameter $\theta$ from some finite-dimensional parameter space $\Theta$ as $\mathcal{N}_\theta$. Then we can define a cost functional $C: \mathcal{U}\times\mathcal{U}\rightarrow\mathbb{R}$ measuring the discrepancy between $\mathcal{N_\theta}$ and the ground truth operator $\mathcal{G}$ and then aim to find the latent parameter to minimize the cost:
\begin{equation}\label{theta}
    \hat{\theta}=\arg\min _{\theta \in \Theta} \mathbb{E}_{a \sim \mu}[C(\mathcal{N_\theta}(a),\mathcal{G}(a))]
    \approx \arg\min _{\theta \in \Theta} \frac{1}{M} \sum_{j=1}^{M}C(\mathcal{N}_\theta(a_j),u_j).
\end{equation}
The cost function is typically defined in terms of the $L^2$ norm.
\subsection{FNO structure}\label{FNO structure}
\begin{figure}[!h]
    \centering
    \includegraphics[width=0.9\textwidth]{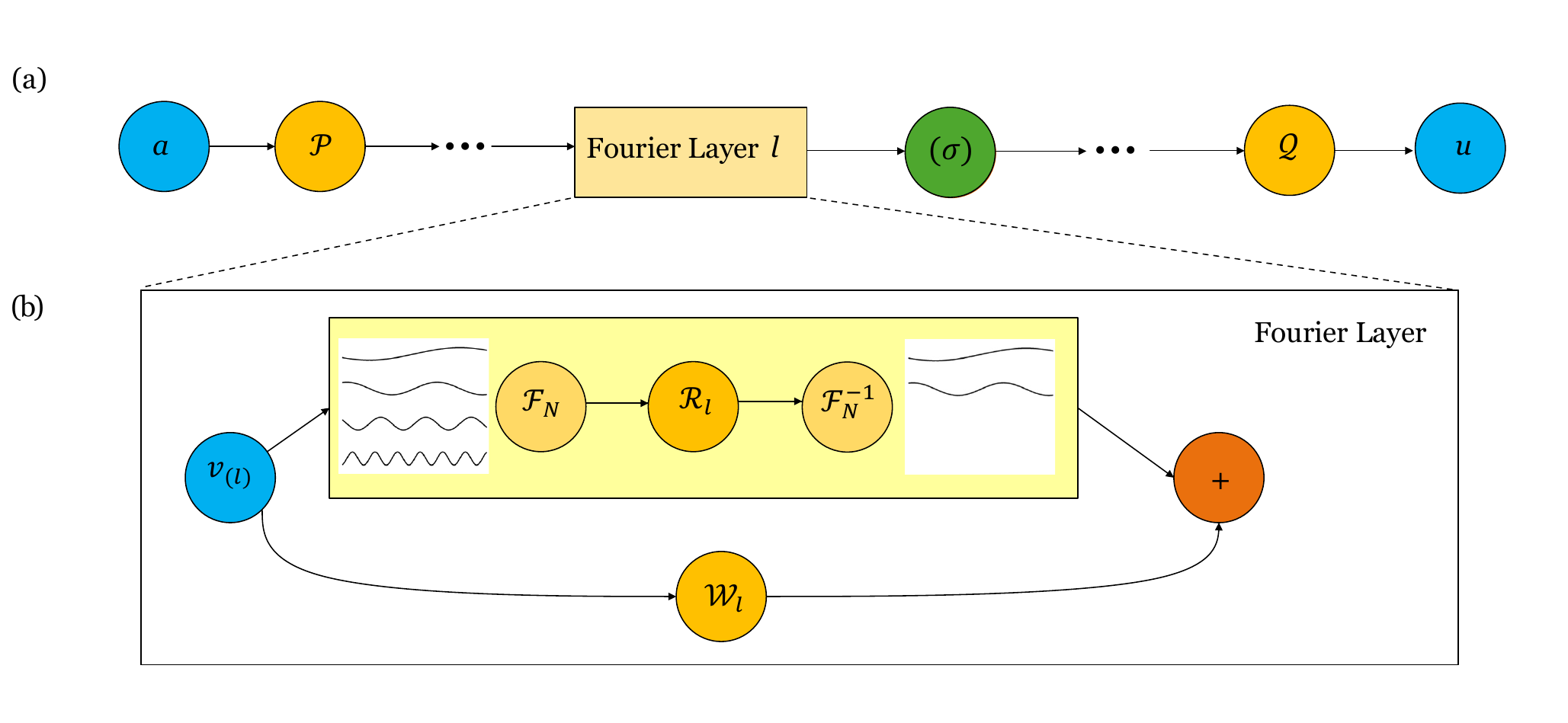}
    \caption{FNO. (a) The architecture of FNO consists of 3 parts: dimension lifting $\mathcal{P}$; a sequence of Fourier layers each followed by a nonlinear activation function except the final one; and dimension reduction $\mathcal{Q}$. (b) Each Fourier layer is given by $\mathscr{F}_N^{-1}[\mathcal{R}_l\circ\mathscr{F}_N[v_{(l)}]]+\mathcal{W}_l[v_{(l)}]$.}
    \label{FNO}
\end{figure}
\cite{lifourier} proposed the FNO structure (see Fig.~\ref{FNO} (a)) of the approximation mapping $\mathcal{N}$: 
\begin{equation*}
    \mathcal{N}: \mathcal{Q} \circ(\mathcal{W}_{L-1}+\mathcal{K}_{L-1}) \circ\dots\circ\sigma\circ(\mathcal{W}_1+\mathcal{K}_1) \circ\sigma\circ(\mathcal{W}_0+\mathcal{K}_0) \circ \mathcal{P},
\end{equation*}
where $L$ denotes the number of Fourier layers and $\sigma$ is pointwisely applied nonlinear activation function. The structure includes three parts \cite{kovachki2023neural}:

\textbf{Dimension lifting}. The dimension lifting  $\mathcal{P}$ to map the $d_a$ dimensional vector field input $a$ to its first $d_v$ dimensional vector field hidden representation $v_0$ with a linear transform:
    \begin{equation}\label{eq_P}
    v_0(x)=\mathcal{P}[a](x)=P_1a(x)+P_21(x), \quad \forall x \in \Omega,
    \end{equation}
    where $P_1\in{\mathbb{R}^{d_v}\times\mathbb{R}^{d_a}}$, $P_2\in{\mathbb{R}^{d_v}}$ are the weight matrix and the bias vector learned during training. Hereafter, $1(x)$ denotes the constant function on $\Omega$. Then the operation $P_21$ produces a constant vector-valued function with the same components as $P_2$. 

\textbf{Dimension reduction}. Similar to Dimension lifting $\mathcal{P}$, the Dimension reduction $\mathcal{Q}$ layer  projects the last hidden representation $v_L$ to the output function $u$:
    \begin{equation}\label{eq_Q}
    u(x)=\mathcal{Q}[v_{(L)}]=Q_1v_{(L)}(x)+Q_21(x), \quad \forall x \in \Omega,
    \end{equation}
    where $Q_1\in\mathbb{R}^{d_u\times d_v}$ and $Q_2\in\mathbb{R}^{d_u}$ are the weight matrix and the bias vector learned during training.
   
    \textbf{Fourier Layers}. For each Fourier layer $l = 0, \dots, L - 2$, the hidden representation $v_{(l)}$ is mapped to $v_{(l+1)}$ by applying the sum of a local linear operator $\mathcal{W}_l$ and a non-local integral kernel operator $\mathcal{K}_l$, followed by a fixed nonlinear activation function $\sigma$. In the final layer ($l = L - 1$), the activation function is omitted.
    
    The linear operator is defined by 
    \begin{equation}\label{eq_W}
    \mathcal{W}_l[v_{(l)}](x)=W_lv_{(l)}(x)+b_l1(x), \quad \forall x\in\Omega,  
    \end{equation}
    where $W_l\in{\mathbb{R}^{{d_v}\times {d_v}}}$ and $b_l\in{\mathbb{R}^{d_v}}$ are the weight matrix and the bias vector learned during training. 
    
    In practice, one can use a convolutional layer with a $1\times 1$ convolution kernel to implement the pointwise linear transforms introduced in \eqref{eq_P}, \eqref{eq_Q} and \eqref{eq_W}. 
    
 Now we turn to the integral kernel operator $
 \mathcal{K}_l$ which has the following formal representation: 
\begin{equation}\label{kernel0}
    \mathcal{K}_l[v_{(l)}](x)=\int_{\Omega}\kappa_l(x,y)v_{(l)}(y) \mathrm{d}y,
\end{equation}
where $\kappa_l$ is the kernel function of the integral operator and of size $d_v\times d_v$. 

If one assumes the kernel is translation-invariant, i.e., $\kappa_l(x,y)=\kappa_l(x-y)$,  the operation in \eqref{kernel0} can be reduced to the following, due to the convolution theorem:
\begin{equation*}\label{kernel1}
\mathcal{K}_l[v_{(l)}](x) =\int_{\Omega}\kappa_l(x-y)v_{(l)}(y) \mathrm{d}y=\mathscr{F}^{-1}\big[\mathscr{F}[\kappa_l]  \cdot \mathscr{F}[v_{(l)}]\big](x),
\end{equation*}
where $\mathscr{F}$ and $\mathscr{F}^{-1}$ represent the Fourier transform and its inverse, respectively.

Now we assume that the operations are conducted in the discretized case, that is to say, $\kappa_l$ and $v_{(l)}$ are no longer continuous functions over the domain of $x$ in $\mathbb{R}^{d_v}$, but instead matrices of size $\mathbb{R}^{d_v\times \nmesh}$, where $\nmesh$ is the size of discretization. Then we approximate the operation \eqref{kernel} with discrete Fourier transform restricted to the first $N$ Fourier modes $\Phi_N:=\{\phi_k\}_{k=1}^N$,
\begin{align}\label{kernel}
    \mathcal{K}_l[v_{(l)}](x)\approx\mathscr{F}_N^{-1}[{\mathcal R}_l \cdot \mathscr{F}_N[v_{(l)}]](x),
\end{align}
where $\mathscr{F}_N$ and $\mathscr{F}^{-1}_N$ represent the Fourier transform  and its inverse with truncation to $\Phi_N$; $\mathcal{R}_l$ denotes the linear operation on  $\Phi_N$, which can be represented by an ${d_v \times d_v\times N\times N}$ tensor $R_{l}$.
 More precisely, we start with ${v}_{(l)}\in\mathbb{R}^{d_v\times \nmesh}$, then  $\hat{v}_{(l)}=\mathscr{F}_N(v_{(l)})\in\mathbb{R}^{d_v\times N}$, and
         \begin{equation}\label{rl_matrix_extend}(\mathscr{F}_N^{-1}\mathcal{R}_l[\hat{v}_{(l)}])_j=\sum_{i=1}^{d_v}\sum_{k=1}^{N}\sum_{m=1}^{N}R_{l;ijkm}(\hat{v}_{(l)})_{ik}\phi_m,\quad \forall j=1,\cdots,d_v. 
    \end{equation}
    
    \begin{remark}The operator $\mathcal{R}_l$ is applied mode-wise in practical implement of FNO. Namely $R_l$ becomes an $d_v\times d_v\times N$ tensor, and \eqref{rl_matrix_extend} becomes,  \begin{equation*}(\mathscr{F}_N^{-1}\mathcal{R}_l[\hat{v}_{(l)}])_j=\sum_{i=1}^{d_v}\sum_{k=1}^{N}R_{l;ijk}(\hat{v}_{(l)})_{ik}\phi_k,\quad \forall j=1,\cdots,d_v. \label{rl_matrix}
    \end{equation*}
    
During the implementation of the later proposed PODNO, we also adapt this simplification. In our proof of universality (Step 1 of Lemma \ref{gn gso}), we utilize the more general formulation \eqref{rl_matrix_extend} to recover the basis in $\Phi_N$ through selection of $R_l$. While in Lemma 7 (Step 1) of the original proof of universality of FNO \cite{kovachki2021universal}, $b_l$ in the linear operator \eqref{eq_W} is generalized to an learnable non-constant function to recover the basis, which may impose difficulties in discretization.
    \end{remark}

In the integral kernel operator $\mathcal{K}_l$ \eqref{kernel} of FNO, to highlight the function's global and essential features, high-frequency components (higher modes in the Fourier transform) are excluded in each iteration, i.e. $N\ll \nmesh$. The discarded information is compensated for by using a linear transformation $W_l$ in \eqref{eq_W}. This design brings two key limitations to FNO: difficulty with high-frequency-dominated problems due to mode truncation in kernel integration and reliance on periodic boundary conditions due to the application of the Fourier transform. Additionally, leveraging the FFT for improved computational efficiency typically requires the domain to be rectangular. To address boundary condition constraints, certain spectral bases have been employed in kernel design to handle non-periodic problems \cite{liu2024render, liu2023spfno}.

However, these variants are not well-suited for problems dominated by high-frequency components and often entail higher computational complexity. Complex dispersive equations, such as nonlinear Schr\"odinger (see Subsection \ref{nls results} and Appendix \ref{nls}) and KP equations (see Subsection \ref{kp results} and Appendix \ref{kp}), exhibit intricate high-frequency patterns, which demand more effective schemes for operator learning. This limitation motivates our proposal for PODNO, which will be discussed in detail in the next section.
\section{Generalized Spectral Operators and PODNO}\label{podno}
In this section, we first propose a more general forward propagation framework inspired by FNO, which we denote as Generalized Spectral Operators (GSO) in Subsection \ref{GSO}. PODNO is subsequently introduced as a special case that leverages the optimality of the POD basis under the $L^2$ metric. Subsection \ref{POD method} explores the construction of the POD basis and the design of the POD layers. The PODNO algorithm is subsequently outlined in Subsection \ref{podno_alg}. Finally, we present the universality of GSO in Subsection \ref{universality}.

\begin{figure}[h!]
    \centering
    \includegraphics[width=0.85\textwidth]{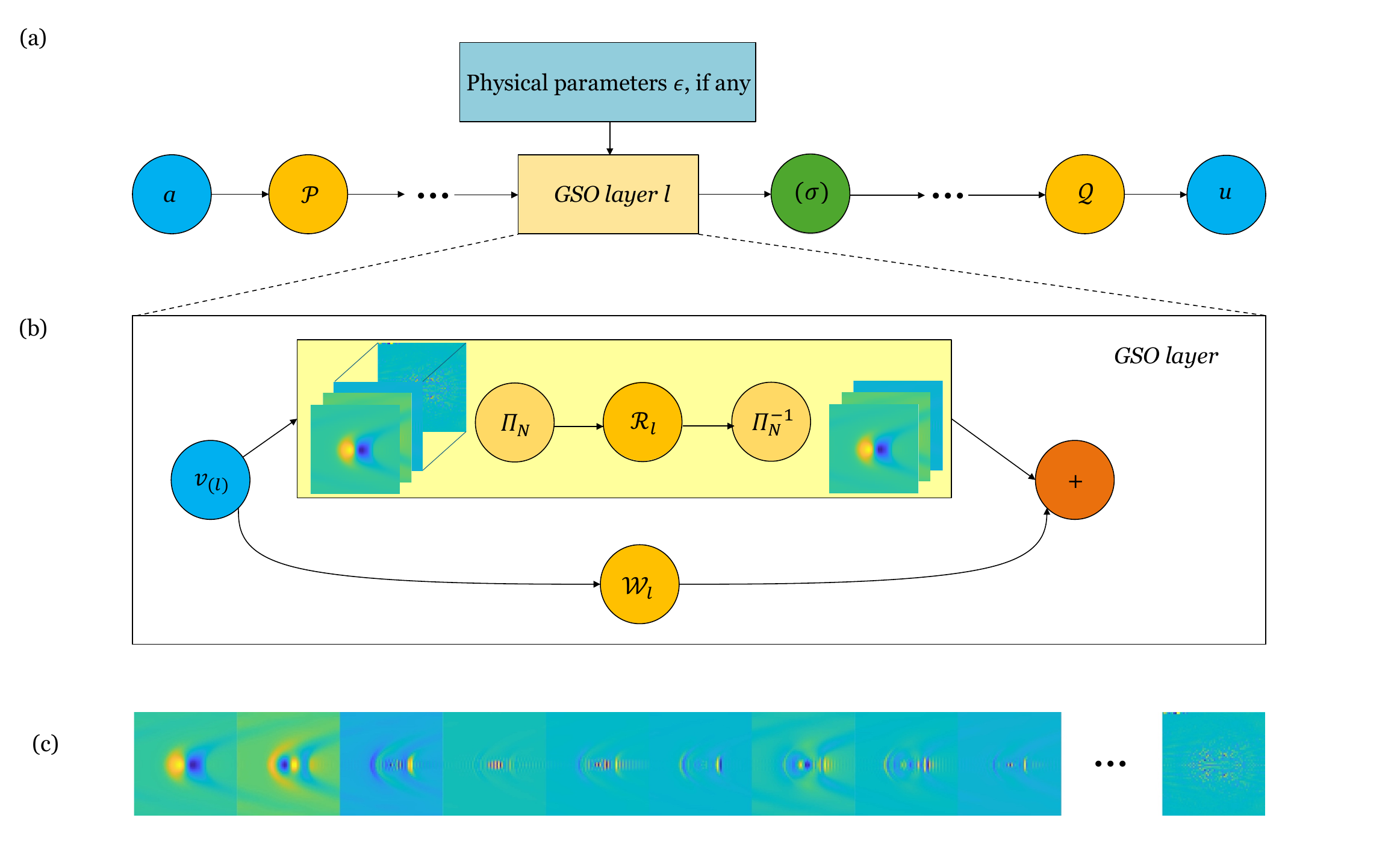}
    \caption{Generalized Spectral Operators (GSO). (a) The architecture of GSO comprises three main components: a dimension lifting $\mathcal{P}$; a sequence of GSO layers $\mathcal{L}_l$ that incorporate physical parameters, and apply a nonlinear activation function $\sigma$ except in the final layer; and a dimension reduction $\mathcal{Q}$. (b) The GSO layers $\mathcal{L}_l=\varPi_N^{-1}[\mathcal{R}_l\cdot\varPi_N [v_{(l)}]]+\mathcal{W}_l[v_{(l)}]$ . (c) The POD basis for a family of KP equations, which will be explicitly described in Subsection \ref{kp results}, provides an orthonormal set that can serve as the basis for the GSO layers.}\label{Generalized FNO}
\end{figure}
\subsection{Generalized spectral operators}\label{GSO}
Given any transform $\varPi_N$ corresponding to a (weighted) truncated orthonormal basis, we can define an essential layer of the form
\begin{equation*}
    v_{(l)}\mapsto\varPi_N^{-1}[\mathcal{R}_l\circ\varPi_N [v_{(l)}]]+\mathcal{W}_l[v_{(l)}],
\end{equation*}
where $\mathcal{R}_l$ denotes the diagonal operator and can be represented by tensors as in \eqref{rl_matrix_extend} in Fourier layers. We will subsequently define neural networks termed GSO based on this structure. This framework encompasses not only operators constructed from spectral transform-based kernel integration, but also those employing alternative orthonormal transforms, such as the POD basis to be introduced later.
The refined architecture of GSO is illustrated in Fig.~\ref{Generalized FNO}~(a)-(b), with detailed information on the forward propagation provided therein:

\textbf{Dimension lifting/ reduction.} The initial function $v_0(x)$ to be processed by the essential GSO layers is obtained by applying a dimension lifting operation to the input function: $v_0(x)=\mathcal{P}[a](x)$, where the lifting operator $\mathcal{P}$, as defined in \eqref{eq_P}, maps the function $a$ into a higher-dimensional space.

After the final GSO layer, the resulting high-dimensional representation $v_{(L)}(x)$ is projected back to a lower-dimensional output using the projection operator $\mathcal{Q}$ as defined in \eqref{eq_Q}:  $u(x)=\mathcal{Q}[v_{(L)}](x)$.

\textbf{GSO layers.} The core of the algorithm consists of $L$ iterative layers. Each layer $l = 0, \dots, L - 1$ involves the following components:
\smallskip
\begin{itemize}
    \item \textbf{Kernel integration operator.} The function $v_{(l)}(x)$ is first transformed into the coefficient space using the orthonormal transform $\varPi_N$, yielding $\varPi_N v_{(l)}$. A diagonal operator $\mathcal{R}_l$, representing the restriction of the orthonormal transform of the kernel to the first $N$ modes, is then applied: $\mathcal{R}_l\circ\varPi_N[v_{(l)}](x)$ like \eqref{rl_matrix_extend} in FNO. The result is subsequently mapped back to the physical space using the inverse transform $\varPi^{-1}_N$, yielding $\varPi^{-1}_N\circ \mathcal{R}_l\circ\varPi_N[v_{(l)}](x)$.
    \medskip 
    \item \textbf{Embedding operator.}  
        Given a physical parameter $\epsilon$, it is first transformed via two consecutive linear layers with an activation in between: 
        \begin{equation*}
        \hat{\epsilon}_l=\hat{\mathcal{W}}^{(2)}_l\sigma (\hat{\mathcal{W}}^{(1)}_l\epsilon), 
        \end{equation*} where $\hat{\mathcal{W}}^{(1)}_l$ and $\hat{\mathcal{W}}^{(2)}_l$ are affine maps and $\sigma$ is the nonlinear activation. The transformed parameter $\hat{\epsilon}_l$ is then broadcast as a constant function: $\hat{\epsilon}_l\mapsto\hat{\epsilon}_l1(x)$. In case when the learning objective operator is not parameter dependent, we take $\hat{\epsilon}_l\equiv0$.
         \medskip
        \item \textbf{Linear operator.} The linear transformation $\mathcal{W}_l[v_{(l)}](x)$ plays the same role as  \eqref{eq_W}.  \medskip
        \item \textbf{Summation}. The outputs from the kernel integration, linear compensation, and embedding operator are aggregated to form the input to the next operations: $[\varPi_N^{-1}[\mathcal{R}_l\circ\varPi_N [v_{(l)}]]+\mathcal{W}_l[v_{(l)}]+\hat{\epsilon}1](x)$, where $1(x)$ denote the constant function taking value $1$.
\end{itemize}
 \medskip

\textbf{Activation function.} For GSO layers $l = 0, \dots, L - 2$, a nonlinear activation function $\sigma$ is applied to the summation as follows 
\begin{equation*}
    v_{(l+1)}(x)=\sigma\big[\varPi_N^{-1}[\mathcal{R}_l\circ\varPi_N [v_{(l)}]]+\mathcal{W}_l[v_{(l)}]+\hat{\epsilon}_l1\big](x),           
\end{equation*}
       and there is no activation for 
       the final layer $l=L-1$. 

         \medskip
% \textbf{Dimension reduction.} After the final GSO layer, the resulting high-dimensional representation $v_{(L)}(x)$ is projected back to a lower-dimensional output using the projection operator $\mathcal{Q}$ as defined in \eqref{eq_Q}:  $u(x)=\mathcal{Q}[v_{(L)}](x)$.

Based on the forward propagation, we define the fundamental GSO as follows:
\begin{defn}[GSO]\label{podno chain}
Define the GSO as the mapping:
\begin{equation}\label{gso formula}
    \mathcal{N}(a):=\mathcal{Q}\circ\mathcal{L}_{L-1}\circ\sigma\circ\mathcal{L}_{L-2}\circ\dots\circ\mathcal{L}_1\circ\sigma\circ\mathcal{L}_0\circ \mathcal{P}[a],
\end{equation}
where $\mathcal{P}$ and $\mathcal{Q}$ are the lifting and projection operators, respectively. The GSO layers $\mathcal{L}_l$ take the following form:
\begin{equation}\label{nonlinear layer}
\begin{aligned}
    \mathcal{L}_l(v_{(l)})(x)&=\mathcal{W}_l[v_{(l)}](x)+\mathcal{K}_l[v_{(l)}](x)+\hat{\epsilon}1(x)\\
    &=W_lv_{(l)}(x)+b_l(x)+\varPi_N^{-1}[\mathcal{R}_l\circ\varPi_N[v_{(l)}]](x)+\hat{\epsilon}1(x).
\end{aligned}
\end{equation}
Here $\mathcal{W}_l[v_{(l)}](x)=W_lv_{(l)}(x)+b_l(x)$ defines a pointwise affine mapping and $\mathcal{R}_l$ defines a non-local linear mapping after the orthonormal transform $\varPi_N$ to the coefficients space, $\hat{\epsilon}1(x)$ helps to embed the physical parameter.
\end{defn}

We note a slight misuse of notation in using the kernel integration operator $\mathcal{K}_l$, which here refers to operations induced by general orthonormal transforms— not limited to the Fourier transform.

Although different bases could be chosen for the kernel integration of GSO (for instance, FNO), we notice that the best orthonormal basis under the $L^2$ metric is the POD basis (to be discussed in Subsection \ref{POD method}). The principal POD modes grasp almost all energy, as we can see in Fig.~\ref{Generalized FNO}~(c), the first basis explicitly captures the peak of the solutions for a family of KP equations (which we describe explicitly in Subsection \ref{kp results}), while the subsequent basis represents other significant features. This is why we aim to drop out of the non-principal modes. Therefore, we propose PODNO, which is based on the POD basis. Unlike the Fourier basis and other predefined spectral bases, the POD basis is data-dependent. Therefore, a key step is to identify the POD basis that enables the construction of the transform $\varPi_N$ and its inverse. 

\subsection{POD basis generation and  POD layers}\label{POD method}
As a model-reduction technique in various fields such as fluid dynamics, signal processing, and data analysis, the POD method \cite{berkooz1993proper} provides the best approximation of a dataset or a field with a given number of basis functions under a specific Sobolev metric, usually the $L^2$ metric.

We start with denoting the solution of interest as a random field $\mathcal{V}\in P(L^2(\Omega))$. The goal of POD methods is to determine an orthonormal basis $\{\phi_k\}^N_{k=1}$ such that the projection of the random field $\mathcal{V}$ onto the subspace spanned by $\{\phi_k\}^N_{k=1}\subset L^2(\Omega)$, $\mathcal{V}_N$, minimizes the error:
\begin{equation}\label{POD_continuous}
    \min_{\phi_1,\phi_2,\dots,\phi_N}\mathbb{E}_P[\|\mathcal{V}-\mathcal{V}_N\|_{L^2}^2].
\end{equation}
The solution to \eqref{POD_continuous} can be obtained by finding the first $N$ generalized eigenfunctions $\{(\lambda_k,\phi_k)\}_{k=1}^N$, in the homogeneous Fredholm integral equation \cite{smithies1958integral}:
\begin{equation}\label{Fredholm integral}
    \int_\Omega \mathcal{C}(x,y)\phi_k(x)\mathrm{d}x=\lambda_k\phi_k(y), %\quad \forall k,
\end{equation}
where $\mathcal{C}(x,y)=\mathbb{E}_P[ \mathcal{V}(x)\mathcal{V}(y)]$ is the covariance kernel and the eigenvalues $\lambda_k$ are ordered as  $\lambda_1 \geq \lambda_2 \geq \dots \geq 0$.

As in the method of snapshots \cite{kunisch2001galerkin}, we consider a set of $M$ realizations (samples) of $\mathcal{V}$: $\{v^{(m)}(x)\}_{m=1}^M$. Then the objective function in \eqref{POD_continuous} can be approximated by
\begin{equation*}
\begin{aligned}
    \mathbb{E}_P[\|\mathcal{V}-\mathcal{V}_N\|_{L^2}^2]& \approx \frac{1}{M} \sum_{m=1}^{M} \int_\Omega \big| v^{(m)}(x) - v^{(m)}_N(x)\big|^2\mathrm{d}x\\
    & =\frac{1}{M} \sum_{m=1}^{M}\int_\Omega\Big|  v^{(m)}(x) - \sum_{k=1}^N \langle v^{(m)}, \phi_k \rangle_{L^2} \, \phi_k(x)\Big|^2\mathrm{d}x.
\end{aligned}
\end{equation*}

Since the POD seeks to minimize the projection error of the data in the $L^2$ norm, the optimal basis functions ${\phi_k}$ lie in the span of the snapshots $\{v^{(m)}\}_{m=1}^M$. Each POD basis function is thus written as a linear combination:
\begin{equation*}
\phi_k(x) = \sum_{m=1}^M w_m^{(k)} v^{(m)}(x), \quad k = 1, \dots, N,    
\end{equation*}
where the coefficients $w_m^{(k)}$ are to be determined. To find these coefficients, we  define the correlation matrix $C \in \mathbb{R}^{M \times M}$ with entries
\begin{equation*}
C_{ij} := \frac{1}{M} \langle v^{(i)}, v^{(j)} \rangle_{L^2}.    
\end{equation*}
The coefficient vectors $w^{(k)} = (w_1^{(k)}, \dots, w_M^{(k)})^\top$ are obtained by solving the finite-dimensional eigenvalue problem
\begin{equation*}
C w^{(k)} = \lambda_k w^{(k)}.    
\end{equation*}
The POD basis functions are then normalized by
\begin{equation*}
\phi_k(x) = \frac{1}{\sqrt{\lambda_k}} \sum_{m=1}^M w_m^{(k)} v^{(m)}(x), \quad k = 1, \dots, N,    
\end{equation*}
where the coefficients ensure that the set $\{\phi_k\}_{k=1}^N$ is orthonormal in $L^2(\Omega)$.

In addition, the truncation error of the expansion can be expressed in terms of the discarded eigenvalues:
\begin{equation*}
    \mathbb{E}[\|\mathcal{V} - \mathcal{V}_N\|^2_{L^2}] \approx \sum_{k=N+1}^{M} \lambda_k,    
\end{equation*}
where $\lambda_k$ are the eigenvalues of the correlation matrix. We further denote the ratio of the captured energy $E_N$ to the total energy $E_M$ by
\begin{equation}\label{rario}
    \rho=\frac{E_N}{E_M}=\frac{\mathbb{E}[\| \mathcal{V}_N\|^2_{L^2}]}{\mathbb{E}[\|\mathcal{V} \|^2_{L^2}]}=\frac{\sum^{N}_{k=1}\lambda_k}{\sum^{M}_{k=1}\lambda_k}.
\end{equation}
Then it is maximized (correspondingly \eqref{POD_continuous} is minimized) when $\lambda$ is sorted in descending order. Therefore, we conclude that the POD basis is the most efficient basis for dimensionality reduction under the $L^2$ metric for approximating data. We leave the construction of the basis in discrete in Appendix~\ref{app:pod} for completeness. 

Similar to FNO, PODNO also incorporates the kernel integral within the GSO layers (referred to here as POD layers):
    \begin{equation}\label{kernel_integration}
        \begin{aligned}
            \mathcal{K}_l[v_{(l)}](x)=\varPi_N^{-1}\left( \mathcal{R}_l\circ\varPi_N [v_{(l)}]\right)(x),
        \end{aligned}
    \end{equation}
where $v_{(l)}(x)\in\R^{d_v}$ represents the input of the $l$-th layer. Here $\varPi_N$ represents the orthonormal transform induced by the POD basis $\Phi_N:=\{\phi_k\}_{k=1}^N$, in which $\mathcal{R}_l$ can be represented by a tensor $R\in\R^{d_v\times d_v\times N\times N}$. 
Denote by $\hat v_{(l)}:=\varPi_N v_{(l)}$, and define the transformed vector as $\tilde v_{(l)}:=\mathcal{R}\hat v_{(l)}$.
That is to say, the kernel acts on $v$ becomes
\begin{equation*}
(\mathcal{K}_l[v_{(l)}])_j(x)=\sum_{k=1}^N(\tilde v_{(l)})_{j k}\phi_k(x),\quad \forall j=1,\ldots,d_v,
\end{equation*}
where $(\tilde v_{(l)})_{j k}$ are the transformed coefficients. Thus, substituting $\tilde v_{(l)}:=\mathcal{R}\hat v_{(l)}$,
    \begin{equation*}
    (\mathcal{K}_l[v_{(l)}])_j(x)=\sum_{i=1}^{d_v}\sum_{k=1}^{N}\sum_{m=1}^{N}R_{l;ijkm}(\hat{v}_{(l)})_{ik}\phi_m,\quad \forall j=1,\ldots,d_v.
    \end{equation*}

\subsection{PODNO algorithm}\label{podno_alg}
Now we present the full PODNO algorithm in Algorithm \ref{podno_alg_parameters}.
\begin{algorithm}[!ht]
\caption{PODNO Algorithm} \label{podno_alg_parameters}
\begin{algorithmic}[1]
\State \textbf{Input:} data triples of input function data, global parameters, and output function data $(a_j, \epsilon_j, u_j)_{j=1}^{M/2}$, number of POD layers $L$, max number of POD modes $N$.
\State Generation of the snapshot matrix by concatenating the $M/2$ pairs of $(a_j, u_j)_{j=1}^{M}$: $X = [a_1, a_2, \dots, a_{M/2}, u_1, u_2, \dots, u_{M/2}]$
\State Perform SVD on the covariance matrix: $XX^T$ and construct the POD transform and its inverse using the POD basis.
\State Decide the max POD mode number $N$ to retain using the energy capture ratio criterion $\rho$ in \eqref{rario}, aiming for more than $99\%$ energy if possible.
\State Initiate the network parameters.
\State Train the model by solving the optimization problem \eqref{theta}.
\State \textbf{Output:} The prediction $(u'_i)_{i=0}^{M'}$ corresponding to new data pairs $(a'_i, \epsilon'_i)_{i=1}^{M'}$.
\end{algorithmic}
\end{algorithm}

Compared with FNO \cite{lifourier}, which relies on the discrete Fourier basis, the underlying PODNO  basis is data-driven. Despite the offline cost in computing of the basis, PODNO has several advantages due to its structure: first, it focuses on modes that preserve most energy of data, instead of the general low frequency modes. Second, as it does not depend on the discrete Fourier transform, it avoids the non-physical padding and hence can be applied to non-periodic problems with ease. It also applies to problems in domains with complex geometry. Lastly, the latent variables of PODNO can be defined as real numbers, which are more efficient in backward propagation.

However, the Fourier basis has an exact analytical expression, whereas the POD basis is approximated using the method of snapshots, which may introduce errors in the space approximation. This could potentially result in a less accurate model than FNO. Therefore, we will access this through numerical experiments in Section \ref{numerical results}. 
\subsection{Universal approximation and error bounds for GSO}\label{universality}
 Now we turn to present the universality of GSO, which provides a rigorous basis for understanding the reliability and performance of our proposed operators, PODNO. 
 
 Our approach is inspired by the constructive proof for the universal approximation theorem of the FNO, established in \cite{kovachki2021universal}. With deliberate reformation and simplification, we demonstrate that, under some regularity constraints on the function spaces, GSO defined in Definition \ref{podno chain} can be constructed to approximate continuous operators in $H^s$ (see Theorem \ref{universal}). 
 
\begin{assu}[Activation function]
    Unless explicitly stated otherwise, 
    the activation function: $\sigma:\mathbb{R}\rightarrow\mathbb{R}$ in \eqref{nonlinear layer} is assumed to be non-polynomial, (globally) Lipschitz continuous, and $\sigma\in C^\infty$.
\end{assu}
\begin{assu}[Regularity of the basis functions]
Let $\Omega \subset \mathbb{R}^d$. The basis functions used to define the orthonormal transform $\varPi_N$ and its inverse $\varPi_N^{-1}$ are assumed to belong to the Sobolev space $H^s(\Omega)$, where $s > d/2$.
\end{assu}
\begin{assu}[Existence of the constant function basis]
The constant function $1(x) \equiv 1$ is included among the basis functions used to define the orthonormal transform $\varPi$ and its inverse $\varPi^{-1}$. This is typically the case in practice. Without loss of generality, we assume that $1(x)$ is taken as the first basis function.
\end{assu}
One can always select $1(x) \equiv 1$ as one of the basis functions and construct the remaining orthonormal basis functions so as to preserve orthonormality.
\begin{thm}[Universal approximation theorem]\label{universal}
Let $\Omega\subset\mathbb{R}^d$, and let $\mathcal{G}: H^s(\Omega; \mathbb{R}^{d_a})\rightarrow H^{s'}(\Omega; \mathbb{R}^{d_u})$ be a continuous operator, where $s'\geq s>d/2$. Suppose $K \subset H^s(\Omega; \mathbb{R}^{d_a})$ is a compact subset. Then, $\forall~\varepsilon>0$, $\exists$ a GSO, $\mathcal{N}: H^s(\Omega; \mathbb{R}^{d_a}) \rightarrow H^s(\Omega; \mathbb{R}^{d_u})$, satisfying
\begin{equation*}
    \sup _{a \in K}\|\mathcal{G}[a]-\mathcal{N}[a]\|_{H^s} < \varepsilon.
\end{equation*}
\end{thm}
The proof for Theorem \ref{universal} is given in Appendix \ref{proof_universal}. 
\begin{remark}
We emphasize that, for the sake of analysis, the orthonormal basis must be chosen from the Sobolev space $H^s$ with $s > d/2$. This requirement guarantees the boundedness of the $L^\infty$ norm of the basis functions, which is essential for transforming the coefficients to the function space. Nonetheless, in the algorithmic implementation and numerical experiments, we adopt the standard approach of using the POD basis under the $L^2$ metric.
\end{remark}

\section{Numerical Results}\label{numerical results}
In this section, we present the numerical performance of the proposed PODNO. First, we validate its accuracy using the publicly available 2D Darcy problem datasets in Subsection \ref{darcy results}, demonstrating its effectiveness in handling problems dominated by low-frequency components. Next, we evaluate PODNO on the Nonlinear Schr\"odinger (NLS) and Kadomtsev–Petviashvili (KP) equations in Subsection \ref{nls results} and \ref{kp results}, respectively, showing that it outperforms FNO in both efficiency and accuracy for high-frequency-dominated problems. We further compare PODNO with the POD-accelerated Lie-Trotter time-splitting scheme in Subsection \ref{pod v.s. podno}, concluding that PODNO surpasses the conventional POD-based splitting approach. Finally, in Subsection \ref{ablation}, we conduct an ablation study by systematically modifying key model components: including the number of modes, snapshots, and spatial resolution for PODNO, as well as the temporal and spatial discretization parameters for the ground truth solver, to assess their impact on overall performance.  
\subsection{Setting and error metrics}\label{setting}
Unless otherwise specified, all numerical results of PODNO are based on the following model configuration:
\begin{table}[H]
    \centering
    \caption{Basic setting of the PODNO forward model\label{setting table}}
    \vspace{-2mm}
    {\small
    \begin{tabular}{l l}
    \toprule
    Parameters & Values \\
    \midrule
    Dimension lifting $\mathcal{P}$ & Linear layer \\
    Dimension reduction $\mathcal{Q}$ & Linear Layer \\
    Width (lifting size) $d_v$& $32$ \\
    Learning rate $\epsilon_1$& $10^{-3}$ \\
    weight decay (regularization) $\epsilon_2$& $10^{-4}$ \\
    Batch size $n_\text{batch}$ & $20$ \\
    Activation function $\sigma$ & Gelu \\
    Number of the kernel integration layers $L$ &$4$\\
    Training set size $n_\text{train}$&$900$\\
    Test set size $n_\text{test}$&$100$\\
    \bottomrule
    \end{tabular}
    }
\end{table}
We train the model by minimizing the following expectation of the relative error:
\begin{equation*}
    \mathbb{E}[C(\mathcal{N}(a_j),u_j)]\approx\frac{1}{n_\text{train}}\sum_{j=1}^{n_\text{train}}\frac{\|\mathcal{N}(a_j)-u_j\|_{L^2}}{\|u_j\|_{L^2}}.
\end{equation*}
where $\mathcal{N}(a_j)$ represents the output from the operator learning model, and $u_j$ refers to the solution computed using conventional numerical methods. We consider $u_j$ as the ground truth.
\subsection{Validation on the Darcy problem}\label{darcy results} 
 As a typical example of the elliptic equation, the Darcy problem has diverse applications, such as modeling subsurface flow pressure, deformation in linearly elastic materials, and electric potential in conductive materials. We consider the following forward Darcy problem
\begin{equation}\label{eqn:darcy}
    \begin{dcases}
        -\nabla\cdot(a(x,y)\nabla u(x,y))=f(x,y), \quad & (x,y)\in\Omega,\\
        u(x,y)=0,\quad  & (x,y)\in\partial\Omega,
    \end{dcases}
\end{equation}
where $\Omega=(0,1)^2$, $a$ denotes the inhomogeneous permeability field, $u$ denotes the corresponding water head and the forcing is constant: $f(x,y)=1$. 

We are interested in learning the operator mapping the permeability to the water head, $\mathcal{G}:\,a\mapsto u$. Similar to \cite{lifourier}, we assume 
\begin{equation}\label{darcy_initial_choice}
a\sim\mu=\psi\circ \mathcal{N}\left(0,(-\Delta+9 I)^{-2}\right), \quad\psi(z)=
\begin{dcases}
    12,\; & \text{if~} z\geq0,\\
    3, \; & \text{if~} z<0,
\end{dcases}
\end{equation}
 with zero Neumann boundary conditions.  The reference solutions $u$ are obtained by using a second-order finite difference scheme on $421\times 421$ grids down-sampled to $85\times 85$ for training and validation.

\begin{figure}[ht]
    \centering

    \includegraphics[width=0.85\textwidth]{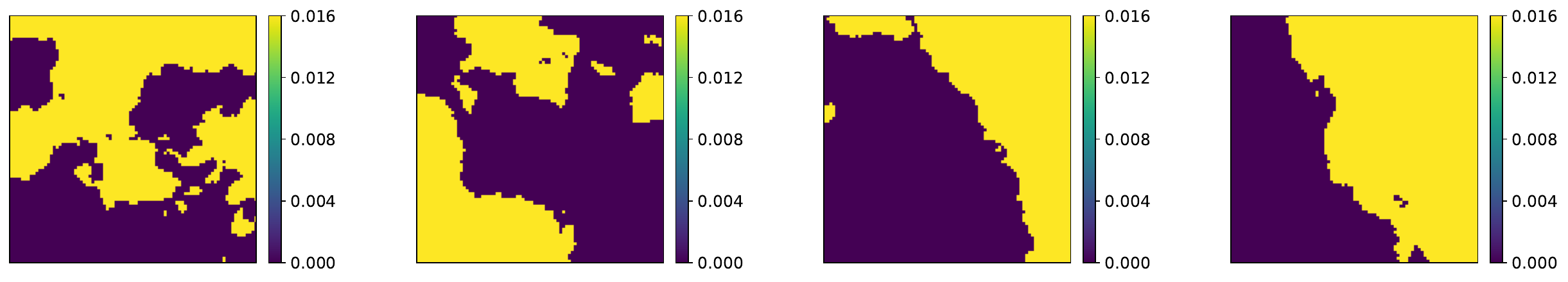}
    \vspace{1mm}

    \includegraphics[width=0.85\textwidth]{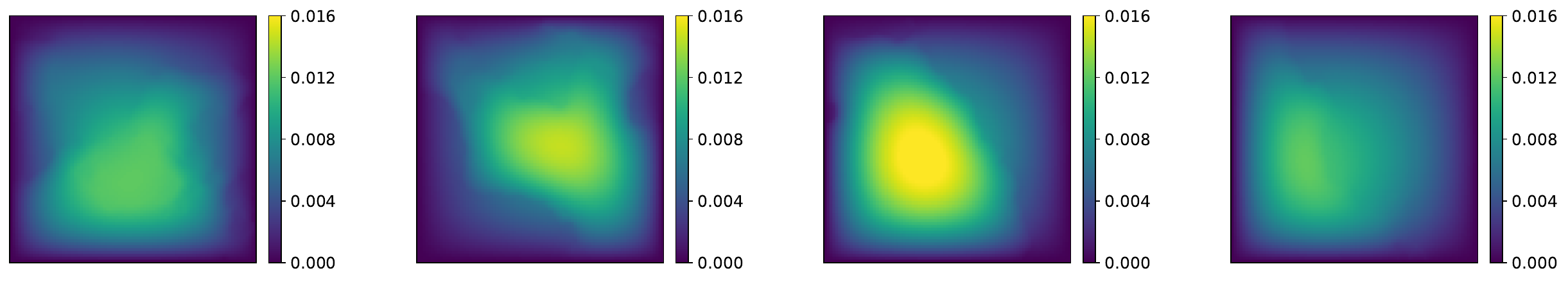}
    \vspace{1mm}

    \includegraphics[width=0.85\textwidth]{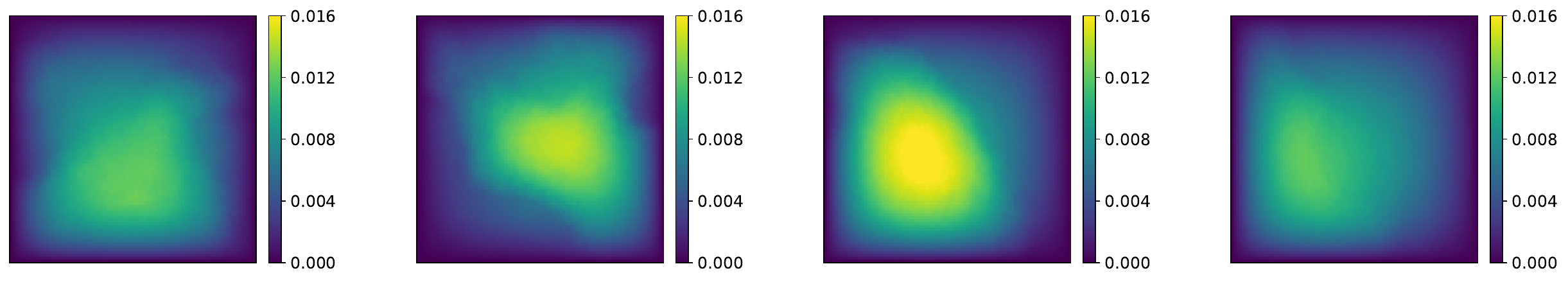}

    \caption{Visualization of the Darcy problem: (Top) Four realizations of random permeability fields $a$ generated according to \eqref{darcy_initial_choice}; (Middle) the corresponding ground truth solutions; (Bottom) the PODNO predictions using the same input fields as in the top row.}
    \label{darcy_merged}
\end{figure}

\begin{figure}[!htbp]
    \centering
    \includegraphics[width=0.4\textwidth]{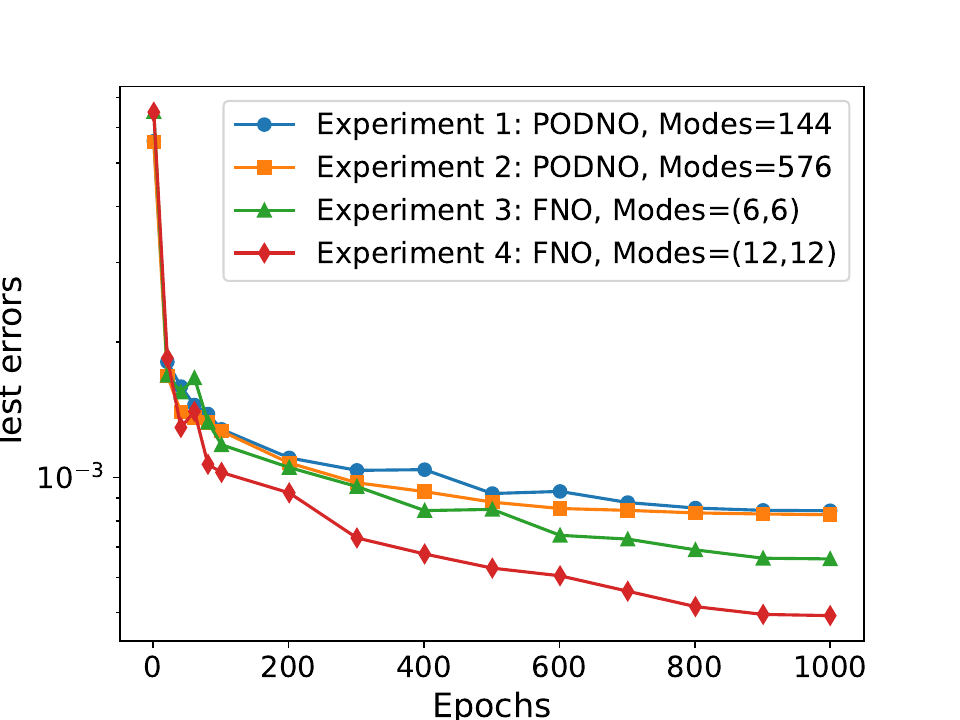}
    \caption{Test errors comparison between FNO and PODNO for the Darcy problem\label{darcy_error_fig}}
\end{figure}

With the neural network settings and parameters outlined in Table \ref{setting table}, we present in Fig.~\ref{darcy_merged} the ground truth solutions and the PODNO predictions after $1000$ training epochs, both corresponding to the permeability field. The results indicate that PODNO achieves a high-quality prediction, comparable with FNO. To further assess its accuracy and efficiency, we compare the quantified performance of PODNO and FNO in Fig.~\ref{darcy_error_fig} and Table \ref{darcy_error_tab}.
\begin{table}[!htbp]
    \centering
    \caption{Comparison of running time and test errors for the Darcy problem} \label{darcy_error_tab}
    \vspace{-2mm}
    \resizebox{\textwidth}{!}{
    \begin{tabular}{{c c c c c c}}
\hline
No.&Experiments & Parameters & Training time (sec/epoch) & Test time (sec/epoch) & Relative errors \\\hline
1&PODNO (modes=$144$) &  $606977$ &  $0.55003$ &  $0.03087$ & $8.4\times 10^{-4}$ \\
2&PODNO (modes=$576$) &  $2376449$ &  $0.71448$ &  $0.03885$ & $8.2\times 10^{-4}$ \\
3&FNO (modes=$(6,6)$) &  $606977$ &  $0.92052$ & $0.04341$ & $6.5\times 10^{-4}$ \\
4&FNO (modes=$(12,12)$) &  $2376449$ &  $1.06654$ & $0.04905$ & $4.9\times 10^{-4}$ \\
\hline
\multicolumn{5}{c}{\text{Note:\, Time for SVD in PODNO generating snapshots is  $8.25932$ seconds}} \\
\hline
\end{tabular}
}  
\end{table}

The use of complex numbers in the Fourier transform makes FNO less time-efficient than PODNO, as it requires doubling the number of real parameters to accommodate complex variables, further impacting efficiency. In addition, while the FFT that is used in FNO is theoretically efficient, it involves global operations and requires reordering of memory (e.g., bit-reversal), all of which introduce computational overhead. On GPUs, these operations are further limited by poor memory locality and non-fused execution. In contrast, PODNO uses a precomputed basis, where the transforms reduce to standard matrix multiplications, which are highly optimized on GPUs. As a result, POD transforms are significantly faster, even when the two models have the same number of parameters. For instance, as shown in Table \ref{darcy_error_tab}, Experiments 1 and 3 have the same parameter size, but PODNO requires significantly less time than FNO, taking $589.16$ seconds compared to $963.93$ seconds. Similarly, Experiment 2 completes the entire training and testing process in 761.59 seconds, while Experiment 4 takes 1115.59 seconds, despite both having $2376$K parameters. Notably, the reported times for PODNO already include the time required for SVD in generating snapshots, which is consistently accounted for in all subsequent discussions.

However, as shown in Fig.~\ref{darcy_error_fig}, while PODNO achieves error magnitudes comparable to FNO, the latter still yields slightly smaller errors for the Darcy problem. This discrepancy arises because FNO's error stems solely from training, as the projection and recovery operations in its kernel integration layers are predefined. In contrast, PODNO’s errors are influenced by both training and the generation of the POD basis. For instance, the number of snapshots, contributes to sampling errors, and thousands of data instances may be insufficient to approximate the entire space. 

Although the accuracy of the POD approximation is not directly measurable, the energy capture ratio $\rho$ in \eqref{rario} is accessible in the discretized space generated by finite snapshots and indicates truncation errors during training. We find PODNO captures $\rho = 95.76\%$ of the energy with 144 modes and $\rho = 99.16\%$ with 576 modes, yet the reduction in error shows only marginal improvement. Therefore, the observed error discrepancy may be attributed to sampling error. Notably, the test errors are significantly smaller than the theoretical truncation errors alone, likely due to the contribution of the linear compensation operation $\mathcal{W}_l$.

Darcy flow is inherently a diffusive process, characterized by the tendency to spread out and smooth sharp gradients, thereby suppressing high-frequency components in the solution. 

Since FNO truncates high-frequency data, we hypothesize that while FNO's accuracy outperforms PODNO in the Darcy problem, it may have limitations in handling high-frequency oscillations. To further investigate, we turn to dispersive systems to evaluate whether PODNO performs better in such scenarios.
\subsection{Nonlinear Schr\"odinger equations} \label{nls results}
In this subsection, we turn to the computation of the nonlinear Schr\"odinger equation (NLS) \cite{kato1987nonlinear}, a dispersive system rich in high-frequency components, where we demonstrate that PODNO outperforms FNO in capturing high-frequency patterns. Specifically, we focus on the following NLS equation with the periodic boundary conditions:
    \begin{equation}\label{eq nls}
    \begin{dcases}
        \mathrm{i}\partial_t u(x,y,t)+\Delta u(x,y,t)+V(x,y)u(x,y,t)=\frac{2}{\epsilon}(|u(x,y,t)|^{\epsilon}-1)u(x,y,t),\\    u(x,y,0)=u_0(x,y),\quad (x,y)\in\Omega=(-1,1)^2,
       \end{dcases} 
       \end{equation}
       where $\epsilon\in(-1,0)\cup(0,\frac{1}{2})$ and the initial condition $u_0$ is generated by 4 Gaussian-type wave packets,
       \begin{equation}\label{eqn:u0nls}
    u_0(x,y)=\sum_{l=1}^4e^{-\frac{\alpha}{2}[(x+\beta_{l,1})^2+(y+\beta_{l,2})^2+i\gamma_{l,1} x+i\gamma_{l,2}y]},
    \end{equation}
 with the following parameters: $\alpha=120$, $\beta_{1,1},~\beta_{1,2},~\beta_{2,2},~\beta_{3,1}\in(-0.6,-0.4)$,\\
 $\beta_{2,1},~\beta_{3,2},~\beta_{4,1}, ~\beta_{4,2}\in(0.4,0.6)$, $\gamma_{1,1}=\gamma_{1,2}=\gamma_{3,1}=\gamma_{3,2}=-2$, $\gamma_{2,1}=\gamma_{2,2}=\gamma_{4,1}=\gamma_{4,2}=2$.
   
The potential $V(x)$ is added to aggregate the high-frequency patterns,
\begin{equation*}
    V(x,y)=V_0\cos(\zeta_1x)\cos(\zeta_2y),
\end{equation*} where we take amplitude $V_0=120$ and high frequency $\zeta_1=\zeta_2=60$ for processing data in a $64\times 64$ resolution.

We aim to learn the operator mapping the initial condition $u_0$ randomly generated by \eqref{eqn:u0nls} and the coefficient $\epsilon$ to the solution at the time $T=0.5$, $\mathcal{G}: (u_0, \epsilon)\mapsto u(\cdot,0.5)$. We randomly sample 30 realizations of $\epsilon$ from the uniform distribution in $(-1,0)\cup(0,\frac{1}{2})$. The remaining neural network settings and parameters are provided in Table \ref{setting table}. Training is performed for 500 epochs, as both the training and test errors plateau beyond this point.
\begin{figure}[ht]
    \centering
    \subfigure{
        \includegraphics[width=0.85\textwidth]{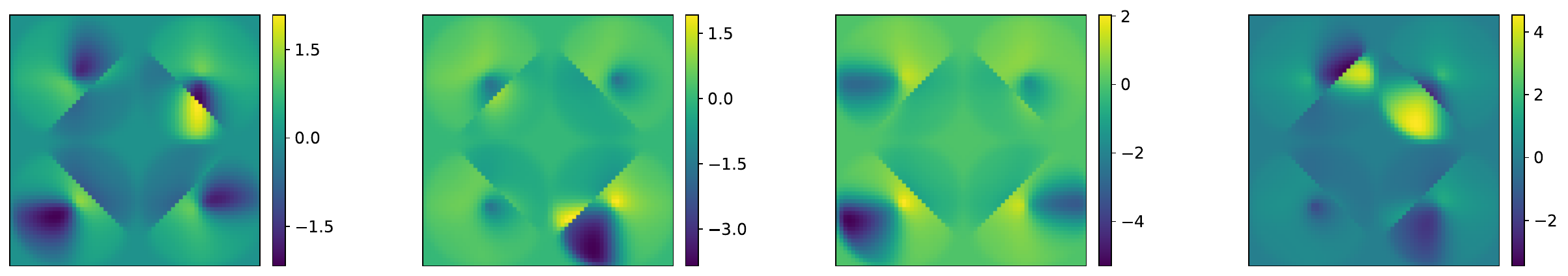}
    }
    \hspace{2mm}
    \subfigure{
        \includegraphics[width=0.85\textwidth]{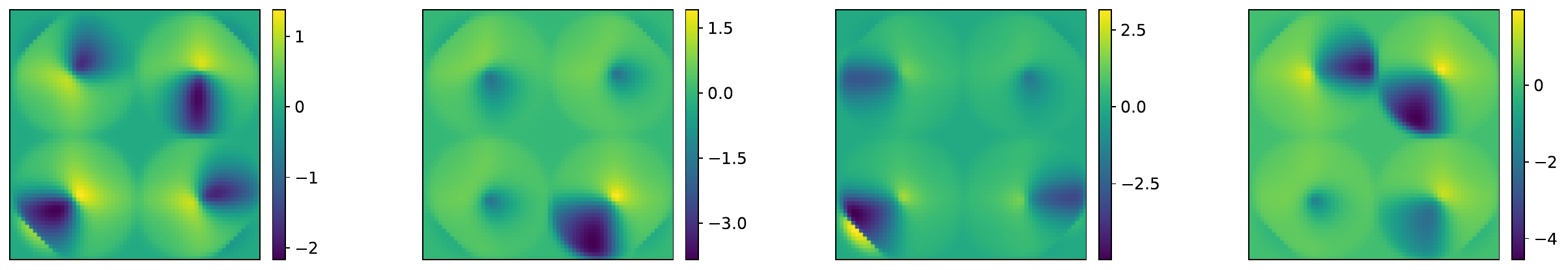}
    }
    \caption{Gaussian-type wave packet initial conditions for the NLS equation, with four randomly generated instances. The top row displays the real part, while the bottom row shows the corresponding imaginary part.}
    \label{nls_initial}
\end{figure}

\begin{figure}[ht]
    \centering
    \subfigure{
        \includegraphics[width=0.85\textwidth]{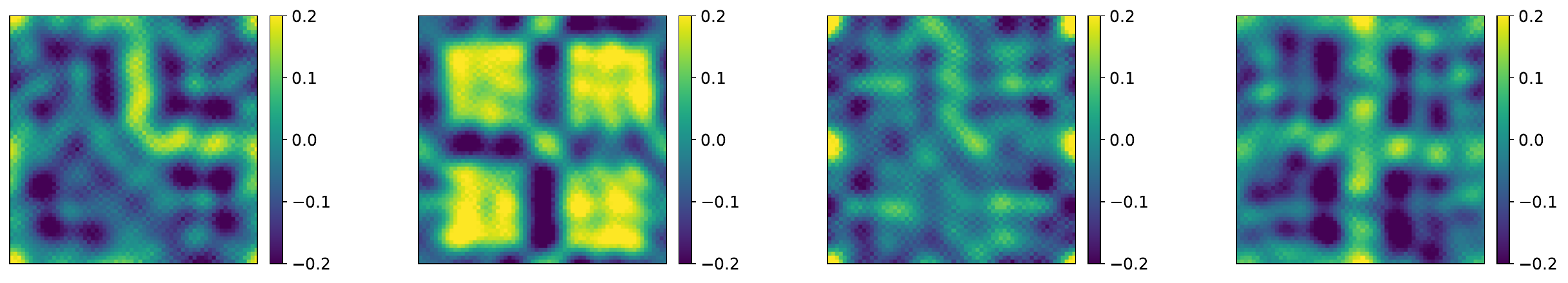}
    }
  
    \subfigure{
        \includegraphics[width=0.85\textwidth]{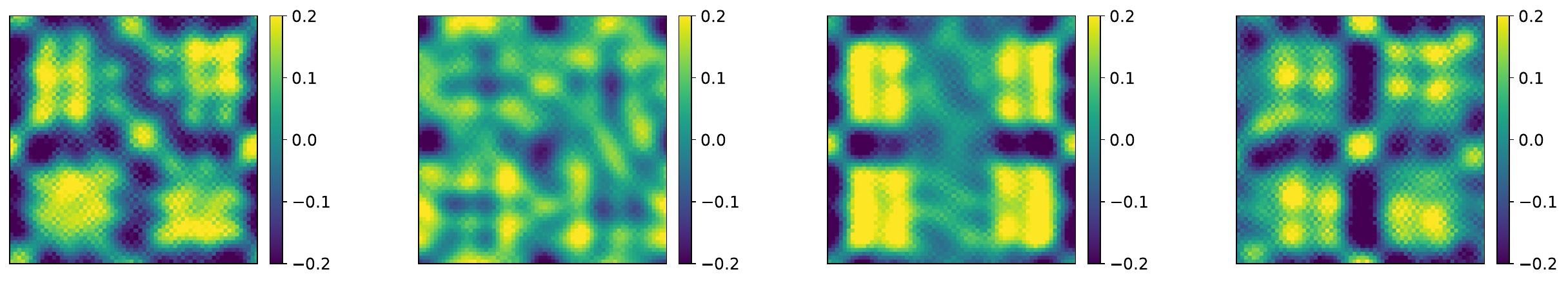}
    }
    \caption{Ground truth solutions of the NLS equation at time $T=0.5$ for the initial conditions in Fig.~\ref{nls_initial} and the coefficients $\epsilon = 0.45$, $-0.81$, $-0.02$, and $0.44$, respectively. The top row displays the real part, while the bottom row shows the corresponding imaginary part.}
    \label{nls_truth}
\end{figure}

\begin{figure}[ht]
    \centering
    \subfigure{
        \includegraphics[width=0.85\textwidth]{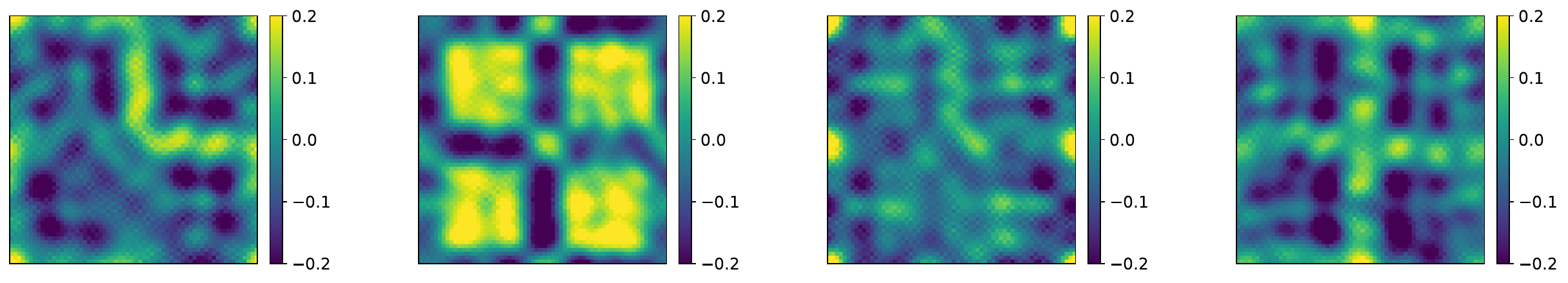}
    }
 
    \subfigure{
        \includegraphics[width=0.85\textwidth]{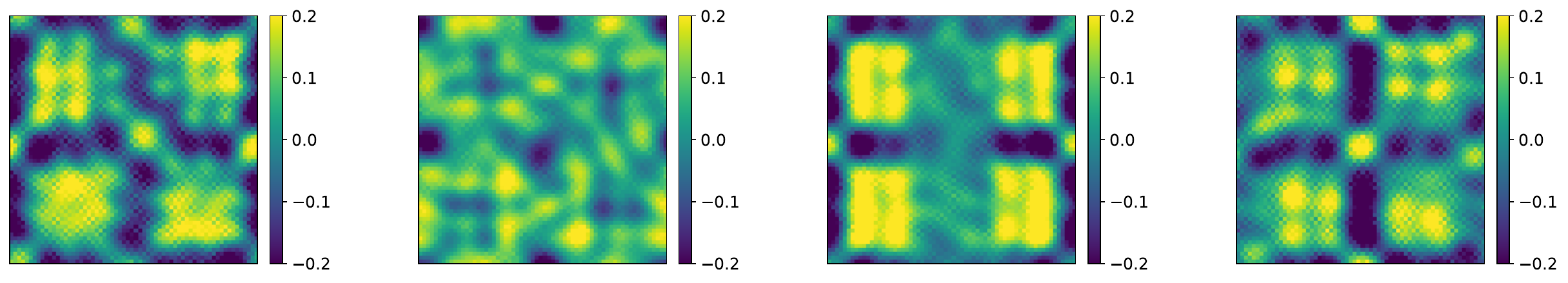}
    }
    \caption{PODNO predictions of the NLS equation at time $T=0.5$ for the initial conditions in Fig.~\ref{nls_initial} and the coefficients $\epsilon = 0.45$, $-0.81$, $-0.02$, and $0.44$, respectively. The top row displays the real part, while the bottom row shows the corresponding imaginary part.}
    \label{nls_solution}
\end{figure}

\begin{figure}[!ht]
    \centering
    \includegraphics[width=0.4\textwidth]{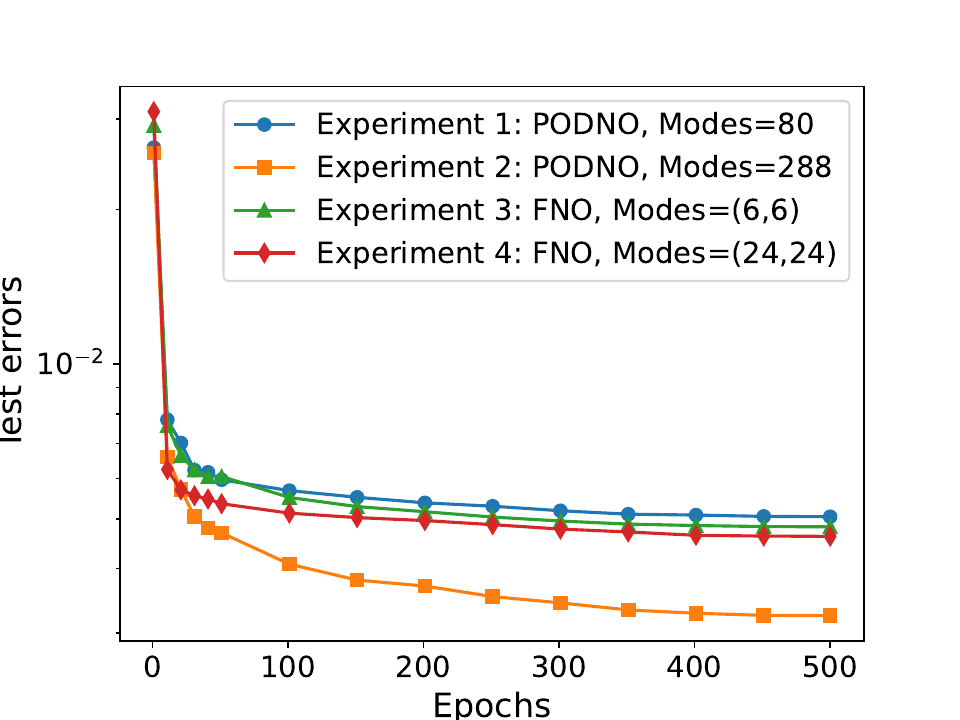}
    \caption{Comparison of numerical errors of two methods for the NLS equation} \label{nls_errors_fig}
\end{figure}

The ground truth solutions (see Fig.~\ref{nls_truth}) and the PODNO predictions (see Fig.~\ref{nls_solution}) corresponding to the initial condition (see Fig.~\ref{nls_initial}) and the coefficients $\epsilon=0.45$, $-0.81$, $-0.02$, and $0.44$ are presented below. The ground truth solutions are computed using the Lie-Trotter time-splitting scheme \cite{zhang2025low} with $1000$ time steps.
Fig.~\ref{nls_errors_fig} and Table \ref{nls_errors_tab} compare the runtime and errors of FNO and PODNO.
\begin{table}[!h]
    \centering
    \caption{Comparison of running time and errors for the NLS equation} \label{nls_errors_tab}
    \vspace{-2mm}
    \resizebox{\textwidth}{!}{
    \begin{tabular}{{c c c c c c}}
\hline
No.&Experiments & Parameters & Training time (sec/epoch) & Test time (sec/epoch) & Relative errors \\\hline
1&PODNO (modes=$80$) &  $366337$ &  $0.57543$ &  $0.02752$ & $5.0\times 10^{-3}$ \\
2&PODNO (modes=$288$) &  $1218305$ &  $0.72661$ &  $0.03932$ & $3.2\times 10^{-3}$ \\
3&FNO (modes=($6,6$)) &  $1218305$ &  $0.98969$ & $0.04292$ & $4.8\times 10^{-3}$ \\
4&FNO (modes=($24,24$)) &  $18913025$ &  $1.11293$ & $0.04003$ & $4.6\times 10^{-3}$ \\
\hline
\multicolumn{5}{c}{\text{Note:\, Time for SVD in PODNO generating snapshots is $14.23452$ seconds}} \\
\hline
\end{tabular}
}  
\end{table}

Even with a parameter count exceeding that of PODNO by more than an order of magnitude (see Experiments 2 and 4, with $1,218K$ and $18,913K$ parameters, respectively), FNO fails to match PODNO’s accuracy. Moreover, as shown in Experiments 3 and 4, increasing the number of modes in FNO results in only marginal accuracy gains, suggesting that FNO is reaching its capacity limit. In contrast, PODNO benefits from a larger number of modes, achieving better results due to an improved energy capture ratio. Specifically, the energy capture ratio increases from $\rho = 98.69\%$ with 80 modes in Experiment 1 to $\rho = 99.98\%$ with $288$ modes in Experiment 2. Furthermore, when using equivalent parameter count (see Experiments 2 and 3), the total training and testing time for FNO ($516.31$ seconds) is significantly higher than that of PODNO ($397.20$ seconds).
\begin{figure}[ht]
    \centering
    \includegraphics[width=0.22\textwidth]
    {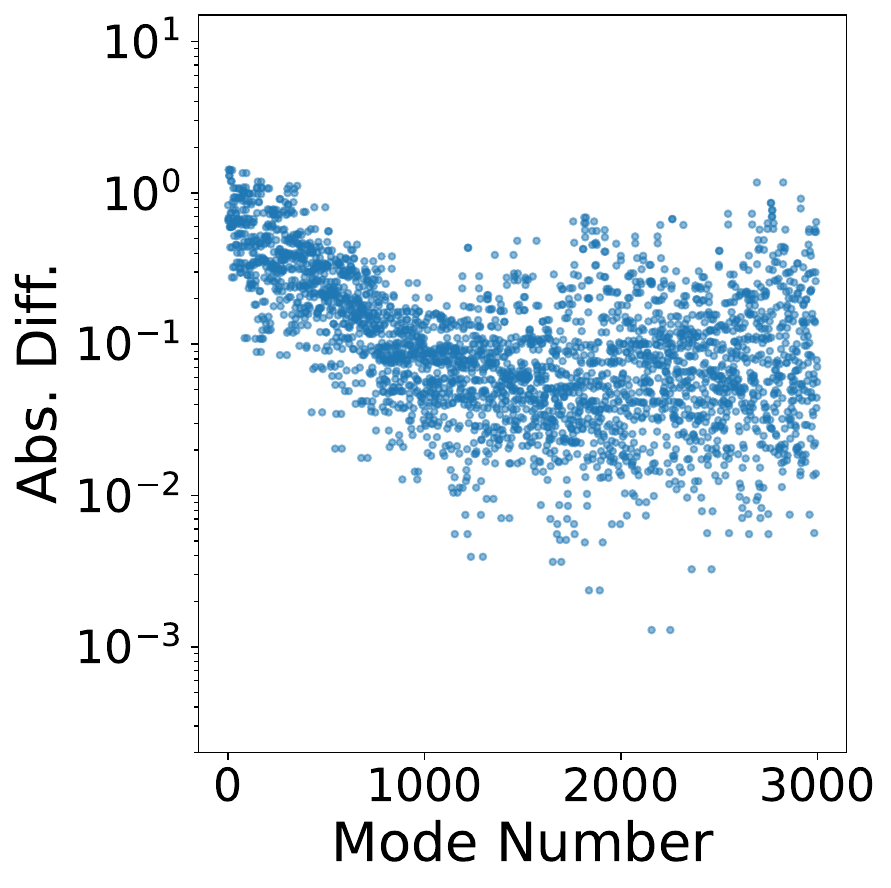}
    \includegraphics[width=0.22\textwidth]
    {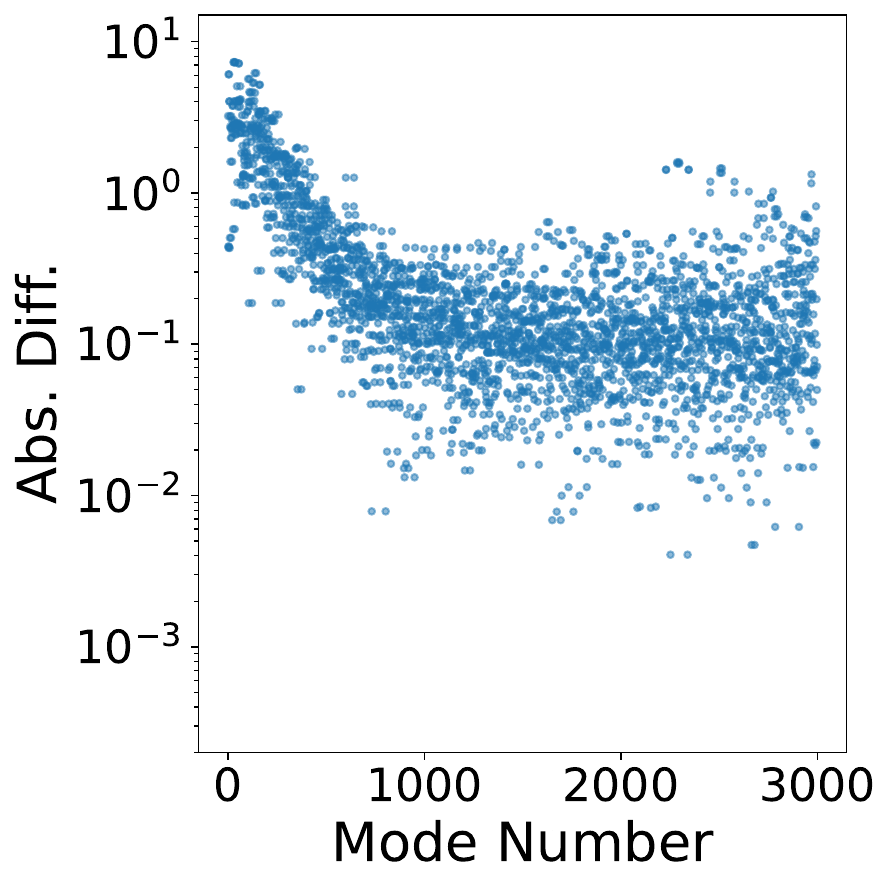}
    \includegraphics[width=0.22\textwidth]
    {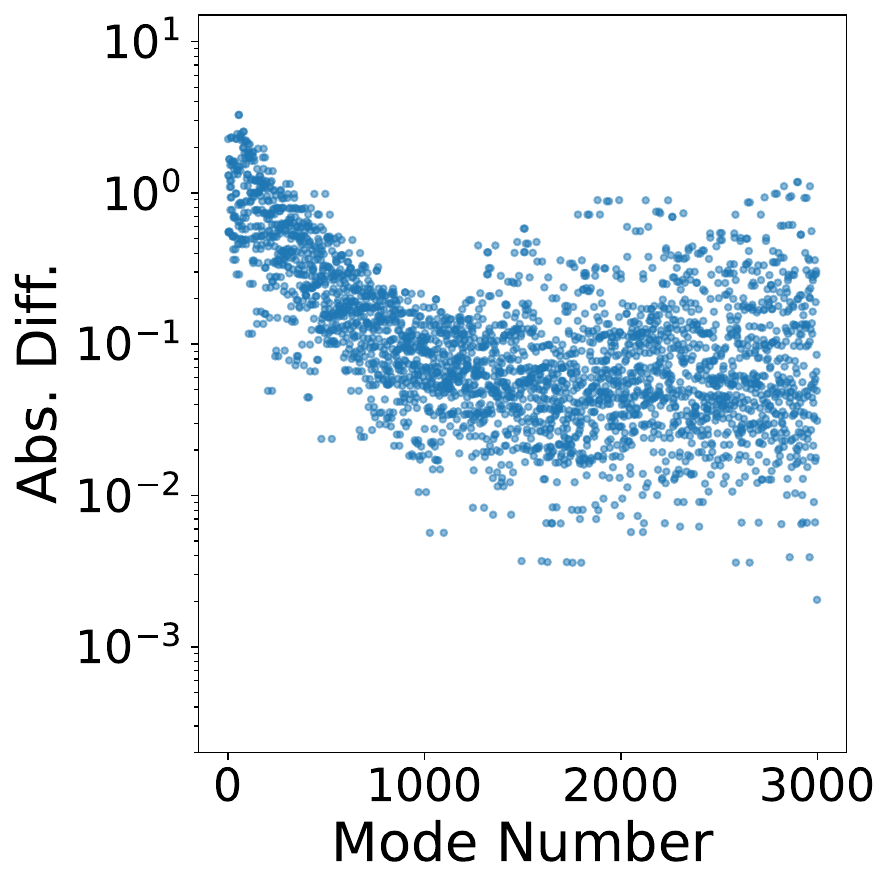}
    \includegraphics[width=0.22\textwidth]
    {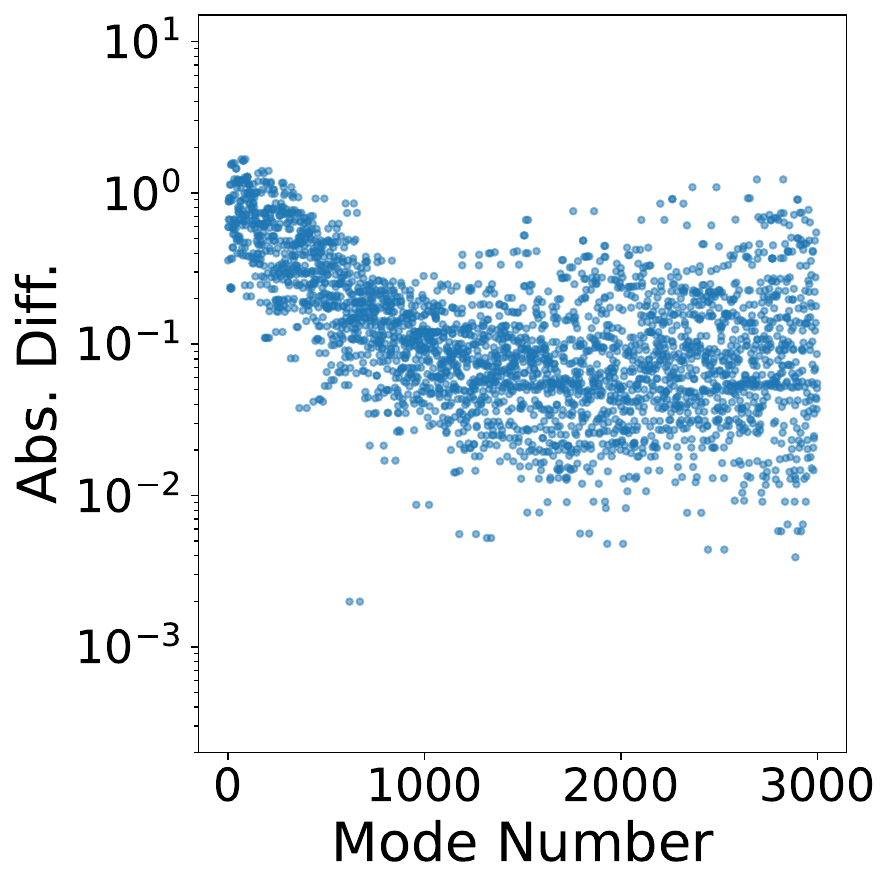}
    \includegraphics[width=0.22\textwidth]
    {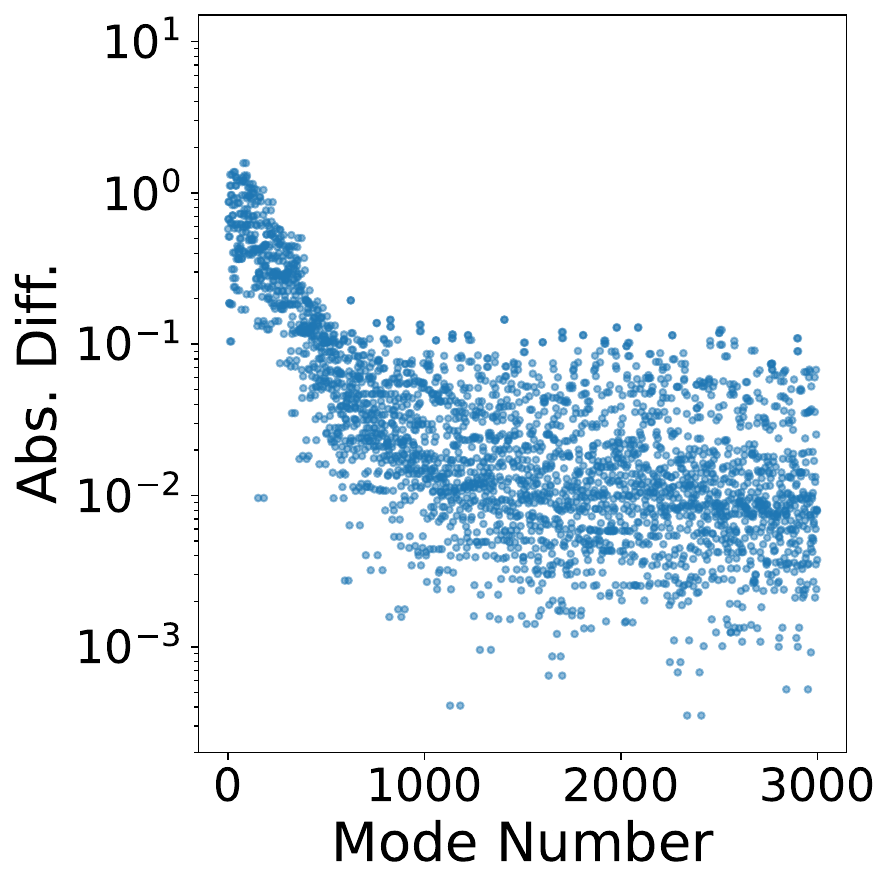}
    \includegraphics[width=0.22\textwidth]
    {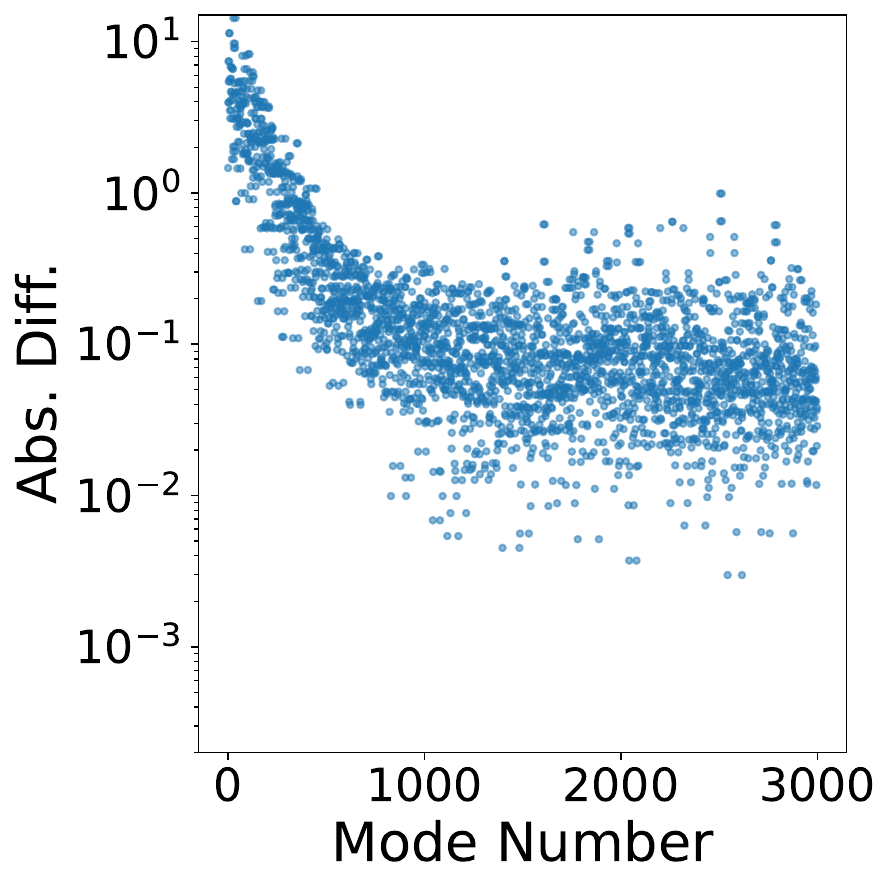}
    \includegraphics[width=0.22\textwidth]
    {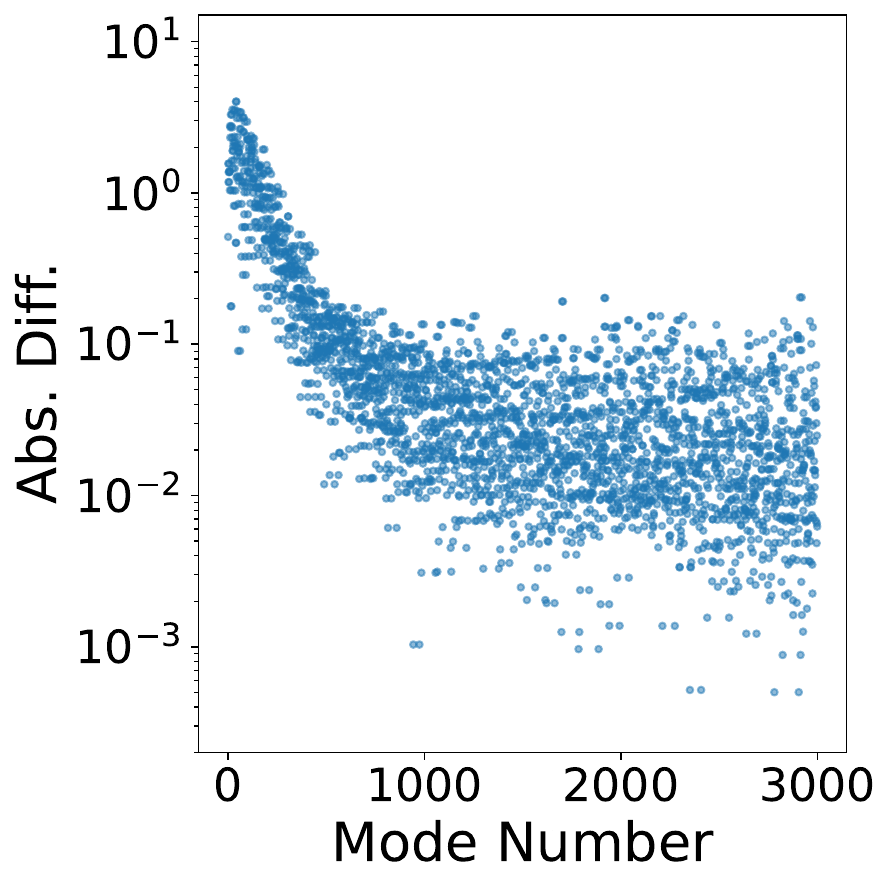}
    \includegraphics[width=0.22\textwidth]
    {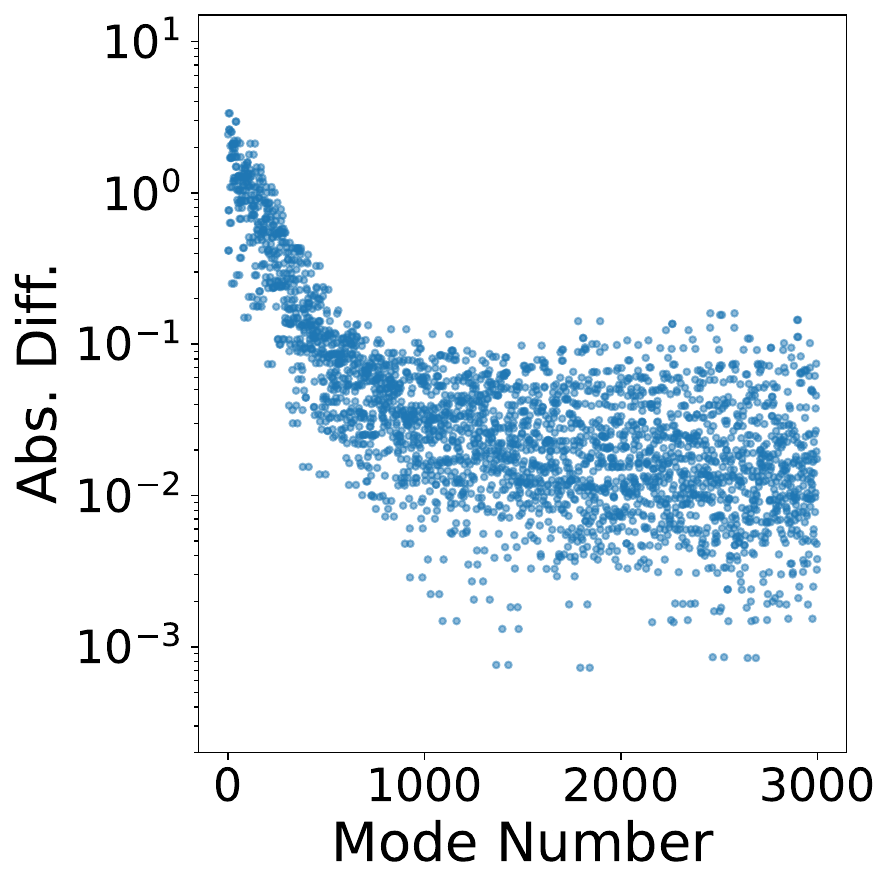}
    \caption{Scatter plot of the absolute difference of Fourier coefficients between the ground truth and the operator prediction for the NLS equation ($y$-axis) against the first $3K$ modes arranged in ascending order of frequency ($x$-axis). The instances correspond to initial conditions in Fig.~\ref{nls_initial}, with coefficients $\epsilon = 0.45$, $-0.81$, $-0.02$, and $0.44$, respectively. The top four subplots present the results for FNO, while the bottom four subplots display the results for PODNO.}
    \label{coefficients_nls}
\end{figure}
To further investigate whether PODNO effectively captures high-frequency patterns while FNO does not, we plot the Fourier coefficients of the absolute difference between the ground truth and the operator's prediction for the four examples above. These coefficients are computed using  \eqref{Fourier coefficients_nls} and then sorted in ascending order based on $k_x+k_y$. We examine the error:
\begin{equation}\label{Fourier coefficients_nls}
    |\mathscr{F}[\mathcal{N}[a_j]-u_j](k_x,k_y)|=\Big|\sum_{i_x=0}^{63}\sum_{i_y=0}^{63}[\mathcal{N}[a_j]-u_j](i_x,i_y)\cdot e^{-2\pi \textrm{i}(\frac{k_xi_x}{64}+\frac{k_yi_y}{64})}\Big|,
\end{equation}
where $k_x$, $k_y$ are the Fourier frequencies, and $i_x$, $i_y$ are the indices of the spatial domain. Since the solutions are complex, we take the modulus for visualization.
We retain the first $3$K modes and discard the higher modes to avoid the numerical inaccuracies caused by the limitations of time steps and resolution. As shown in Fig.~\ref{coefficients_nls}, the absolute error of PODNO exhibits an $L$-shaped pattern, indicating smaller coefficients for high-frequency modes error. In contrast, FNO's error coefficients form a $U$-shaped pattern, signifying an increase of the errors at higher modes. This example validates our claim that PODNO effectively captures high-frequency patterns, whereas FNO does not.

\subsection{Kadomtsev-Petviashvili (KP) equation}\label{kp results}
Another typical example of the dispersive system is the KP equation \cite{kadomtsev1970stability}, which is highly sensitive to small perturbations in the initial condition and often exhibits chaotic behavior. This example will help assess whether PODNO increasingly outperforms FNO as high-frequency components become more pronounced. We consider the KP-I model \cite{klein2007numerical}:
\begin{equation}\label{eq kp}
    \begin{dcases}
    \partial_x\left(\partial_t u+u \partial_x u+\epsilon^2 \partial_{x x x} u\right)- \partial_{y y} u=0,  & (x,y)\in\Omega,\\
    u(x,y,0)=u_0(x,y), & (x,y)\in \bar \Omega,
    \end{dcases}
\end{equation}
with the periodic boundary conditions, $\Omega=(-\pi,\pi)^2$, and $\epsilon=0.02$.
\begin{figure}[H]
    \centering

    \includegraphics[width=0.85\textwidth]{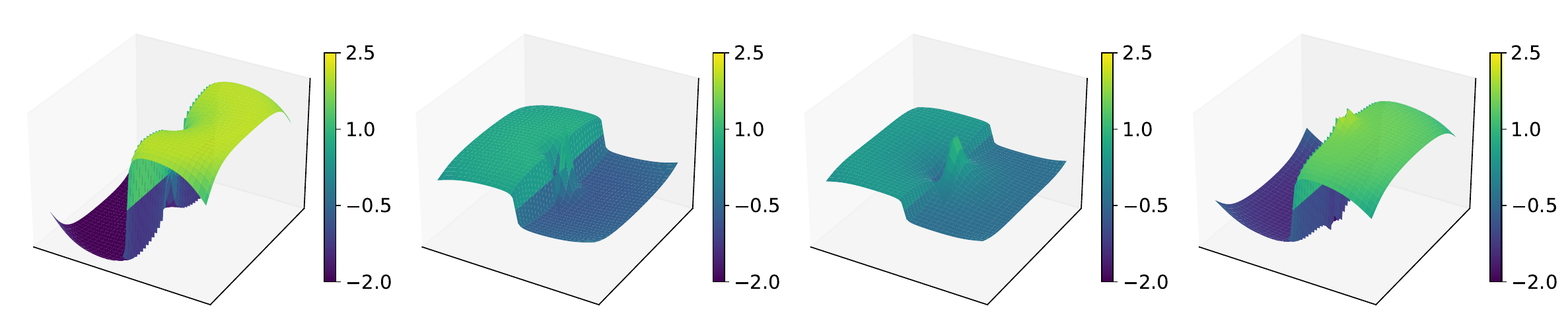}
    \vspace{1mm}
    
    \includegraphics[width=0.85\textwidth]{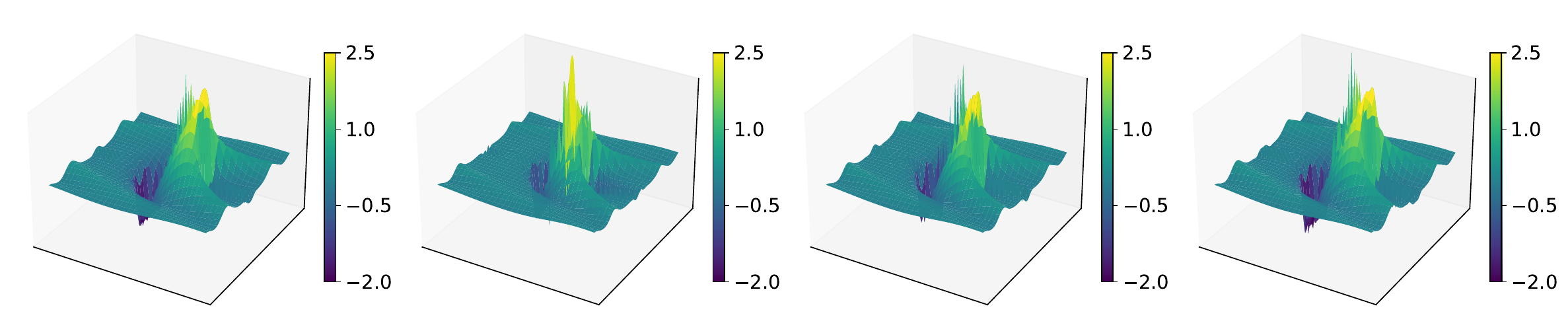}
    \vspace{1mm}
    
    \includegraphics[width=0.85\textwidth]{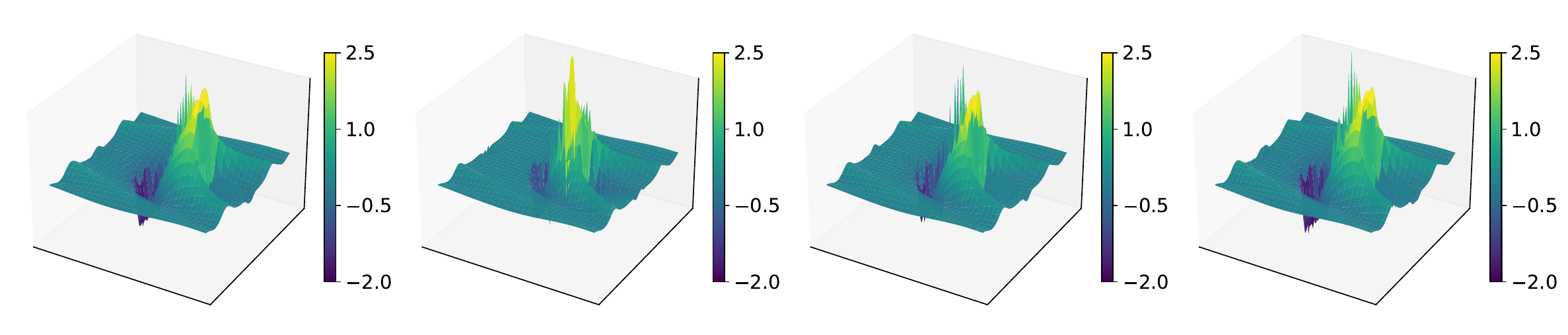}

    \caption{Visualization of the KP equation: (Top) Four randomly generated initial conditions according to \eqref{kp_initial_choice}; (Middle) corresponding ground truth solutions at time $T = 0.3$; (Bottom) PODNO predictions at time $T = 0.3$.}
    \label{kp_merged}
\end{figure}
The initial condition is chosen as
\begin{equation}\label{kp_initial_choice}
u_0(x,y)=\frac{8x\tanh{(\theta_1 r)}}{r+\delta}\sech^2{(\theta_2 r)}, \quad r=\sqrt{x^2+y^2}, 
\end{equation}
where $\theta_1,\theta_2\in(\frac{3}{2},\frac{5}{2})$ are chosen randomly from a uniform distribution. 
In the ground truth solver, we add a minimal $\delta=10^{-10}$ to the denominator of the initial condition to avoid zero. 
We aim to learn the operator mapping from the initial condition to the solution at the time $T=0.3$, $\mathcal{G}: u_0\mapsto u(\cdot, 0.3)$. The resolution is $64\times 64$ and the ground truth solver, Exponential Time Differencing Fourth-Order Runge-Kutta (ETD4RK) \cite{cox2002exponential}, generates the dataset with 1000 time steps. The other settings and parameters are given in Table \ref{setting table}. We run 1000 epochs for training the parameters of the model.
The initial conditions as well as the corresponding ground truth solutions and the PODNO predictions (see Fig.~\ref{kp_merged}) are given above, respectively. 
\begin{table}[htbp]
    \centering
    \caption{Comparison of time and errors  for the KP equation} \label{kp_errors_tab}
    \vspace{-2mm}
    \resizebox{\textwidth}{!}{
    \begin{tabular}{c c c c c c}
\toprule
No.&Experiments& Parameters & Training time (sec/epoch) & Test time (sec/epoch) & Relative errors \\\midrule
1&PODNO(modes$=9$)&	$54017$&	$0.29999$&$0.01676$&$1.01\times 10^{-3}$\\
2&PODNO(modes$=144$)&$606977$&	$0.32231$	&$0.01835$	&$3.16\times 10^{-4}$\\
3&FNO(modes$=(6,6)$)&$606977$&$0.66920$&$0.02794$&$8.93\times 10^{-4}$\\
4&FNO(modes$=(12,12)$)&$2376449$&$0.68592$	&$0.02810$	&$4.10\times 10^{-4}$\\
\midrule
\multicolumn{5}{c}{\text{Note:\, Time for SVD in PODNO generating snapshots is $3.21382$ seconds}} \\
\bottomrule
\end{tabular}
}  
\end{table}
\begin{figure}[htbp]
    \centering
    \includegraphics[width=0.4\textwidth]{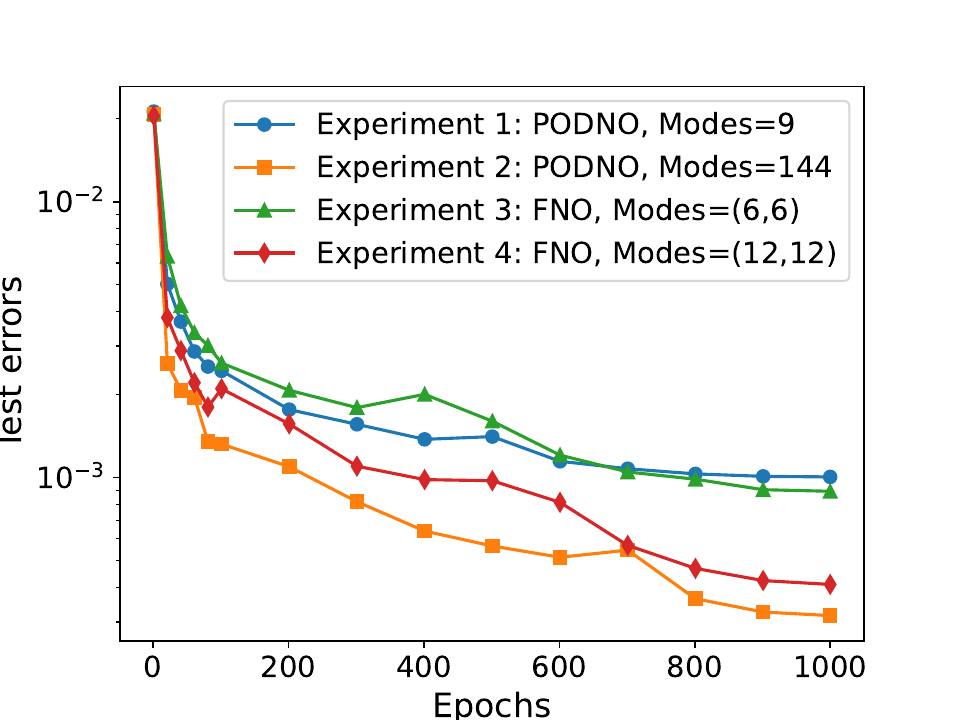}
    \caption{Test errors comparison between FNO and PODNO for the KP equation.\label{kp_errors_fig}}
\end{figure}

Due to the setting of the initial conditions, most of the differences in the solutions come from the sharp oscillations in a small domain locally while in most of the regions the solutions are smooth. This explains why, as shown in Fig.~\ref{kp_errors_fig} and Table \ref{kp_errors_tab}, PODNO with 9 modes and only $54K$ parameters (see Experiment 1, with an error of $1.01\times 10^{-3}$) could yield good results comparable to FNO which has more than 10 times the number of parameters (see Experiment 3, with $607K$ parameters and an error of $8.93\times 10^{-4}$). In addition, when PODNO and FNO have the same parameter count, PODNO produces an error that is largely smaller than FNO’s (see Experiments 2 and 3, with errors of $3.16\times 10^{-4}$ and $8.93\times 10^{-4}$, respectively). Moreover, the runtime for Experiment 2 is $343.87$ seconds, while Experiment 3 takes $697.14$ seconds. Despite having the same number of parameters, FNO takes more than twice as long to run as PODNO. This is because, in this problem, the energy is concentrated in a very small number of POD modes, leading to an extremely sparse transform matrix that significantly accelerates computation.
\begin{figure}[ht]
    \centering
    \includegraphics[width=0.22\textwidth]{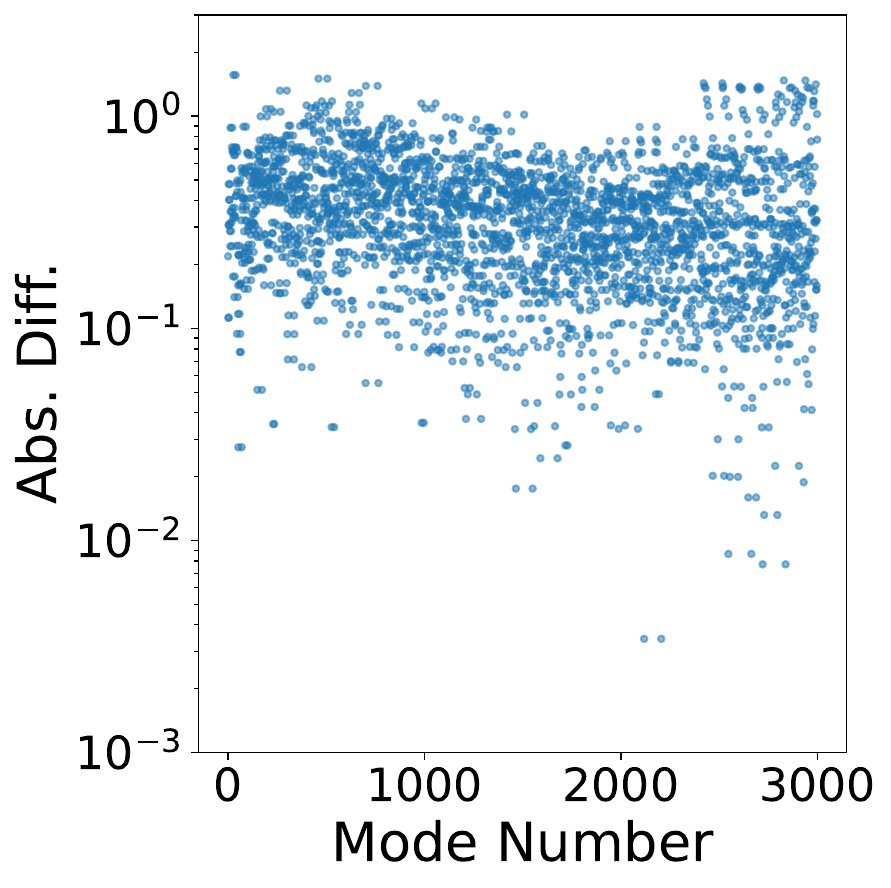}
    \includegraphics[width=0.22\textwidth]{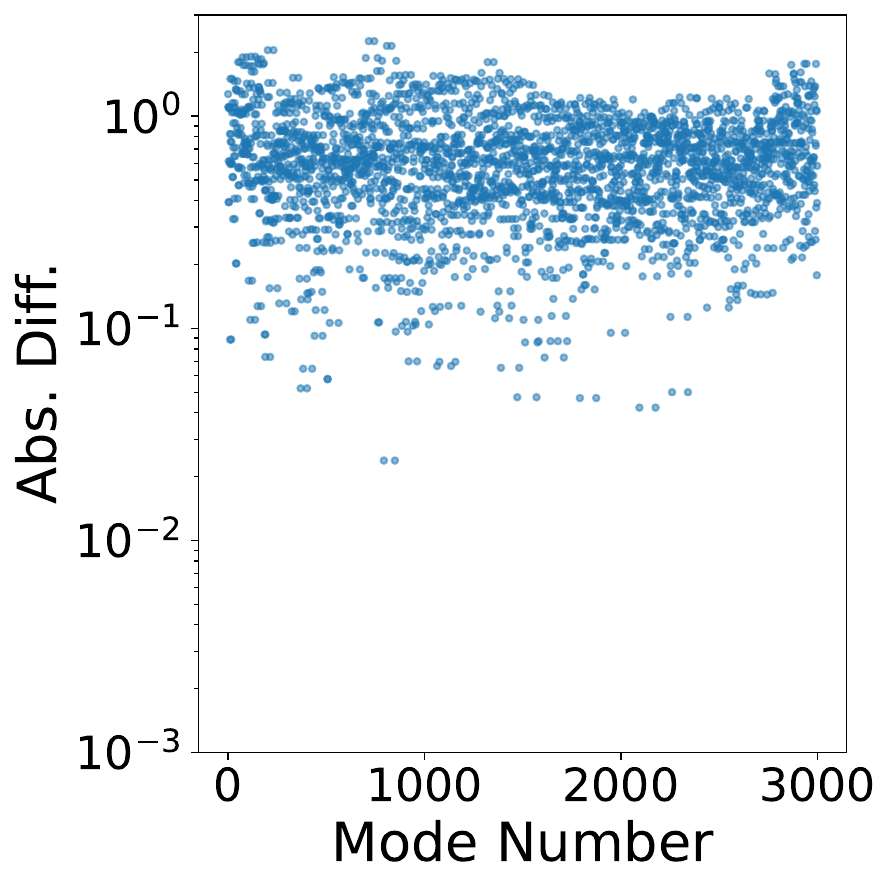}
    \includegraphics[width=0.22\textwidth]{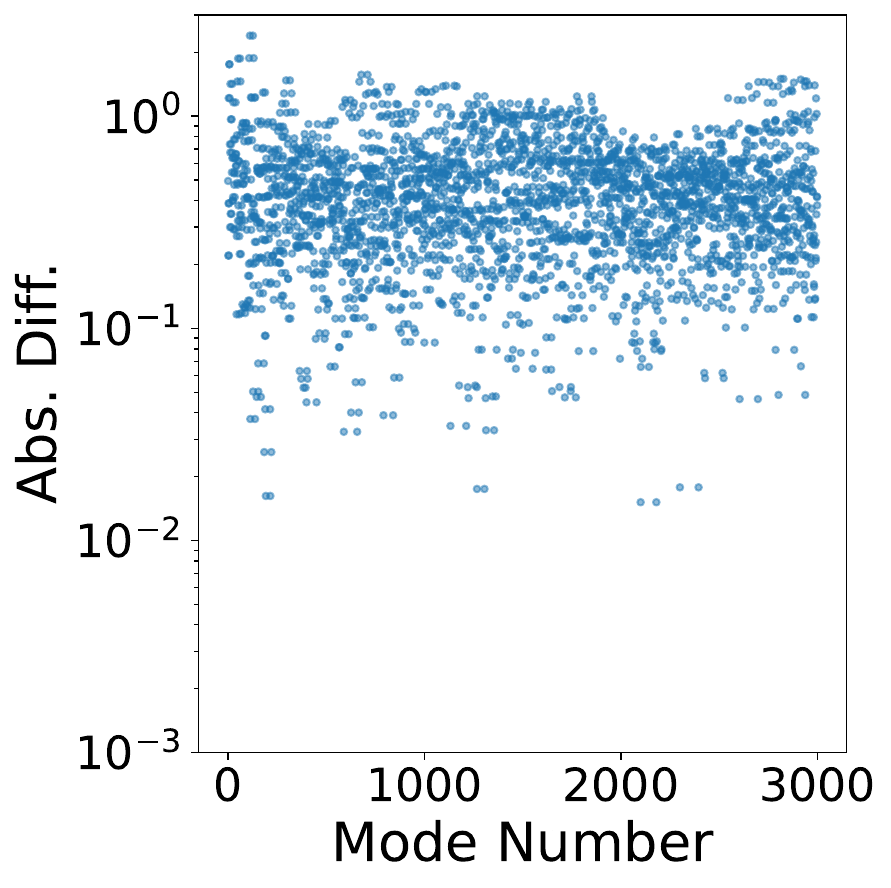}
    \includegraphics[width=0.22\textwidth]{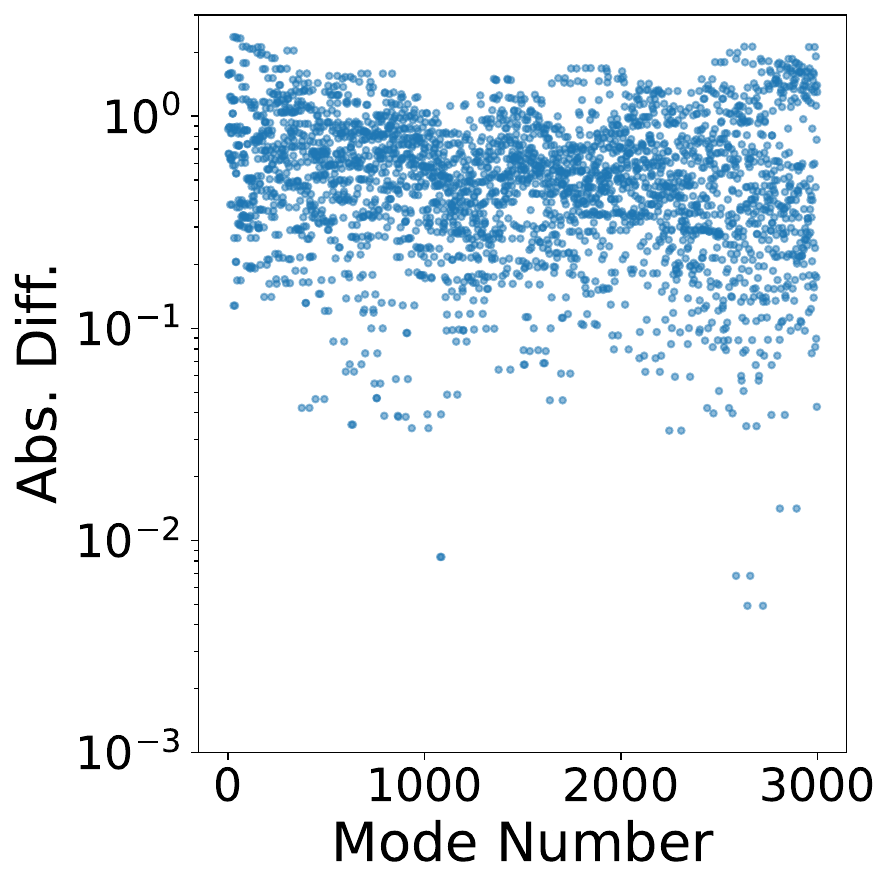}
    
    \includegraphics[width=0.22\textwidth]{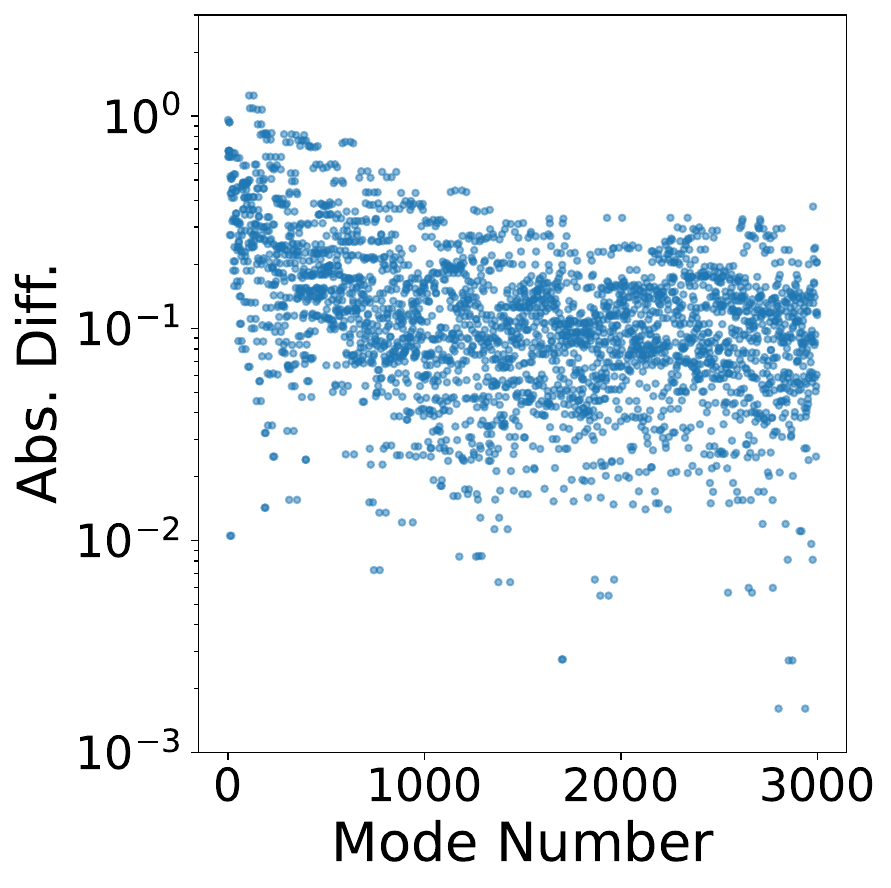}
    \includegraphics[width=0.22\textwidth]{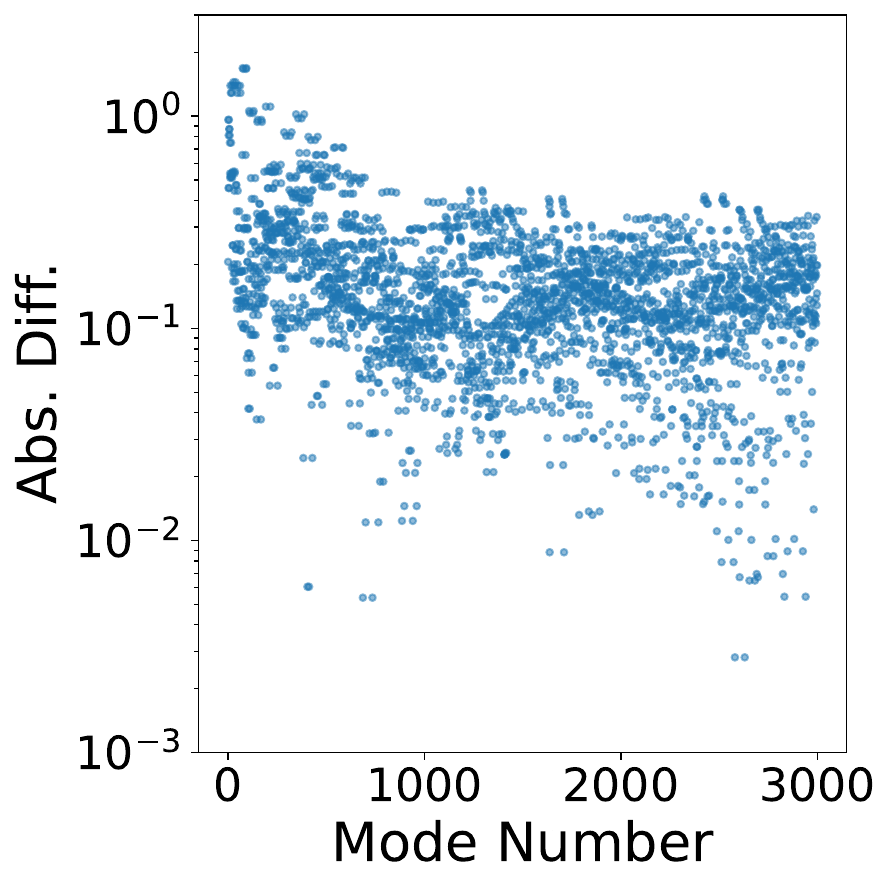}
    \includegraphics[width=0.22\textwidth]{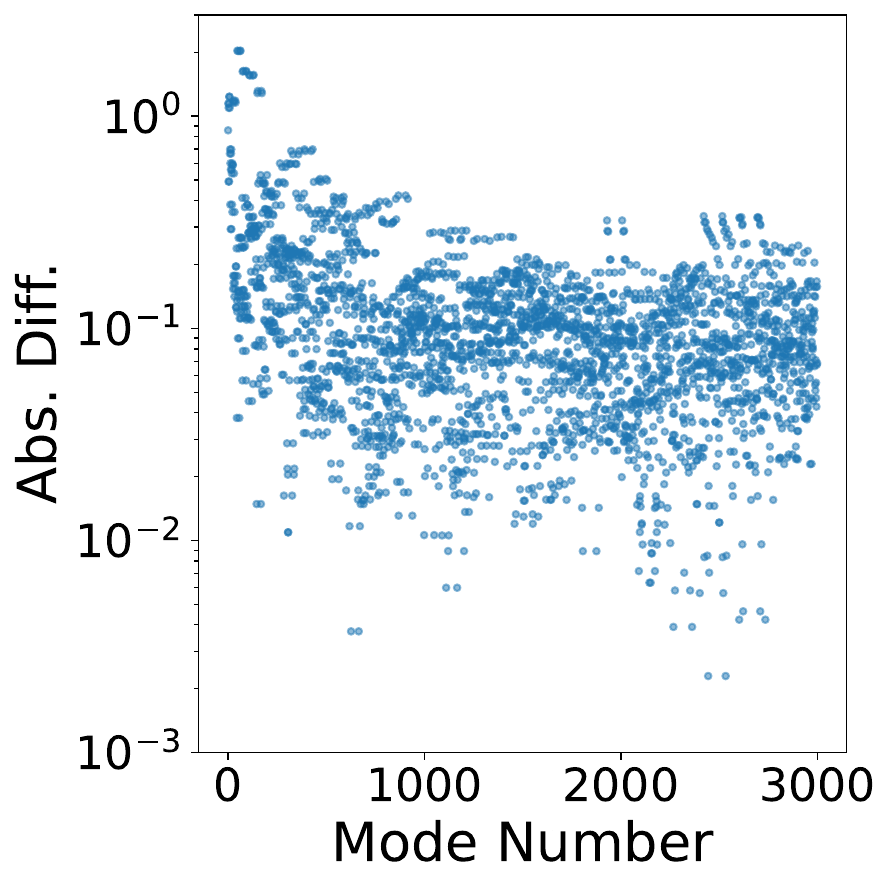}
    \includegraphics[width=0.22\textwidth]{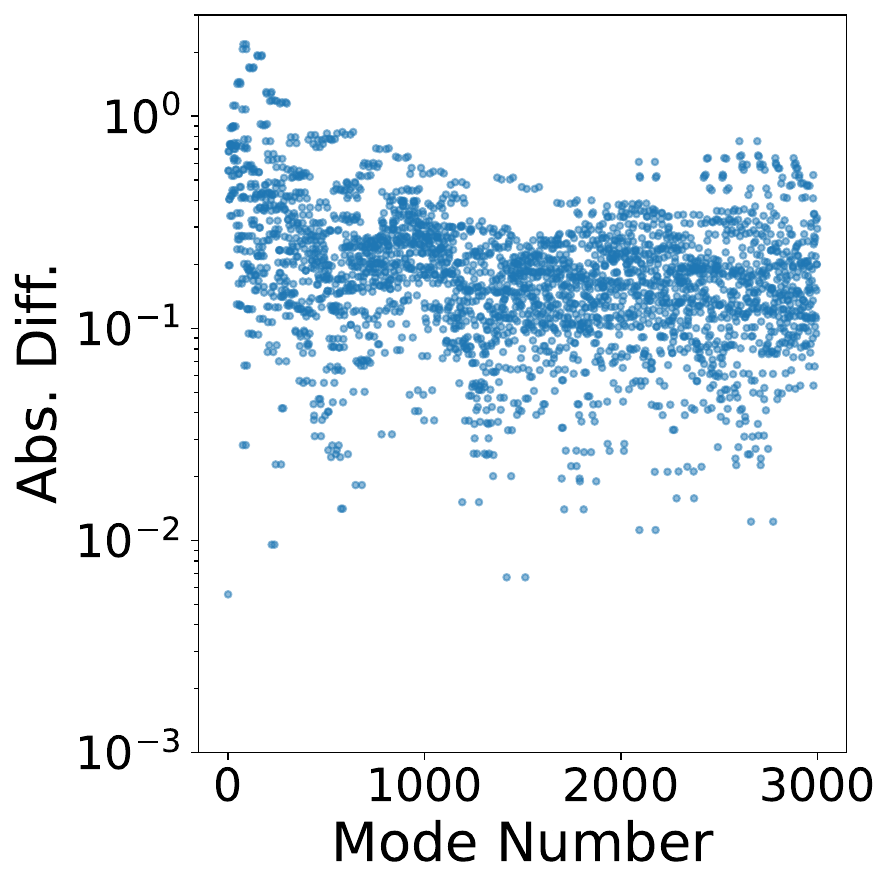}

    \caption{Scatter plot of the absolute difference of Fourier coefficients between the ground truth and the operator prediction for the KP equation ($y$-axis) against the first $3K$ modes arranged in ascending order of frequency ($x$-axis). The top row presents the results for FNO, while the bottom row shows the results for PODNO. The eight instances correspond to initial conditions in Fig.~\ref{kp_merged}, with coefficients computed using \eqref{Fourier coefficients_nls}.}
    \label{coefficients_kp}
\end{figure}

Next, we investigate whether the superiority of PODNO in solving the KP equation also stems from its ability to learn high-frequency components, and thereby achieving smaller errors. Similar to the NLS equation, we take the $64 \times 64$ resolution, which yields 4096 Fourier modes. The Fourier coefficients are also computed using \eqref{Fourier coefficients_nls} and then sorted in ascending order based on $k_x+k_y$.

Fig.~\ref{coefficients_kp} also presents the first $3K$ Fourier coefficients for the learning errors, demonstrating that PODNO effectively captures high-frequency modes, as its Fourier coefficients remain stable in the high-frequency range. In contrast, FNO exhibits pronounced error spikes at higher frequencies, aligning with our observations for the NLS equation.

In conclusion, while FNO performs well on smooth PDEs, it struggles with strong oscillations and high-frequency modes, whereas PODNO excels in these cases.

\subsection{Comparison of PODNO with time-splitting spectral methods accelerated by a POD basis}\label{pod v.s. podno}
We usually use the traditional FFT-based splitting spectral method to solve the NLS equation, the algorithm is listed in the Algorithm \ref{fft splitting} of Appendix \ref{app. spliting}. Rather than using PODNO, we consider whether it is possible to use the POD method to accelerate the splitting spectral method. We list this new approach in the Algorithm in Appendix \ref{pod-lts}.

\begin{figure}[ht]
    \centering
\includegraphics[width=0.32\textwidth]{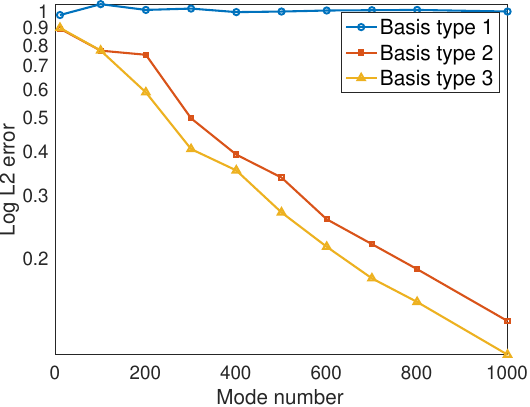} \qquad 
\includegraphics[width=0.32\textwidth]{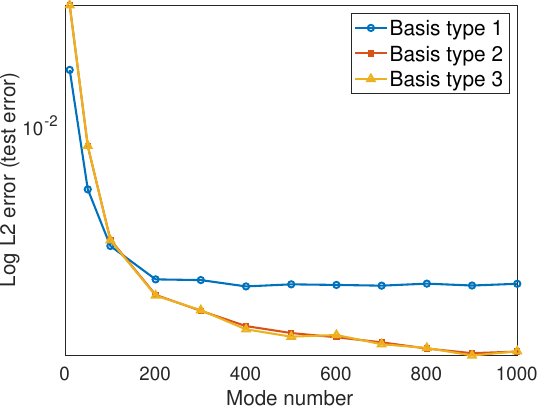}
    \caption{\small Convergence plots for POD-accelerated Lie-Trotter time-splitting solver (left) and PODNO with 900 snapshots and 500 epochs training (right). The smallest mode number is $10$ in the plot.}
    \label{comparison PODNO & POD splitting}
\end{figure}

The left side of Fig.~\ref{comparison PODNO & POD splitting} presents the convergence plot of the POD-accelerated Lie-Trotter time-splitting solver with different choices of POD mode numbers. We consider three snapshot preparation cases (see Algorithm \ref{pod-lts}): including only the solutions at the beginning and end, which proves ineffective (basis type 1); or also including the solution at each time step (basis type 2); or also including both the solution and its differences at each time step (basis type 3), which two yield better results. However, nearly all the modes are required to achieve acceptable errors in the POD-accelerated time-splitting solver. In contrast, PODNO (shown on the right side of the figure) achieves an error magnitude of $10^{-3}$ using only no more than 100 modes. In addition, generating a basis using only snapshots of the initial condition and final solution is sufficient for operator learning. This basis type 1 approach is particularly effective when using a small number of modes, as it achieves a higher energy capture ratio compared to PODNO which employs basis types 2 and 3.

\subsection{Ablation of the algorithm}\label{ablation}
In this subsection, we conduct an ablation study on both the ground truth solvers and PODNO. We examine how the resolution and time steps affect the stability of the ground truth solvers, and how the resolution, snapshots, and mode count impact the performance of PODNO.
\subsubsection{Stability of the ground truth solver} 
\begin{figure}[htbp]
    \centering

    \includegraphics[width=0.23\textwidth]{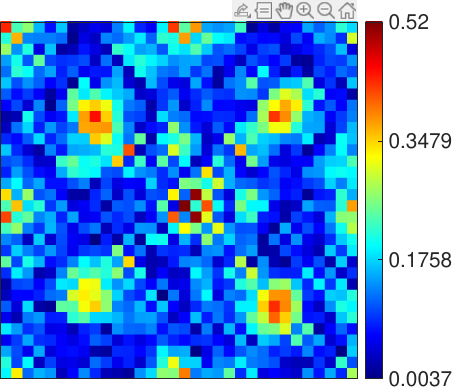}
    \includegraphics[width=0.23\textwidth]{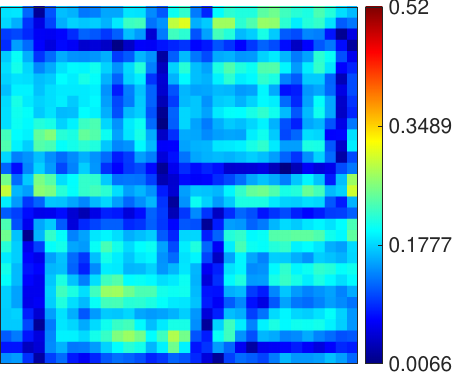}
    \includegraphics[width=0.23\textwidth]{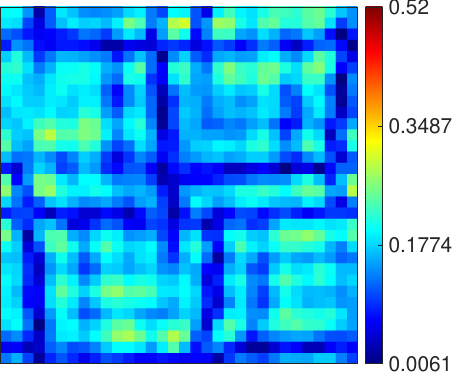}
    \includegraphics[width=0.23\textwidth]{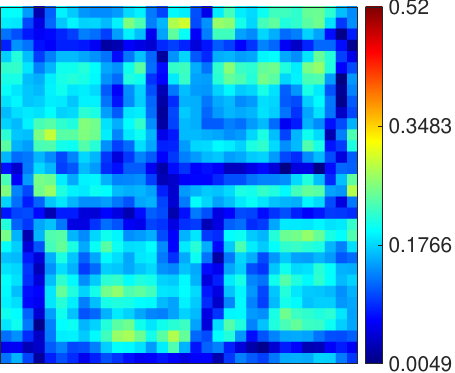}\\
\vspace{2mm}
    \includegraphics[width=0.23\textwidth]{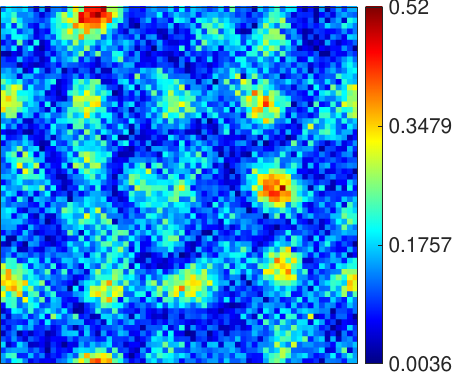}
    \includegraphics[width=0.23\textwidth]{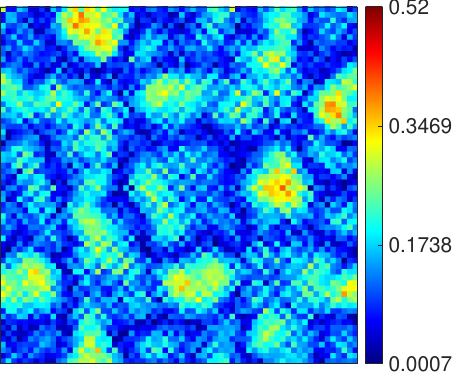}
    \includegraphics[width=0.23\textwidth]{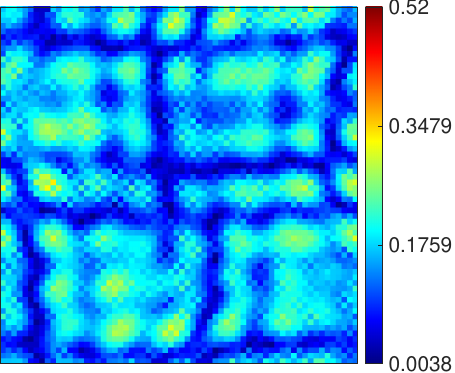}
    \includegraphics[width=0.23\textwidth]{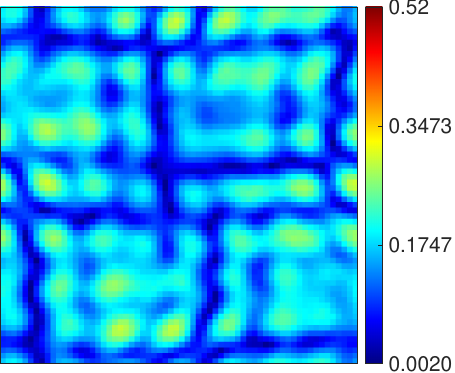}

    \caption{Solutions at time $T = 0.5$ obtained using the Lie-Trotter time-splitting scheme for the NLS equation, with a fixed initial condition and constant $\epsilon$, but varying time step and mesh configurations. The top row shows results on a $32 \times 32$ mesh, and the bottom row on a $64 \times 64$ mesh. From left to right, the columns correspond to $250$, $500$, $1000$, and $2000$ time steps, respectively.}
    \label{LTSFPS_ablation}
\end{figure}

\begin{figure}[htbp]
    \centering

    \includegraphics[width=0.23\textwidth]{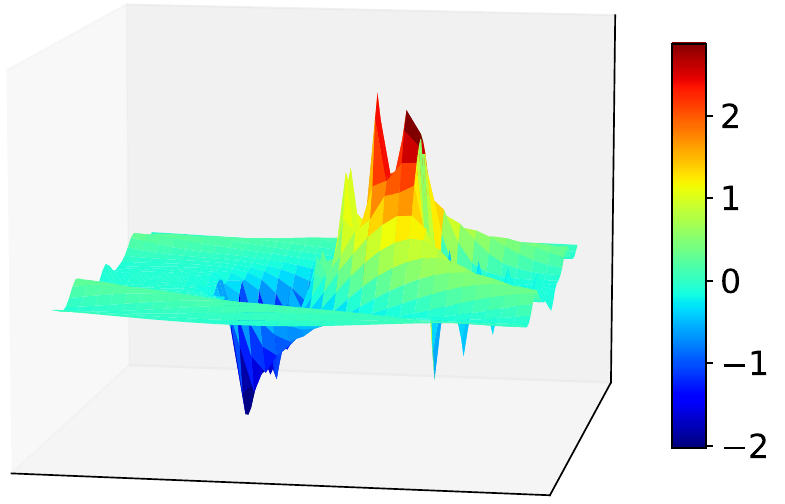}
    \includegraphics[width=0.23\textwidth]{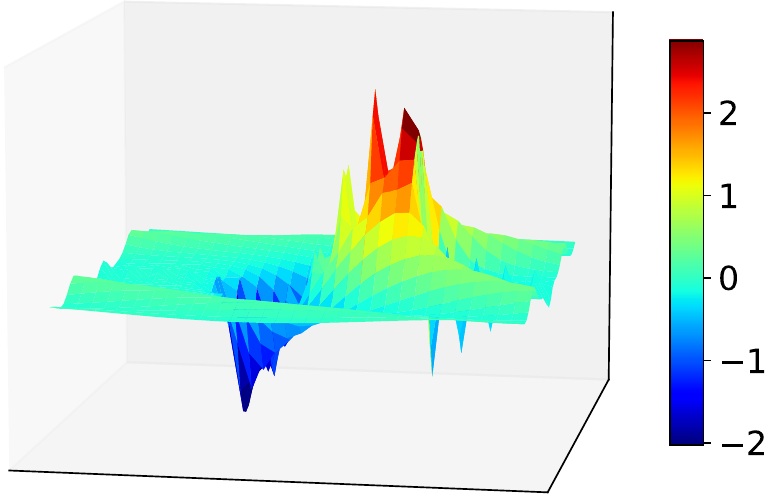}
    \includegraphics[width=0.23\textwidth]{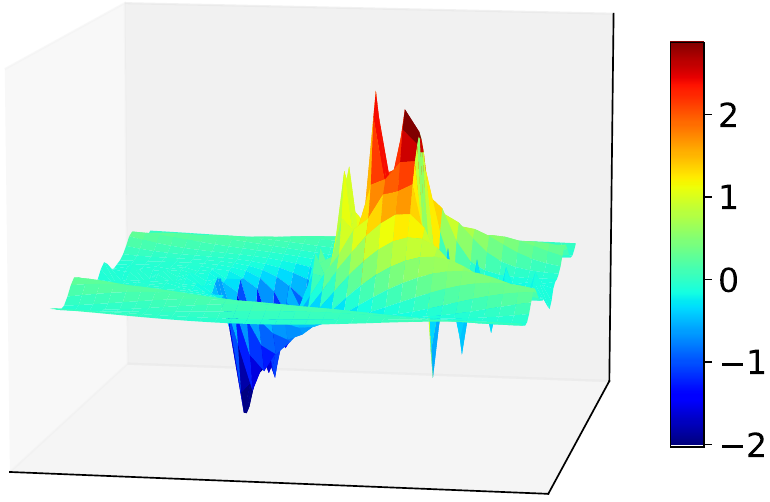}
    \includegraphics[width=0.23\textwidth]{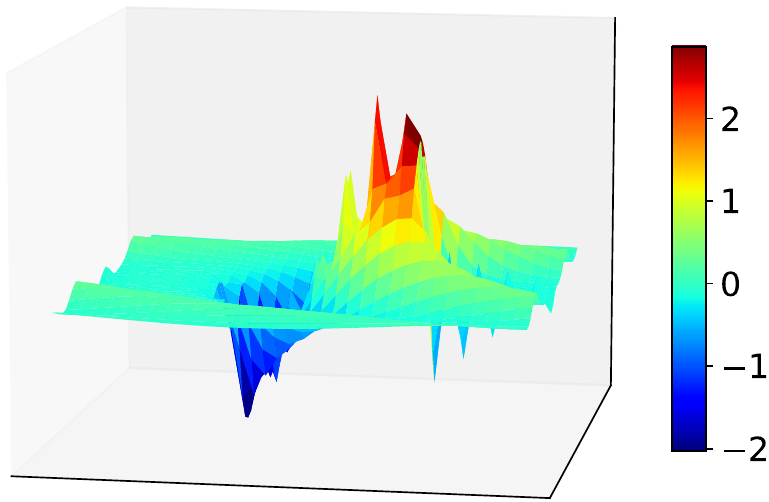}

    \includegraphics[width=0.23\textwidth]{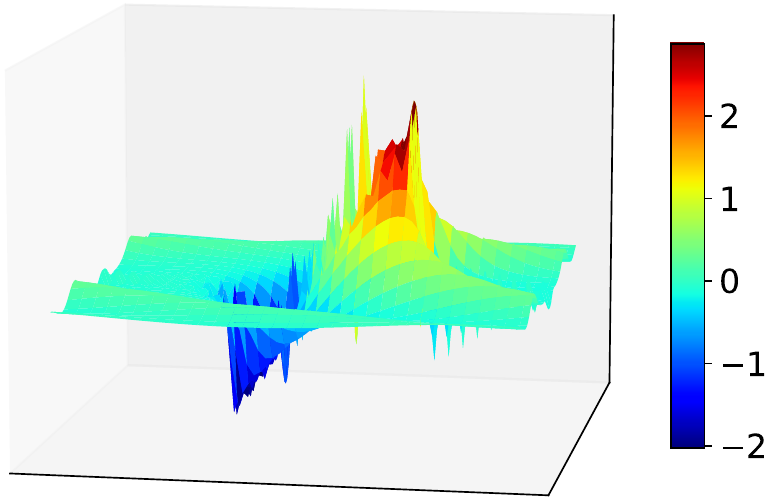}
    \includegraphics[width=0.23\textwidth]{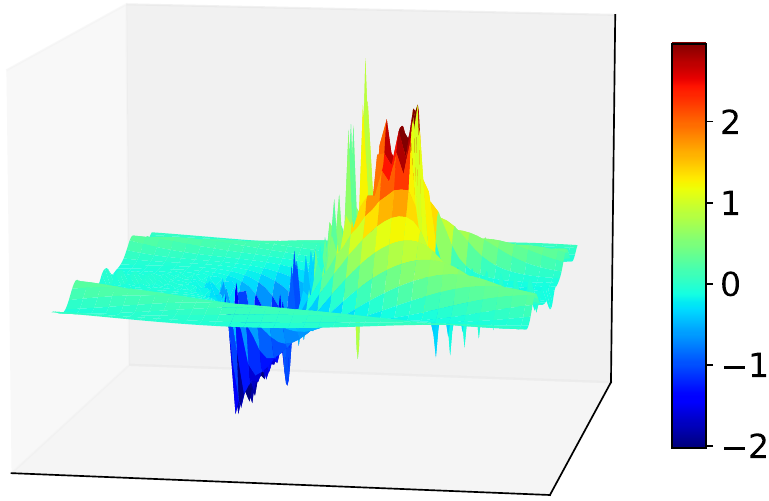}
    \includegraphics[width=0.23\textwidth]{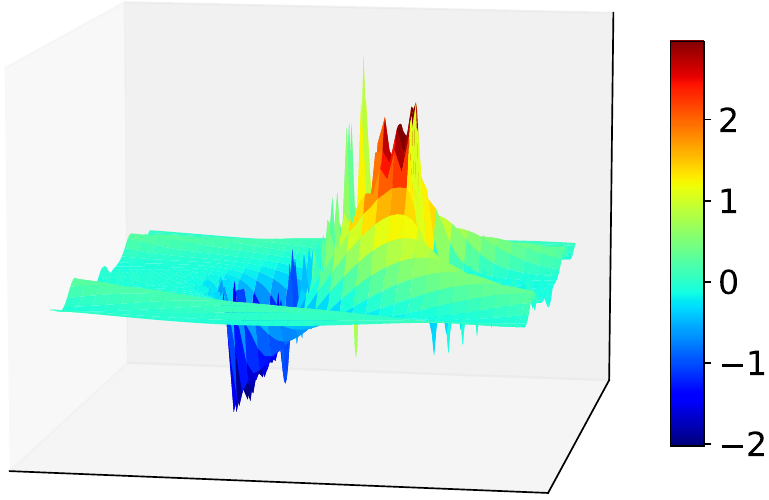}
    \includegraphics[width=0.23\textwidth]{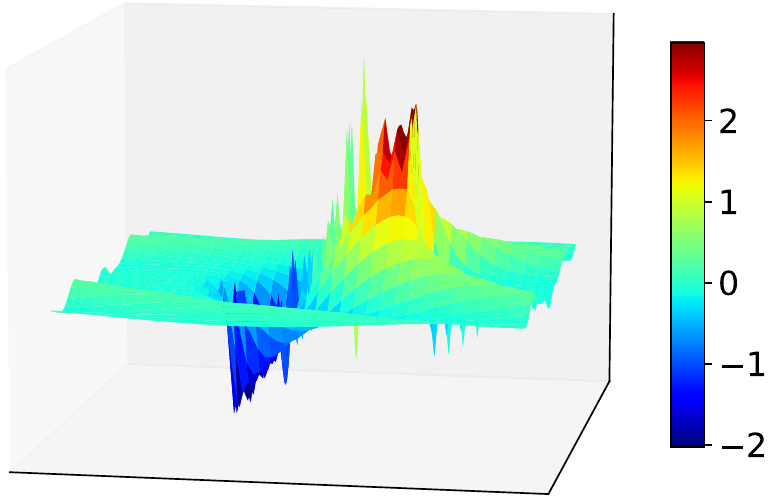}

    \caption{Solutions at time $T=0.3$ obtained using Exponential Time Differencing Fourth-Order Runge-Kutta (ETD4RK) for the KP equation. The top row shows results on a $32 \times 32$ mesh, and the bottom row on a $64 \times 64$ mesh. From left to right, the columns correspond to $10$, $100$, $500$, and $1000$ time steps, respectively.}
    \label{ETD4RK_ablation}
\end{figure}

As shown in  \cite{zhang2025low}, the $L^2$ error of the Lie-Trotter time-splitting Scheme is of the order $\mathcal{O}((\delta_t^{\frac{s}{2}}+\nmesh^{-s})\ln{\nmesh})$, where $\delta_t$ is the size of each time step, $s$ is the Sobolev regularity, $\nmesh$ is the size of resolution. This error expression highlights that the accuracy of the Lie-Trotter time-splitting methods depends on both temporal and spatial resolution. As shown in Fig.~\ref{LTSFPS_ablation}, excessively large size of the time step degrades the algorithm's performance. However, larger meshes exhibit greater robustness to variations in time step size. Due to the high-frequency components of the solution to the NLS equation, larger meshes may capture the overall behavior while dampening local oscillations. Consequently, we conclude that smaller meshes with more time steps lead to better results, which can be taken as the ground truth.

Fig.~\ref{ETD4RK_ablation} compares the solutions of the KP equation under different mesh sizes and time steps. In our setup, high-frequency variations are confined to a small local domain, while the remaining larger region exhibits uniform behavior. As a result, the algorithm only fails in fewer than tens of time steps for a $T=0.3$ time span. Similar to the case of the NLS equation, the severe oscillations in the solution necessitate a finer resolution to capture singular points accurately. Although $32 \times 32$ meshes can capture most of the energy, they may not fully resolve localized high-frequency features.
\subsubsection{Effect of spatial discretization on PODNO}
The discretized version of the POD method encounters an inherent challenge when applied to very large discretized datasets because of the use of SVD. In high-frequency-dominated PDE solutions, the energy is typically concentrated in a small number of low-frequency modes, with higher-frequency components contributing minimally. This leads to a rapid decay of singular values, where many higher-order singular values approach zero. At higher resolutions, this decay becomes even more pronounced, further exacerbating the issue of ill-conditioning. That is why we can only consider at most $64\times 64$ meshes in the prediction of the NLS and KP equations.
\begin{figure}[ht]
    \centering
\includegraphics[width=0.4\textwidth]{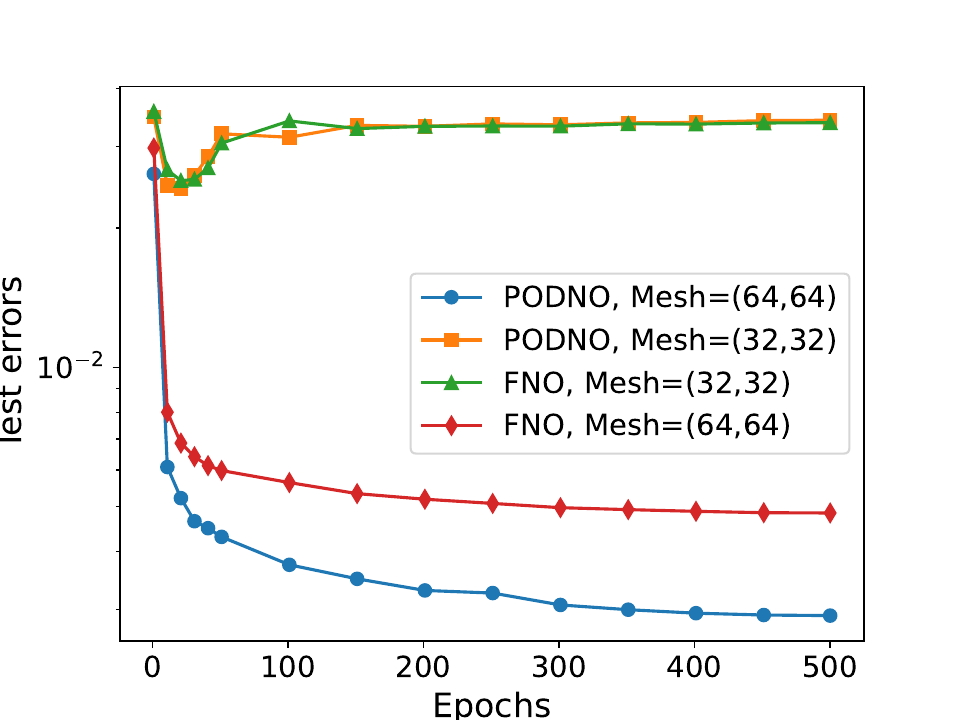}
~
\includegraphics[width=0.4\textwidth]{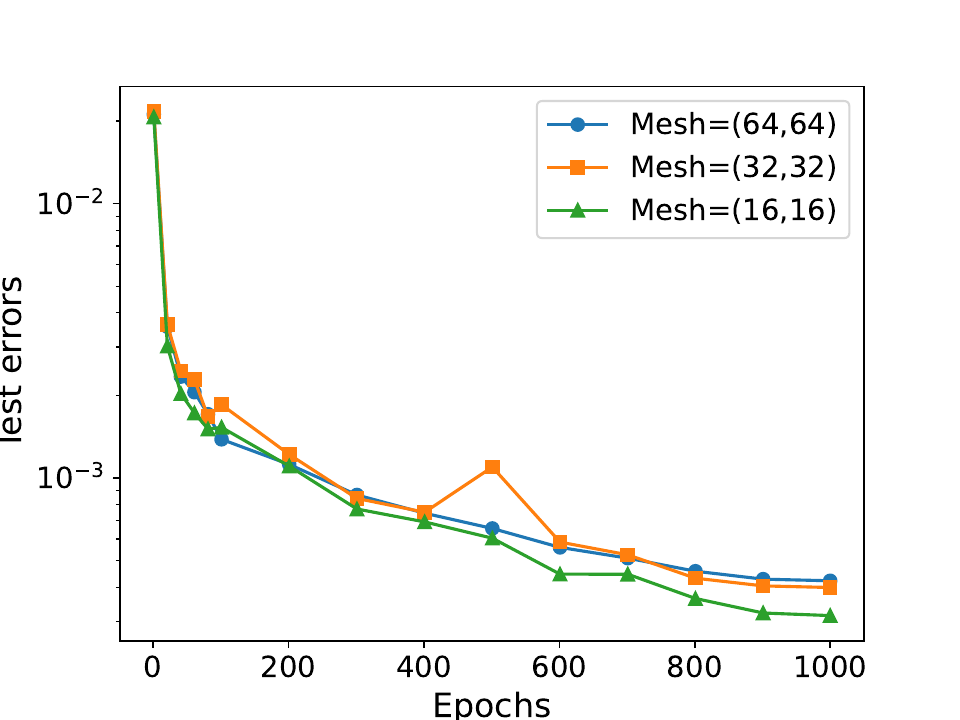}
    \caption{Error comparison between different spatial discretization for the NLS equations (left) and the KP equations (right).}
    \label{PODNO_ablation_mesh}
\end{figure}

Fig.~\ref{PODNO_ablation_mesh} shows the error comparison for different spatial discretizations (at most $64\times 64$) applied to both the NLS equation and the KP equation. For the NLS equation, PODNO fails to learn effectively with $32 \times 32$ meshes. This failure is attributed to the inability to capture high-frequency singularities, which leads to missing important dynamics in the ground truth. In contrast, for the KP equation, PODNO produces nearly the same small error across all three mesh sizes. Interestingly, the $16\times 16$ mesh even exhibits a slightly lower error, as the high-frequency component errors in small regions are smoothed out by the larger mesh. However, this does not necessarily indicate that the larger mesh is superior, as our primary focus is on high-frequency components, and the reduced error does not capture the essential dynamics of the mechanism.

We conclude that while PODNO is resolution-invariant, it still requires a high-resolution solver to effectively capture energy patterns in problems dominated by high frequencies.
\subsubsection{Effect of snapshot numbers on PODNO}
We noted that the error in PODNO arises from both the approximation of the solution space and the learning process. The number of snapshots primarily affects the approximation of the solution space, with a larger number of snapshots generally leading to a better approximation.
\begin{figure}[!ht]
    \centering
  \includegraphics[width=0.4\textwidth]{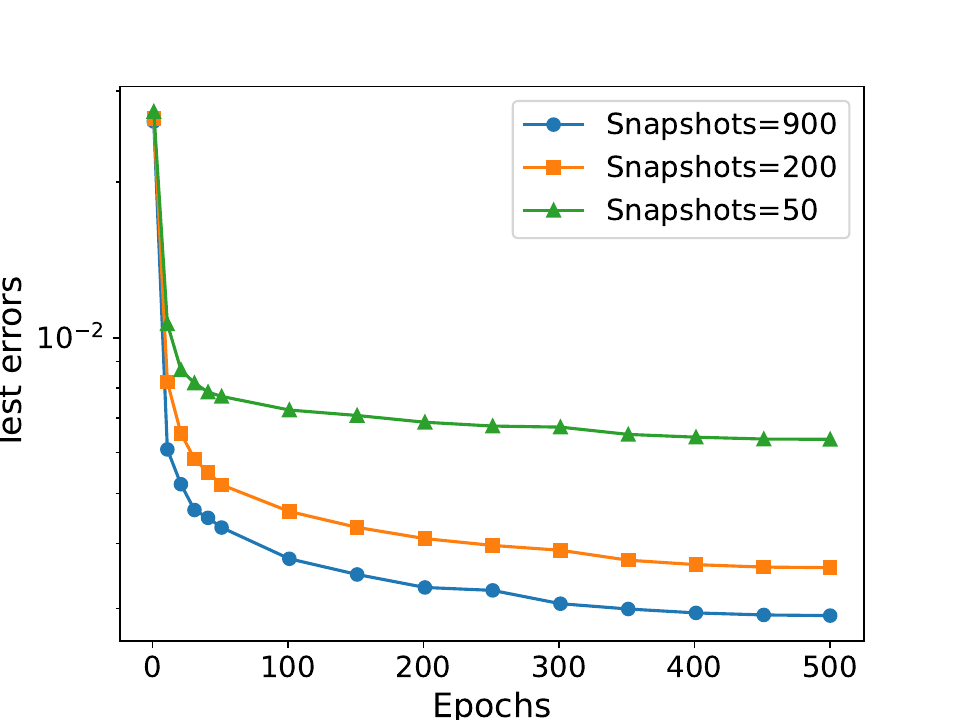} ~ 
    \includegraphics[width=0.4\textwidth]{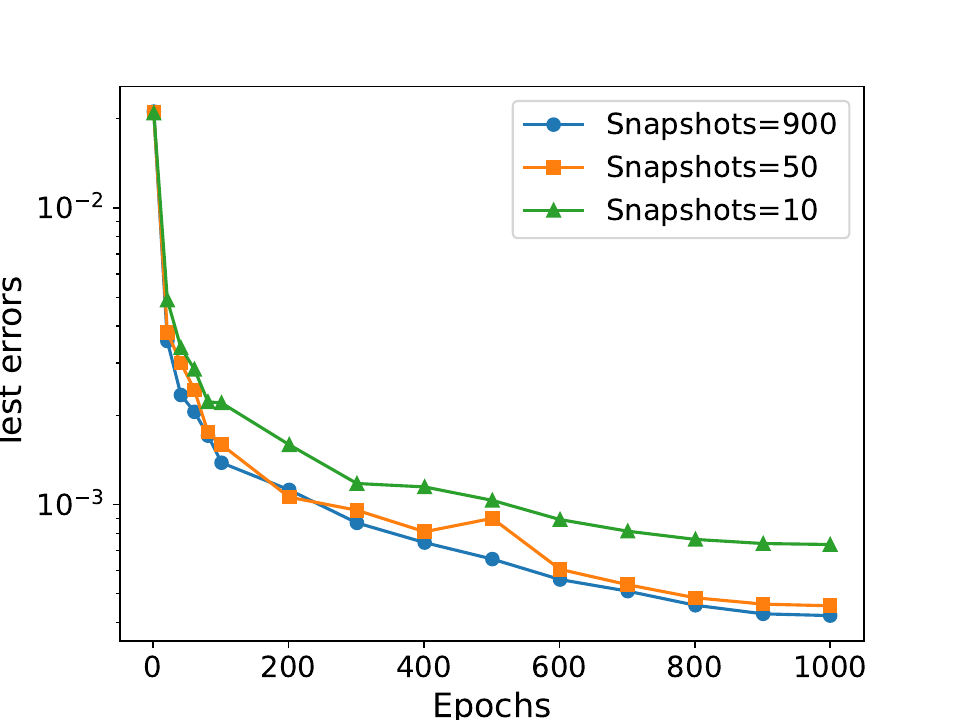}
    \caption{Error comparison between different snapshot numbers for the NLS equations (left) and the KP equations (right).}
    \label{PODNO_ablation_snapshots}
\end{figure}
Fig.~\ref{PODNO_ablation_snapshots} demonstrates that, for both the NLS and KP equations, more snapshots result in smaller errors. Moreover, we would like to highlight the result for the KP equations: although the errors for $10$ and $900$ snapshots show only a small difference, the learning process has already been dampened when using only 10 snapshots. This is because, under our initial condition setup, the region of nonuniform oscillatory solutions is small and localized, unlike the global solutions in the NLS equation. As a result, even a slight change in the numerical error could mean a significant deviation within this localized region.
\subsubsection{Effect of mode numbers on PODNO}
\begin{figure}[!ht]
    \centering
    \includegraphics[width=0.4\textwidth]{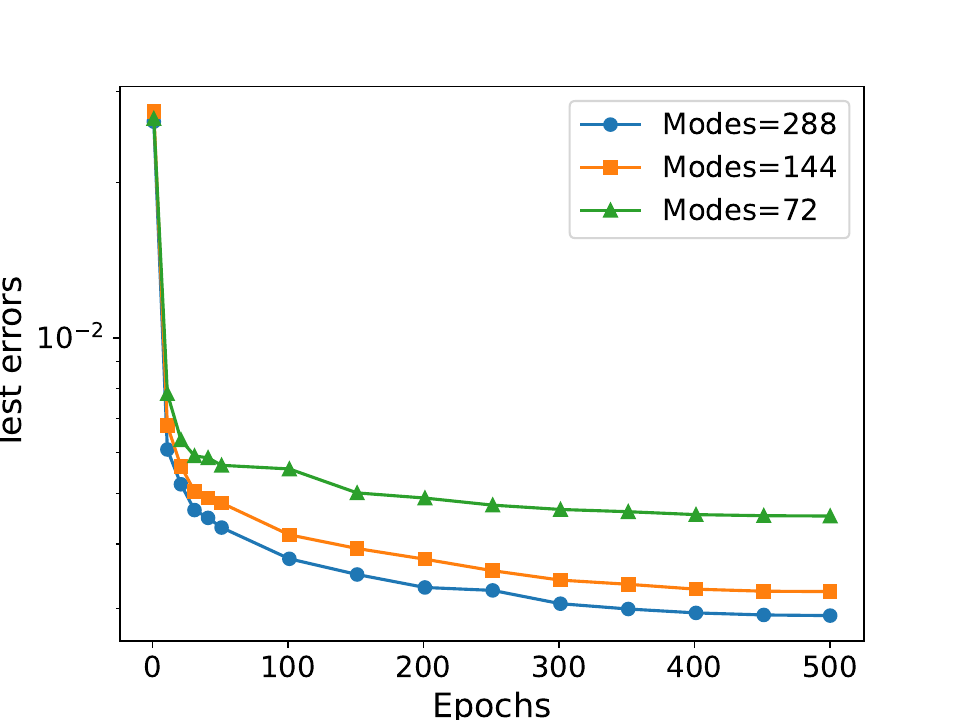}
    ~
 \includegraphics[width=0.4\textwidth]{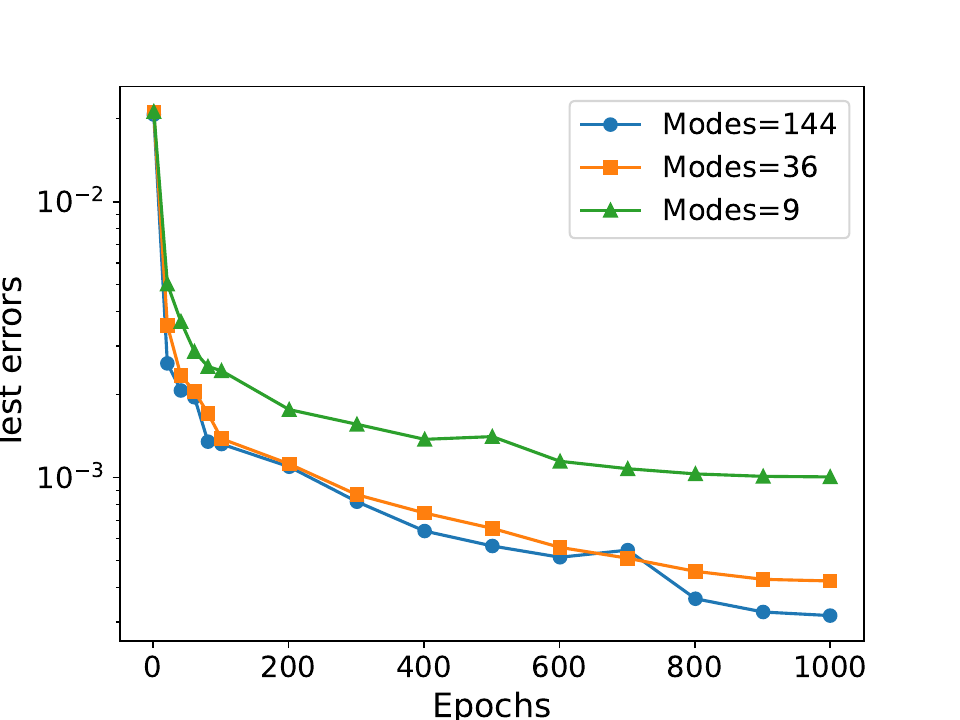}
    \caption{Error comparison between different modes for NLS equations (left) and KP equations (right).}
    \label{PODNO_ablation_modes}
\end{figure}
A portion of the error in the learning process arises from the choice of mode numbers, as all modes contribute to capturing the energy to some extent. Larger POD mode numbers capture more energy in the partial differential equations, which can lead to better predictions in theory. As shown on the left side of Fig.~\ref{PODNO_ablation_modes}, PODNO with $288$ modes or $144$ modes yields almost the same errors, as the energy capture ratios $\rho$ approach $100\%$. However, PODNO with only $72$ modes performs worse, with an energy capture ratio of $\rho=99.6\%$.

For the KP equations, shown on the right, even with only 9 modes, PODNO captures $99.8\%$ energy but still performs worse compared to the larger number of modes that capture more energy. As mentioned earlier, the focus is on the small oscillatory region, where even a small error can indicate a failure in the training. We conclude that the results for both equations align with theoretical expectations, showing that a larger mode count leads to better learning performance.

\section{Conclusion}
In this paper, we leverage the fact that the normal POD basis is the optimal truncated basis in the $L^2$ metric to design a data-driven operator learning methodology, PODNO. 
Comparing with FNO, the PODNO is better tailored for high-frequency-dominated problems. For instance, our experiments on the NLS and KP equations demonstrate that PODNO outperforms FNO in both efficiency and accuracy in this setting, as illustrated in Fig.~\ref{nls_errors_fig} and Fig.~\ref{kp_errors_fig}, as well as Tables \ref{nls_errors_tab} and \ref{kp_errors_tab}. Additionally, PODNO surpasses the POD-based operator splitting scheme, as shown in Fig.~\ref{comparison PODNO & POD splitting}. 

We also explore the universality of PODNO.  When $s'\geq s>d/2$, where $d$ is the dimension of the domain $\Omega\in\mathbb{R}^d$, we establish in Theorem \ref{universal} that there exists a GSO $\mathcal{N}: H^s(\Omega; \mathbb{R}^{d_a})\rightarrow H^{s}(\Omega; \mathbb{R}^{d_u})$ to approximate a continuous operator $\mathcal{G}: H^s(\Omega; \mathbb{R}^{d_a})\rightarrow H^{s'}(\Omega; \mathbb{R}^{d_u})$ with small $H^s$ error. Thus, PODNO, as a representative example of GSO, exhibits the same universality.

Furthermore, we examine the impact of spatial discretization, the number of snapshots, and the number of retained modes on PODNO's performance. While PODNO is mesh-invariant, an overly fine mesh in high-frequency-dominated problems can lead to ill-conditioning in the singular value decomposition as the singular values approach zero; a mesh that is too coarse may lose essential information, preventing the network from learning meaningful patterns (Fig.~\ref{PODNO_ablation_mesh}). Increasing the number of snapshots helps mitigate sampling error, as snapshots represent only a subset of the full dataset and may not fully capture the system’s variability and complexity (Fig.~\ref{PODNO_ablation_snapshots}). Similarly, increasing the number of modes reduces truncation error by retaining more principal components, thereby minimizing global error (Fig.~\ref{PODNO_ablation_modes}).

A key challenge for future research is ensuring a sufficient number of representative snapshots to accurately approximate the solution space, which currently limits PODNO’s ability to outperform existing baseline models in low-frequency-dominated problems. Alternatively, exploring an adaptive POD \cite{singer2009using} method may introduce greater flexibility and enhance performance.

\bibliographystyle{siam}
\bibliography{Full}
\clearpage
\begin{appendix}
\section{Proof for the Universal Approximation Theorem}\label{proof}
In Appendix \ref{proof}, we begin by providing the preliminaries for the proof, which include the mathematical notation, some definitions, and some propositions. We then proceed to present the proof for the universal approximation theorem for the GSO, following the establishment of a preparatory lemma.
\subsection{Preliminaries}\begin{table}[!htbp]
    \centering
    \caption{Glossary of mathematical notation}\label{Notation}

    {\footnotesize  % Adjust font size to small
    \begin{tabular}{l l}
    \toprule
    Symbol & Description \\
    \midrule
    $d$&dimension of the domain $\Omega\in\mathbb{R}^d$\\
    $s$&degree of smoothness (or regularity)\\
    $N$&number of the modes retained\\
    $\phi_k$ &$k$-th orthonormal basis, $\Phi_k:=\{\phi_k\}_{k=1}^N$\\
    $v$, $a$ & function in the phisical space\\
    $\hat v_k$, $\hat a_k$ & $k$th-expansion coefficient of $v$ or $a$ in $\phi_k$\\
    $P_N$&$H^s$-orthonormal projection onto $H^s_N$, see \eqref{pn}\\
    $\mathcal{G}$&continuous operator to be learned by the GSO\\
    $\mathcal{G}_N$& projection of the continuous operator $\mathcal{G}$ acting on $H_N^s$, see \eqref{gn}\\
    $\hat{\mathcal{G}}_N$&orthonormal conjugate/ the orthonormal dual operator of $\mathcal{G}_N$, see \eqref{gn decomposition}\\
    $\varPi_N$& truncated orthonormal transform that transfers the functions in $H^s$ \\
    &to coefficient space with regard to the basis $\Phi_N$\\
    $\varPi_N^{-1}$& inverse of $\varPi_N$ that transfers the coefficients back to physical space\\
    $\varPi_N^\dagger$ &operator that maps a function to a vector of constant functions representing\\
    &its coefficients\\
    $(\varPi_N^{-1})^\dagger$ &operator that reconstructs a function in the physical space by \\
    &summing the basis functions weighted by constant coefficients functions\\
    $\hat{\mathcal{G}}_N^\dagger$ &operator that applies $\hat{\mathcal{G}}_N$ pointwise to vector of constant functions\\
    $\mathcal{N}$&GSO\\
    $\mathcal{P}$&lifting operator in $\mathcal{N}$\\
    $\mathcal{Q}$&projection operator in $\mathcal{N}$\\
    $\mathcal{L}$&kernel integration layer without an activation, $\mathcal{L}=\varPi_N^{-1}\circ\mathcal{R}\circ\varPi_N+\mathcal{W}$\\
    $\mathcal{R}$&orthonormal multiplier; linear transformation applied to the coefficients\\
    $\mathcal{W}$&pointwise linear transformation, $\mathcal{W}[v](x)=Wv(x)+b(x)$\\
    $\mathcal{M}_f$, $\mathcal{M}_{\hat{\mathcal{G}}_N}$ &ordinary neural networks that approximates the function $f$ or $\hat{\mathcal{G}}_N$\\
    $\mathcal{M}$ &composition of activations and GSO layers that lifts ordinary neural networks\\
    &to act on function space\\
    % $W$&the linear transformation of the pointwise affine mapping\\
    % $b$&the bias of the pointwise affine mapping\\ 
    % $\varPi$ &the orthonormal projection operator\\
    % $\varPi^{-1}$ &the inverse of the orthonormal projection operator\\
    $\sigma$ &activation function\\
    $I^N_r$ & N-dimensional hypercube $I^N_r = [-r, r]^N \subset \mathbb{R}^N$\\
    $1(x)$ & constant function taking value 1 in $\Omega$\\
    \bottomrule
    \end{tabular}
    }
\end{table}

\begin{defn}\label{defn pn}
Let $\Omega\in\mathbb{R}^d$, $\phi_k$ are the orthonormal basis in $H^s(\Omega)$:
\begin{equation*}
    \langle \phi_i,\phi_j\rangle_{H^s}=\delta_{ij}=
    \begin{cases}
        1,\quad\text{if~}i=j,\\
        0,\quad\text{if~}i\neq j.
    \end{cases}
\end{equation*}
We denote the projection to the space spanned by the first $N$ basis $\{\phi_k\}_{k=1}^N$, $H_N^s\left(\Omega\right)$, by
\begin{equation*}
P_N: H^s(\Omega) \rightarrow H_N^s(\Omega), \quad v \mapsto P_N v,    
\end{equation*}
the $H^s$-orthogonal projection onto $H_N^s(\Omega)$; or more explicitly, given $v\in H^s(\Omega)$, let $\hat v_k=\langle v_k,\phi_k\rangle_{H^s}$, then
\begin{equation}\label{pn}
P_N(v)=\sum_{1\leq k\leq N} \hat{v}_k \phi_k(x).  
\end{equation}
\end{defn}

\begin{defn}\label{defn gn}
We define the following operator,
\begin{equation}\label{gn}
\mathcal{G}_N: H^s(\Omega) \rightarrow H^s(\Omega), \quad \mathcal{G}_N[a]:=P_N \mathcal{G}[P_N a],    
\end{equation}
with $P_N$ being the projection defined in Definition \ref{defn pn}.

We can also define an orthonormal conjugate or orthonormal dual operator of $\mathcal{G}_N$ in the form
$\hat{\mathcal{G}}_N: \mathbb{R}^{N} \rightarrow \mathbb{R}^N$ such that
\begin{equation}\label{gn decomposition}
\mathcal{G}_N[a]=\varPi_N^{-1} \circ \hat{\mathcal{G}}_N \circ \varPi_N[P_N a],  
\end{equation}
holds for all real-valued $a \in H^s(\Omega)$. Here, $\varPi_N$ is the orthonormal transform that transfers the function from physical space $H^s(\Omega)$ to coefficients space $\mathbb{R}^N$, and $\varPi_N^{-1}$ is the inverse orthonormal transform that transfers the coefficients back to the physical space, respectively.
\end{defn}
If we denote $\mathcal{G}_N[a] := a'$, the transformation can be expressed as the following chain:
\begin{equation*}
    a \xrightarrow{P_N}\sum_{1\leq k\leq N}\hat{a}_k\phi_k\xrightarrow{\varPi_N} (\hat{a}_k)_{k=1}^N\xrightarrow{\hat{\mathcal{G}}_N} (\hat{a}'_k)_{k=1}^N\xrightarrow{\varPi_N^{-1}}a'.
\end{equation*}
Alternatively, we can define a chain where each operator maps between functions:
\begin{equation*}
    a \xrightarrow{P_N}\sum_{1\leq k\leq N}\hat{a}_k\phi_k\xrightarrow{\varPi_N^\dagger} (\hat{a}_k1(x))_{k=1}^N\xrightarrow{\hat{\mathcal{G}}_N^\dagger} (\hat{a}'_k1(x))_{k=1}^N\xrightarrow{(\varPi_N^{-1})^\dagger}a'.
\end{equation*}
In this formulation, instead of mapping functions to their coefficient vectors, we use $\varPi_N^\dagger$ to map them to vectors of constant functions, where each component is a constant function whose value matches the corresponding coefficient. The operator $\hat{\mathcal{G}}_N^\dagger$ acts on these vectors of constant functions, mapping them to new vectors of constant functions. Finally, $(\varPi_N^{-1})^\dagger$ maps the result back to the physical space. This gives the alternative decomposition of $\mathcal{G}_N$, which we show later can be approaximated by GSO,
\begin{equation}\label{gn dagger}
\mathcal{G}_N(a) = (\varPi_N^{-1})^{\dagger} \circ \hat{\mathcal{G}}_N^\dagger \circ \varPi_N^\dagger[P_N a].
\end{equation}
Here we list several results in general approximation theory, which will be used in the proof.
\begin{prop}[Sobolev embedding theorem \cite{evans2022partial}]\label{embedding}
Let $\Omega\in\mathbb{R}^d$, and let $s>d/2$, we have a compact embedding $H^s(\Omega)\hookrightarrow C(\Omega)$ into the space of continuous functions. In particular, there exists a constant $C_1=C(s,d,\Omega)>0$, such that
\begin{equation*}
\|v\|_{L^{\infty}} \leq C_1\|v\|_{H^s}, \quad \forall~v \in H^s(\Omega) .    
\end{equation*}
\end{prop}
\begin{prop}[Sobolev product estimate \cite{kato1988commutator}]\label{product estimate}
Let $\Omega \subset \mathbb{R}^d$, and let $ s>d/2$. Suppose that $f, g \in H^s(\Omega) \cap L^\infty(\Omega)$. Then there exists a constant $C_2 = C(s, d, \Omega) > 0$ such that:
\begin{equation*}
\|fg\|_{H^s(\Omega)} \leq C_2 \left( \|f\|_{L^\infty(\Omega)} \|g\|_{H^s(\Omega)} + \|f\|_{H^s(\Omega)} \|g\|_{L^\infty(\Omega)} \right)
\end{equation*}
    
\end{prop}
\begin{prop}[Universal approximation theorem for ordinary neural networks \cite{hornik1991approximation}]\label{ordinary uat} 
$\forall~f\in C(\mathbb{R}^n;\mathbb{R})$, suppose the activation function $\sigma$ is continuous, bounded, and non-constant, $\forall~\varepsilon>0$, there exists  an ordinary neural network $\mathcal{M}_f$
\begin{equation*}
    \|f(x)-\mathcal{M}_f(x)\|_{L^\infty (K_0)}<\epsilon,
\end{equation*}
where $K_0$ is a compact domain in $\mathbb{R}^n$.
\end{prop}

For any continuous operator $\mathcal{G}$, there exists a corresponding operator $\mathcal{G}_N$ as defined in Definition~\ref{defn gn}. With the above preparations in place, we proceed to show the following: in Lemma~\ref{gn gso}, we demonstrate that a GSO $\mathcal{N}$ can be constructed to approximate $\mathcal{G}_N$; and in Theorem~\ref{proof_universal}, we prove that $\mathcal{G}_N$ can be used to approximate $\mathcal{G}$. By combining these results, we conclude that $\mathcal{G}$ can be approximated by a suitable GSO $\mathcal{N}$.
\subsection{GSO approximations of $\mathcal{G}_N$}
We begin by presenting the lemma stating that $\mathcal{G}_N$ could be approximated by a GSO.
\begin{lemma}\label{gn gso}
Let $\Omega\subset\mathbb{R}^d$ and $s>d/2$. Suppose $K\subset H^s(\Omega)$ is a compact set. Then $\forall~\varepsilon>0$, $\exists$ 
a GSO $\mathcal{N}: H^s(\Omega)\mapsto H^s(\Omega)$ such that
\begin{equation*}
    \sup_{v\in K}\|\mathcal{G}_N[a]-\mathcal{N}[a]\|_{H^s}<\varepsilon,
\end{equation*}    
where the operator $\mathcal{G}_N$ is defined in Definition \ref{defn gn}.
\end{lemma}
\begin{proof}
We consider the decomposition $\mathcal{G}_N(a) = (\varPi_N^{-1})^{\dagger} \circ \hat{\mathcal{G}}_N^\dagger \circ \varPi_N^\dagger[P_N a]$ as in \eqref{gn dagger}, and aim to construct a GSO that simulates this decomposition.

\textbf{(Step 1)} At first, we show that $\hat{\mathcal{G}}_N^\dagger\circ\varPi_N^\dagger\circ P_N[a]$ can be approximated by the composition of kernel integration layers with activation and a lifting layer, $\mathcal{L}_{L-2}\circ\dots\circ\mathcal{L}_0\circ\mathcal{P}$.

We define the row selection matrices $E_k^{(1)}\in\mathbb{R}^{N\times N}$ such that the $(i,j)$-th entry of $E_k^{(1)}$ is given by
\begin{equation*}
(E_k^{(1)})_{ij} = \begin{cases}
1,& \text{if }~i=1~\text{and}~j=k,\\
0,& \text{otherwise}.
\end{cases}    
\end{equation*}
Then we could define a tensor $E^{(1)}\in\mathbb{R}^{N\times N\times N\times N}$ such that 
\begin{equation*}
    \label{E1}E^{(1)}_{ijkm}=\begin{cases}
(E_k^{(1)})_{ij},& \text{if }~k=m,\\
0,& \text{otherwise}.
\end{cases}  
\end{equation*}
We notice that $\forall~a\in K$, the following equation holds:
\begin{equation*}
\varPi_N^{-1}E_k^{(1)}\varPi_N(a(x))_{k=1}^N=\hat{a}_k1(x),   
\end{equation*}
where $\hat{a}_k$ denotes the coefficient when projecting the function $a$ to the $k$-th orthonormal basis $\phi_k$.
From this observation, we could construct a composition of a GSO layer without activation and a lifting layer $\mathcal{L}_0^{(1)}\circ\mathcal{P}$ to approximate $\varPi_N^\dagger\circ P_N$, with 
\begin{equation*}
\begin{aligned}
    \mathcal{P}:H^s(\Omega;\R)&\rightarrow H^s(\Omega;\R^N)\\
    a&\mapsto (a)_{k=1}^N.
\end{aligned}
\end{equation*}
In the GSO layer $\mathcal{L}_0^{(1)}$ without the activation, the weight $W_0^{(1)}=0$, the bias $b_0^{(1)}=0$, the orthonormal multiplier $\mathcal{R}_0^{(1)}$ represented by $E^{(1)}$. 

However, to construct a (not final) GSO layer, we lack the activation function, thus, we need more processing, and we notice that we may consider approximating $\hat{\mathcal{G}}_N^\dagger\circ\varPi_N^\dagger\circ P_N$ instead of approximating $\varPi_N^\dagger\circ P_N$.

By the Riesz basis property, $\exists~C_3>0$, such that 
\begin{equation*}
    \|\|\varPi_N^\dagger[P_Na]\|_{L^\infty}=\|(\hat{a}_k1(x))_{k=1}^N\|_{L^\infty}=\max_{1\leq k\leq N}|\hat{a}_k|<C_3.
\end{equation*}

By the universal approximation theorem for ordinary neural networks (Proposition \ref{ordinary uat}), $\forall~\varepsilon>0$, $\exists$ an ordinary neural network $\mathcal{M}_{\hat{\mathcal{G}}_N}$ such that
\begin{equation*}
    \sup_{y\in I_{C_3}^N}|\hat{\mathcal{G}}_N(y)-\mathcal{M}_{\hat{\mathcal{G}}_N}(y)|<\varepsilon.
\end{equation*}
The $l$-th layer of $\mathcal{M}_{\hat{\mathcal{G}}_N}$ is of the form
\begin{equation}\label{weights and bias}
    y_{(l)}\mapsto \sigma (W_l^\mathcal{M}y_{(l)}+b_l^\mathcal{M}).
\end{equation}

Since $\mathcal{M}_{\hat{\mathcal{G}}_N}$ provides a uniform approximation for $\hat{\mathcal{G}}_N$ on the compact set $I_{C_3}^N$, we extend this approximation pointwise to vectors of constant functions in $H^s(\Omega; I_{C_3}^N)$. Specifically, for any $a^*:=\varPi_N^\dagger[P_Na] \in H^s(\Omega; I_{C_3}^N)$, we can construct the composition of activations and GSO layers $\mathcal{M}:=\sigma\circ\mathcal{L}_{L-2}^{(2)}\circ\dots\circ\sigma\circ\mathcal{L}_1^{(2)}\circ\sigma\circ\mathcal{L}_0^{(2)}$ to approximate $\hat{\mathcal{G}}_N^\dagger$, and each layer $\mathcal{L}_l^{(2)}$ is of the form
\begin{equation}\label{mid layer}
    v_{(l)}(x)\mapsto W_l^\mathcal{M}v_{(l)}(x)+b_l^\mathcal{M}(x),
\end{equation}
if we set $v_{(0)}:=a^*$.

Since we notice that 
\begin{equation*}
W_0^\mathcal{M}\varPi_N^{-1}E^{(1)}_k\varPi_N=\varPi_N^{-1}E_k^{(1)}W_0^\mathcal{M}\varPi_N,    
\end{equation*}
we could define a new GSO layer $\mathcal{L}_0:=\mathcal{L}_0^{(2)}\circ\mathcal{L}_0^{(1)}$ of the form:
\begin{equation*}
(a)_{k=1}^N(x)\mapsto\sigma[\varPi_N^{-1}\circ\mathcal{R}_0\circ\varPi_N\circ P_N[(a)_{k=1}^N]+W_{0}[(a)_{k=1}^N]+b_{0}](x),
\end{equation*} 
where $b_{0}=b_{0}^\mathcal{M}$ and $\mathcal{R}_0$ represented by the tensor $\tilde{E}^{(1)}$ such that
\begin{equation*}
    \tilde{E}^{(1)}_{ijkm}=\begin{cases}
(E_k^{(1)}W_0)_{ij},& \text{if }~k=m,\\
0,& \text{otherwise}.
\end{cases}  
\end{equation*}

In addition, taking the following kernel integration layers $\mathcal{L}_l=\mathcal{L}_l^{(2)}$ of the form \eqref{mid layer}, we can then approximate $\hat{\mathcal{G}}_N^\dagger\circ\varPi_N^\dagger\circ P_N$ by $\sigma\circ\mathcal{L}_{L-2}\circ\dots\circ\sigma\circ\mathcal{L}_1\circ\sigma\circ\mathcal{L}_0\circ\mathcal{P}$.

We recall  $\mathcal{M}:=\sigma\circ\mathcal{L}_{L-2}^{(2)}\circ\dots\circ\sigma\circ\mathcal{L}_1^{(2)}\circ\sigma\circ\mathcal{L}_0^{(2)}$ and $v=\varPi_N^\dagger\circ P_N[a]=\mathcal{L}_0^{(1)}\circ\mathcal{P}[a]$. From this, we obtain
\begin{equation*}
\begin{aligned}
&\sup_{a\in K}\|\hat{\mathcal{G}}_N^\dagger\circ\varPi_N^\dagger\circ P_N[a]-\sigma\circ\mathcal{L}_{L-2}\circ\dots\circ\sigma\circ\mathcal{L}_1\circ\sigma\circ\mathcal{L}_0\circ\mathcal{P}[a]\|_{L^2}\\
=&\sup_{a\in K}\|\hat{\mathcal{G}}_N^\dagger\circ\varPi_N^\dagger\circ P_N[a]-\sigma\circ\mathcal{L}_{L-2}^{(2)}\circ\dots\circ\sigma\circ\mathcal{L}_1^{(2)}\circ\sigma\circ\mathcal{L}_0^{(2)}\circ\mathcal{L}_0^{(1)}\circ\mathcal{P}[a]\|_{L^2}\\
=&\sup_{a\in K}\|\hat{\mathcal{G}}_N^\dagger\circ\varPi_N^\dagger\circ P_N[a]-\mathcal{M}\circ\mathcal{L}_0^{(1)}\circ\mathcal{P}[a]\|_{L^2}\\
=&
\sup_{a^*\in H^s(\Omega, I_{C_3}^N)}\|\hat{\mathcal{G}}_N^\dagger[a^*]-\mathcal{M}[a^*]\|_{L^2}\\
=&\sup_{a^*\in H^s(\Omega, I_{C_3}^N)}\big( \int_{\Omega} |\hat{\mathcal{G}}_N^\dagger[a^*](x)-\mathcal{M}[a^*](x)|^2 dx \big)^{1/2}\\
\leq&\big( \int_{\Omega} \sup_{a^*(x)\in I_{C_3}^N}|\hat{\mathcal{G}}_N(a^*(x))-\mathcal{M}_{\hat{\mathcal{G}}_N}(a^*(x))|^2 dx \big)^{1/2}\\
<&\big( \int_{\Omega} \varepsilon^2 dx \big)^{1/2}\\
=&|\Omega|^{1/2}\varepsilon
\end{aligned}   
\end{equation*}
where $|\Omega|$ denotes the Lebesgue measure of $\Omega$.

Since ${\mathcal{G}}_N^\dagger[a^*]$ and $\mathcal{M}[a^*]$ are just vectors of constant functions, we can show that the errors of each component are also bounded,
\begin{equation*}
\begin{aligned}
\|(\hat{\mathcal{G}}_N^\dagger[a^*])_k-(\mathcal{M}[a^*])_k\|_{H^s}&=\|(\hat{\mathcal{G}}_N^\dagger[a^*])_k-(\mathcal{M}[a^*])_k\|_{L^2}\\
&\leq\|\hat{\mathcal{G}}_N^\dagger[a^*]-\mathcal{M}[a^*]\|_{L^2}\\
&=\|\hat{\mathcal{G}}_N^\dagger[a^*]-\mathcal{M}[a^*]\|_{H^s}\\
&<|\Omega|^{1/2}\varepsilon,    
\end{aligned} 
\end{equation*}
where $1\leq k\leq N$.

\textbf{(Step 2)} Next, we show that $(\varPi_N^{-1})^\dagger$ can be approximated by the composition of a GSO layer without activation and a projection layer $\mathcal{Q}\circ\mathcal{L}_{L-1}$.

We define different row selection matrices $E_k^{(2)}$ such that the $(i,j)$-th entry of $E_k^{(2)}$ is given by
\begin{equation*}
(E_k^{(2)})_{ij} = \begin{cases}
1,& \text{if }~j=1~\text{and}~i=k,\\
0,& \text{otherwise}.
\end{cases}    
\end{equation*}
Then we could define a tensor $E^{(1)}\in\mathbb{R}^{N\times N\times N\times N}$ such that 
\begin{equation*}
    E^{(2)}_{ijkm}=\begin{cases}
(E_k^{(2)})_{ij},& \text{if }~k=m,\\
0,& \text{otherwise}.
\end{cases}  
\end{equation*}
We notice that $\forall$ vector of constant functions $\eta(x)=(\eta_k1(x))_{k=1}^N$,
\begin{equation*}
    (\varPi_N^{-1})^\dagger[\eta](x)=\varPi_N^{-1}\mathrm{diag}(E_k^{(2)})\varPi_N\eta(x)=(\eta_k\phi_k(x))_{k=1}^N=(\eta_k1(x)\phi_k(x))_{k=1}^N
\end{equation*}
From this observation, we define the layer $\mathcal{L}_{L-1}$ with the form
\begin{equation*}
    v\mapsto\varPi_N^{-1}\circ\mathcal{R}_{L-1}\circ\varPi_N\circ[v]
\end{equation*}
where $\mathcal{R}_{L-1}$ is represented by $E^{(2)}$.
Thus $\mathcal{Q}\circ\mathcal{L}_{L-1}=(\varPi_N^{-1})^\dagger$ with
\begin{equation*}
\begin{aligned}
    \mathcal{Q}:H^s(\Omega;\R^N)&\rightarrow H^s(\Omega;\R)\\
    (h_k)_{k=1}^N &\mapsto \sum_{k=1}^N h_k.
\end{aligned}
\end{equation*}

(\textbf{Step 3}) To summarize, we construct a GSO $\mathcal{N}$ composed of the following components.
The first component is the lifting layer, defined as
\begin{equation*}
\begin{aligned}
    \mathcal{P}:H^s(\Omega;\R)&\rightarrow H^s(\Omega;\R^N)\\
    a&\mapsto (a)_{k=1}^N.
\end{aligned}
\end{equation*}

This is followed by the composition of activations and GSO layers.
For the first layer $l=0$, $\mathcal{R}_0$ is represented by $\tilde{E}^{(1)}$,
    the weight $W_0=0$ and the bias $b_0(x)=b_0^\mathcal{M}$ is the same as \ref{weights and bias}.
    
For the intermediate layers $l=1,2,\dots,L-2$, $\mathcal{R}_l\equiv 0$, the weight $W_{l}=W_{l}^\mathcal{M}$ and the bias $b_{l}(x)=b_{l}^\mathcal{M}$ the same as \ref{weights and bias}.

For the final layer $l=L-1$, $\mathcal{R}_{L-1}$ is represented by $E^{(2)}$, the weight $W_{L-1}=0$ and the bias $b_{L-1}(x)=0$.
Finally, the projection layer maps back to the target dimension:
\begin{equation*}
\begin{aligned}
    \mathcal{Q}:H^s(\Omega;\R^N)&\rightarrow H^s(\Omega;\R)\\
    (h_k)_{k=1}^N &\mapsto \sum_{k=1}^N h_k.
\end{aligned}
\end{equation*}
We now estimate the approximation error between $\hat{\mathcal{G}}_N$ and $\mathcal{N}$. There exists a constant $C>0$ such that
\begin{equation*}
\begin{aligned}
&\sup_{a\in K}\|\hat{\mathcal{G}}_N[a]-\mathcal{N}[a]\|_{H^s}\\
=& \sup_{a\in K}\|(\varPi_N^{-1})^\dagger\circ\hat{\mathcal{G}}_N^\dagger\circ\varPi_N^\dagger[P_Na]-\mathcal{Q}\circ\mathcal{L}_{L-1}\circ\dots\circ\sigma\circ\mathcal{L}_{1}\circ\sigma\circ\mathcal{L}_{0}\circ\mathcal{P}[a]\|_{H^s}\\
=&\sup_{a^*\in H^s(\Omega, I_{C_3}^N)}\|(\varPi_N^{-1})^\dagger\circ\hat{\mathcal{G}}_N^\dagger[a^*]-\mathcal{Q}\circ\mathcal{L}_{L-1}\circ\mathcal{M}[a^*]\|_{H^s}\\
=&\sup_{a^*\in H^s(\Omega, I_{C_3}^N)}\|\sum_{k=1}^N\phi_k(\hat{\mathcal{G}}_N^\dagger[a^*])_k-\sum_{k=1}^N\phi_k(\mathcal{M}[a^*])_k\|_{H^s}\\
\leq &N\sup_{a^*\in H^s(\Omega, I_{C_3}^N)}\max_{1\leq k\leq N}\|\phi_k(\hat{\mathcal{G}}_N^\dagger[a^*])_k-\phi_k(\mathcal{M}[a^*])_k\|_{H^s}\\
\leq& NC\sup_{a^*\in H^s(\Omega, I_{C_3}^N)}\max_{1\leq k\leq N}\|(\hat{\mathcal{G}}_N^\dagger[a^*])_k-(\mathcal{M}[a^*])_k\|_{H^s}\\
<&NC|\Omega|^{1/2}\varepsilon,
\end{aligned}
\end{equation*}
the second-to-last line follows from the Sobolev embedding theorem \ref{embedding} and the Sobolev product estimate \ref{product estimate}. 

Since $\varepsilon>0$ is arbitrary, we conclude the claim of this lemma follows.
\end{proof}
\begin{remark}
    It is worth mentioning that, the proof of Lemma \ref{gn gso} is constructive and hence the size of the tensors (for instance \eqref{E1}) may not be optimal. Due to the complexity of training neural networks, the optimal construction may be discretization-dependent.
\end{remark}
\subsection{Proof of Theorem \ref{universal}}\label{proof_universal}
% \noindent \textbf{Proof of Theorem \ref{universal}.}\label{proof_universal}
Here we recall the statement of the universal approximation theorem for ease of reading.
\begin{thm}[Theorem \ref{universal}]
Let $\Omega\in\mathbb{R}^d$, and let $$\mathcal{G}: H^s(\Omega; \mathbb{R}^{d_a})\rightarrow H^{s'}(\Omega; \mathbb{R}^{d_u})$$ be a continuous operator, where $s'\geq s>d/2$. Suppose $K \subset H^s(\Omega; \mathbb{R}^{d_a})$ is a compact subset. Then, $\forall~\varepsilon>0$, $\exists$ a GSO, $\mathcal{N}: H^s(\Omega; \mathbb{R}^{d_a}) \rightarrow H^s(\Omega; \mathbb{R}^{d_u})$, satisfying
\begin{equation*}
    \sup _{a \in K}\|\mathcal{G}[a]-\mathcal{N}[a]\|_{H^s} < \varepsilon.
\end{equation*}
\end{thm}
\begin{proof}
Throughout the proof, we set $d_a=d_u=1$ for convenience. The general case with $d_a,~d_u>1$ follows analogously. 

We break the proof into two parts. In the first part, we will show that $\exists~N \in \mathbb{N}$, such that
\begin{equation*}
\|\mathcal{G}[a]-\mathcal{G}_N[a]\|_{H^s} < \varepsilon, \quad \forall~a \in K,  
\end{equation*}
where the operator $\mathcal{G}_N$ is defined in \eqref{gn}; In the second part, we show that there exists a GSO: $\mathcal{N}: H^s(\Omega;\mathbb{R}) \rightarrow H^s(\Omega;\mathbb{R})$, such that
\begin{equation*}
\sup _{a \in K}\|\mathcal{G}_N[a]-\mathcal{N}[a]\|_{H^s}<\varepsilon.    
\end{equation*}
This result has already been established in Lemma \ref{gn gso}.

Then we could provide an approximation of $\mathcal{G}$, such that
\begin{equation*}
\sup _{a \in K}\|\mathcal{G}[a]-\mathcal{N}[a]\|_{H^s} <2\varepsilon.    
\end{equation*}
Since $\varepsilon>0$ was arbitrary, the claim follows from this.

Here we focus on the first part.
We note that since $K \subset H^s(\Omega;\mathbb{R})$ is compact, $\tilde{K}$ defined by
\begin{equation*}
    \tilde{K}:=K \cup \bigcup_{N \in \mathbb{N}} P_N K,
\end{equation*}
is also compact. 
Let $\{u_n\}$ be a sequence in $\tilde K$.
\begin{itemize}
    \item If infinitely many of the $u_n$ lie in $K$, then since $K$ is compact, we can extract a convergent subsequence that converges to some $u\in K\subset\tilde K$.
    \item If infinitely many of the $u_n$ lie in $P_N K$ for a $N\in \mathbb{N}$, then since $P_N K$ is compact, we can extract a convergent subsequence that converges to some $u\in P_N K\subset\tilde K$.
    \item Suppose that for each $N \in \mathbb{N}$, only finitely many terms of the sequence $\{u_n\} \subset \tilde K$ lie in $P_N K$. We construct a subsequence $\{u_{n'_N}\} \subset \{u_n\}$ by selecting, for each $N$, an element $u_{n'_N} \in P_N K \setminus \bigcup_{j < N} P_j K$. Since $u_{n'_N} \in P_N K$, there exists $k_N \in K$ such that $u_{n'_N} = P_N k_N$. The sequence $\{k_N\} \subset K$ admits a convergent subsequence $\{k_{N_j}\} \to k \in K$ by compactness of $K$. Moreover, since $P_N \to I$ strongly, it follows that $P_{N_j} k_{N_j} \to k$. Hence,
\begin{equation*}
  u_{n'_{N_j}} = P_{N_j} k_{N_j} \to k \in K \subset \tilde K.
\end{equation*}
  Therefore, the original sequence $\{u_n\}$ has a subsequence converging in $\tilde K$.
\end{itemize}
% \zw{what if none of the statement 'infinitely many ... lies in ...' happen. You cannot easily say  an infinite union of compact set is compact...}
To conclude, $\tilde K=K \cup \bigcup_{N \in \mathbb{N}} P_N K$ is compact.

Since $\mathcal{G}$ is continuous, its restriction to $\tilde{K}$ is uniformly continuous, i.e. there exists a modulus of continuity $\omega:[0, \infty) \rightarrow[0, \infty)$, such that
\begin{equation*}
\|\mathcal{G}[a]-\mathcal{G}([a^{\prime}])\|_{H^s} \leq \omega(\|a-a^{\prime}\|_{H^s}), \quad \forall a,~a^{\prime} \in \tilde{K} .
\end{equation*}
From the definition of the projection $\mathcal{G}_N$ in Definition \ref{defn gn}, we have
\begin{equation*}
\begin{aligned}
\|\mathcal{G}[a]-\mathcal{G}_N[a]\|_{H^s}
&\leq\|\mathcal{G}[a]-P_N \mathcal{G}[a]\|_{H^s}+\|P_N \mathcal{G}[a]-P_N \mathcal{G}[P_N [a]]\|_{H^s} \\
& \leq\|\mathcal{G}[a]-P_N \mathcal{G}[a]\|_{H^s}+\|\mathcal{G}[a]-\mathcal{G}[P_N a]\|_{H^s} \\
& \leq \sup _{v \in \mathcal{G}(\tilde{K})}\|(1-P_N) v\|_{H^s}+\omega(\sup _{a \in \tilde{K}}\|(1-P_N) a\|_{H^s}),
\end{aligned}    
\end{equation*}
where
\begin{equation*}
    \limsup_{N \rightarrow \infty} \sup _{a \in \tilde{K}}\|(1-P_N) a\|_{H^s}=0.
\end{equation*}
Since $\tilde{K}$ is compact and $\mathcal{G}$ is a continuous operator, the image $\mathcal{G}(\tilde{K})$ is also compact,
\begin{equation*}
\limsup_{N \rightarrow \infty} \sup _{u \in \mathcal{G}(\tilde{K})}\|(1-P_N) v\|_{H^s}=0.
\end{equation*}
This limitation leads to,
$\forall~\varepsilon>0$, $\exists~N_0\in\mathbb{N}$, $\forall~N>N_0$,
\begin{equation*}
\left\|\mathcal{G}[a]-\mathcal{G}_N[a]\right\|_{H^s} < \varepsilon, \quad \forall a \in K \subset \widetilde{K}.    
\end{equation*}
\end{proof}
\section{POD method in discrete}\label{app:pod}
In this section, we list the POD method in discrete for completeness, where 
we employ the method of snapshots \cite{kunisch2001galerkin}. Specifically, we organize both the discretized function value data from the input space $\mathcal{A}$ and the solution space $\mathcal{U}$ into a snapshot matrix, where each column represents a snapshot vector. Assuming the domain is discretized into a resolution of $\nmesh$ and we have $M$ snapshots, the resulting snapshot matrix $X$ is of size $\nmesh\times M$. The goal is to approximate the eigenspectrum of the dataset by using the covariance matrix $C=\frac{1}{M}XX^T$. 
% Disregarding the constant, $C=XX^T$ captures the correlations between the spatial points across the snapshots.
Consequently, the Fredholm integral \eqref{Fredholm integral} can be expressed in its discretized form as:
\begin{equation}\label{discrere Fredholm integral}
        C\Phi=\Phi\Lambda,
\end{equation}
where $\Lambda$ represents the diagonal matrix of eigenvalues of $C$, and $\Phi$ denotes the corresponding orthonormal basis matrix, which can be approximated by the principal left singular vectors obtained via SVD. Assuming we are interested in the $N$ principal left singular vectors, we can partition the singular vector matrices $U$,$V$, and the singular value matrix $\Sigma$ into 2, 2, and 4 blocks as follows:
\begin{equation}\label{SVD}
    X=U\Sigma V^T \approx\begin{bmatrix}
    U_N&U_{\nmesh-N}
    \end{bmatrix}\begin{bmatrix}
    \Sigma_N&0\\
    0&0
    \end{bmatrix}\begin{bmatrix}
    (V_N)^T\\
    (V_{M-N})^T
    \end{bmatrix}=U_N\Sigma_N(V_N)^T.
\end{equation}

Then \eqref{discrere Fredholm integral} has an approximating equivalent form:
\begin{equation*}
\begin{aligned}
&MC=XX^T\approx(U_N\Sigma_N (V_N)^T)(V_N\Sigma _N(U_N)^T)=U_N(\Sigma_N)^2(U_N)^T\\\Rightarrow&MCU_N\approx U_N(\Sigma_N)^2.    
\end{aligned}
\end{equation*}
The POD basis matrix $\Phi_N=[\phi_1,\phi_2,\dots, \phi_N]\approx U_N\in\mathbb{R}^{\nmesh\times N}$ retains the first $N$ singular vectors and thus captures the most significant energy for the $N$ modes. Data from the physical space can then be projected onto the subspace spanned by the columns of $H_N$, enabling an efficient representation of the data with reduced dimensionality.
\section{Splitting spectral methods}\label{app. spliting}
We use the splitting spectral method to solve the NLS equation of the form
\begin{equation*}
    \mathrm{i}\partial_t u+\Delta u+Vu=\frac{2}{\epsilon}(|u|^{\epsilon}-1)u,\\
\end{equation*}
which can be reformulated (see \cite{zhang2025low}) in the form of  
\begin{equation*}
    \mathrm{i}\partial_t u=\mathcal{A} [u]+\mathcal{B} [u],
\end{equation*}
with $\mathcal{A} [u]:=-\Delta u$, $\mathcal{B} [u]:=\frac{2}{\epsilon}(|u|^{\epsilon}-1)u-Vu$. We then separate the equation into two separated flows:
\begin{subequations}
\begin{align}
    &\begin{cases}\label{split nls a}
        \mathrm{i}\partial_t z(x,y,t)=-\Delta z(x,y,t),\\
        z(x,y,0)=z_0(x,y).
    \end{cases}\\
    &\begin{cases}\label{split nls b}
        \mathrm{i}\partial_t \omega(x,y,t)=\frac{2}{\epsilon}(|\omega(x,y,t)|^{\epsilon}-1)\omega(x,y,t)-V(x,y)\omega(x,y,t),\\
        \omega(x,y,0)=\omega_0(x,y).
    \end{cases}
\end{align}
\end{subequations}
Their exact solutions in terms of the flow maps are expressed as follows:
\begin{subequations}
\begin{align}
    &\Phi^t_\mathcal{A}[z_0](x,y)=e^{\mathrm{i}t\Delta}z_0(x,y),\\
    &\Phi^t_\mathcal{B}[\omega_0](x,y)=e^{-\mathrm{i}t(\frac{2}{\epsilon}(|\omega|^{\epsilon}-1)-V)}\omega_0(x,y).
\end{align}
\end{subequations}

Let $\delta_t> 0$ be the time-stepping size and $T> 0$ be the given final time; for $0\leq {i_t}\leq N_T:= [T/\delta_t]$, we can find the approximations of the solutions through
\begin{equation*}
    u_{i_t}(x)=(\Phi_\mathcal{A}^{\delta_t}\Phi_\mathcal{B}^{\delta_t})^{i_t}[u_0](x)
\end{equation*}

We do such a separation because it is time-efficient to solve \eqref{split nls a} in the Fourier space. If we denote by
\begin{equation*}
    \mathscr{F}[u](k_x,k_y,t)=\int_{\mathbb{R}^2} u(x,y,t)e^{-\mathrm{i}(k_xx+k_yy)}\mathrm{d}x\mathrm{d}y,
\end{equation*}
the Laplacian $\Delta$ in Fourier space is only $-(k_x^2+k_y^2)$. For \eqref{split nls b}, we have to perform the FFT on the final solution instead. The transformation between Fourier space and the physical space is an isometry, thus guaranteeing the algorithm's accuracy.
\begin{algorithm}
\caption{FFT-accelerated Lie-Trotter time-splitting scheme~\cite{zhang2025low}}\label{fft splitting}
\begin{algorithmic}[1]
\State \textbf{Inputs:} initial condition $u_0$, total time steps $N_T$, final time $T$.
\For{$i_t=1,2,\dots,N_T$}
    \State Calculate the Laplacian in the Fourier space: $\Delta_\mathscr{F}=\mathscr{F}\Delta\mathscr{F}^{-1}$
    \State Map the solution of \eqref{split nls b} from the physical space to the Fourier space: $(e^{-\mathrm{i}\delta_t\cdot \mathcal{B}[u_{{i_t}-1}]}u_{{i_t}-1})_\mathscr{F}=\mathscr{F}(e^{-\mathrm{i}\frac{T}{N_T}\cdot \mathcal{B}[ u_{{i_t}-1}]}\cdot u_{{i_t}-1})$
    \State Compose the solution to the NLS equation in the Fourier space: \newline$(u_{i_t})_\mathscr{F}=(e^{i\delta_t\Delta_\mathscr{F}})\cdot(e^{-\mathrm{i}\delta_t\cdot \mathcal{B}[ u_{{i_t}-1}]}u_{{i_t}-1})_\mathscr{F}=(e^{\mathrm{i}\frac{T}{N_T}\Delta_\mathscr{F}})\cdot\mathscr{F}(e^{-\mathrm{i}\frac{T}{N_T}\cdot \mathcal{B}[ u_{{i_t}-1}]}u_{{i_t}-1})$
    \State Turn back to the physical space: $u_{i_t}=\mathscr{F}^{-1}((u_{i_t})_\mathscr{F})$
\EndFor
\State \textbf{Output:} $u_{N_T}$ at the final time step $N_T$
\end{algorithmic}
\end{algorithm}

However, since the Fourier transform is a no-loss transform, the calculation of the nonlinear part is highly complex. Since the POD basis is the best basis in the $L^2$ space, rather than using a neural operator, it is natural to consider POD acceleration in Algorithm \ref{pod-lts} as a variation of the FFT-accelerated Lie-Trotter time-splitting scheme. Even though we do not have a simple expression of the Laplacian in the POD space, only considering the principal POD modes will largely contribute to the complexity of solving \eqref{split nls b}.
\begin{algorithm}
\caption{POD-accelerated Lie-Trotter time-splitting scheme}\label{pod-lts}
\begin{algorithmic}[1]
\State \textbf{Inputs:} initial condition $u_0$, total time steps $N_T$, final time $T$, POD modes $k\ll n^2$.
\State \textbf{Snapshots preparation:}
    \Statex \quad Case 1: $X_{\text{Basis type}~1} = [u_0, u_{N_T}]$;
    \Statex \quad Case 2: $X_{\text{Basis type}~2} = [u_0, u_1, u_2, \dots, u_{N_T}]$;
    \Statex \quad Case 3: $X_{\text{Basis type}~3} = [u_0, u_1, \dots, u_{N_T}, u_1-u_0, u_2-u_1, \dots, u_{N_T}-u_{N_T-1}]$.
\State \textbf{POD space generation:} Generate the orthonormal basis by SVD of $XX^T$.
\State Truncate to obtain the POD basis $\Phi_N$ and thus define the (truncated) POD transform $\varPi_N$ and its inverse $\varPi_N^{-1}$ as shown in Appendix \ref{app:pod}.
\For{$i_t = 1,2,\dots,N_T$}
    \State Map $e^{\mathrm{i}\delta_t\Delta}$ to the coefficients space: $(e^{\mathrm{i}\delta_t\Delta})_\varPi=\varPi_N e^{\mathrm{i}\frac{T}{N_T}\Delta}\varPi_N^{-1} $
    \State Map the solution of \eqref{split nls b} to the coefficients space: \newline
    $(e^{-\mathrm{i}\delta_t\cdot \mathcal{B}[ u_{{i_t}-1}]}u_{{i_t}-1})_\varPi = \varPi_N\left(e^{-\mathrm{i}\frac{T}{N_T}\cdot \mathcal{B}[ u_{{i_t}-1}]}\cdot u_{{i_t}-1}\right)$
    \State Compose the solution in the POD space:\newline
    $(u_{i_t})_\varPi=(\varPi_N e^{\mathrm{i}\delta_t\Delta} \varPi_N^{-1}) \cdot \varPi_N(e^{-\mathrm{i}\frac{T}{N_T} \cdot \mathcal{B}[ u_{{i_t}-1}]} u_{{i_t}-1})$
    \State Return to the physical space:$u_{i_t} = \varPi_N^{-1}((u_{i_t})_\varPi)$
\EndFor
\State \textbf{Output:} $u_{N_T}$ at the final time step $N_T$
\end{algorithmic}
\end{algorithm}
\section{Properties of dispersive equations and numerical algorithms}\label{dispersive}
In dispersive equations, different frequency components of a wave travel at different speeds, causing wave packets (superpositions of different frequencies) to spread out or disperse over time. As a result, solutions to dispersive equations often exhibit oscillatory behavior. Since the kernel of FNO filters out the high-frequency components of the solution, it may struggle to accurately learn dispersive equations. This limitation motivates our proposal of PODNO. Below, we outline the two types of dispersive equations studied in our numerical experiments.
\subsection{Nonlinear Schr\"odinger (NLS) Equation\label{nls}}
Here we consider the NLS equation with power law nonlinearity and the potential:
    \begin{equation*}
    \begin{cases}
        \textrm{i}\partial_t u+\Delta u+V(x)u=\frac{2}{\epsilon}(|u|^{\epsilon}-1)u,\\    u(x,0)=u_0,\quad x\in\Omega.
    \end{cases}
    \end{equation*}
    
The nonlinear term $\frac{2}{\epsilon}(|u|^{\epsilon}-1)u$ \cite{galati2013nonlinear} is an approximation of the nonlinear term of the Logarithmic Schr\"odinger equation, $\lambda \ln|u|^2u$. As $\epsilon\rightarrow 0$, we can expand $|u|^\epsilon$ using the Taylor series:
\begin{equation*}
\begin{split}
    &|u|^\epsilon=e^{\epsilon\ln|u|}\approx1+\epsilon\ln|u|\\
    \Rightarrow& \frac{2}{\epsilon}(|u|^\epsilon-1)\approx \ln|u|^2,
\end{split}
\end{equation*}
where $O(\epsilon^2)$ represents higher-order terms. This expansion is valid for small $\epsilon$.

Small $\epsilon$ has a great role in the strength of the nonlinearity. The term $\frac{2\lambda}{\epsilon}(|u|^{\epsilon}-1)u$ becomes highly sensitive to variations in $|u|$. This amplifies the effect of nonlinearity, potentially making the equation more unstable or introducing sharp features such as steep gradients or localized spikes in $u$. For long-time dynamics, the solutions corresponding to small $\epsilon$ may exhibit rapid phase transitions, stronger dispersion relations, or collapse phenomena depending on the initial conditions and the specific form of the potential. These behaviors arise due to the dominance of the nonlinear term. In addition, the nonlinearity promotes the formation of highly localized structures such as solitons or blow-up solutions. These solutions are more sensitive to perturbations, making them less stable and more dynamic.

In Bose-Einstein Condensates (BECs), the nonlinear term describes interactions between particles in a condensate; in nonlinearity optics, the power-law nonlinearity describes the intensity-dependent refractive index of materials.

The potential $V(x)$ determines the environment in which the quantum system evolves. There are some common potential functions used in various contexts:
\begin{itemize}
    \item $V(x)=0$ describes a particle not influenced by any external force.
    \item Box potential (or infinite potential well) \cite{da2020information} 
    \begin{equation*}
    V(x)=
        \begin{cases}
            0,\quad \text{if}~0< x< L,\\
            \infty,\quad \text{otherwise},
        \end{cases}
    \end{equation*}
    models a particle confined within a box of length $L$, with zero probability of existing outside the box.
    \item Finite potential well \cite{veselov1984finite}
        \begin{equation*}
    V(x)=
        \begin{cases}
            -V_0,\quad \text{if}~-1\leq x\leq 1,\\
            0,\quad \text{otherwise},
        \end{cases}
    \end{equation*}
    models a particle confined within a finite region of potential depth $V_0$.
    \item Harmonic Oscillator Potential \cite{marsiglio2009harmonic} $V(x)=\frac{1}{2}m\omega^2x^2$ represents a quantum harmonic oscillator, such as vibrations in molecules.
    \item Periodic potential (Kronig-Penney Model)\cite{greiner2001bose, seaman2005nonlinear} $V(x)=V_0\cos(\zeta x)$ represents a particle in a crystalline solid where the potential repeats periodically.
\end{itemize}

In the periodic potential, the frequency of the oscillations in the cosine function is 
$\zeta$. By increasing $\zeta$, the period of oscillation decreases, resulting in higher-frequency components. Therefore, the periodic potential could be a good choice for high-frequency-dominated problems.\\
For optical Lattices \cite{morsch2006dynamics, christodoulides2003discretizing}, the periodic potentials can be created by standing waves of laser light, with tunable depth $V_0$ and spacing $\frac{2\pi}{\zeta }$; For magnetic Lattices \cite{luttinger1951effect}, magnetic field variations create periodic trapping potentials.

The dispersive term (the Laplacian $\Delta$) acts to spread out wave packets. With the nonzero potential, the effect of the Laplacian is modulated so the dispersive behavior can be confined, scattered, or distorted. With large $\epsilon$, dispersion primarily governs the dynamics. In some cases, especially when the equation is strongly nonlinear or the potential is highly complex, the system can exhibit quasi-periodic or chaotic behavior, leading to irregular, high-frequency oscillations that are difficult to predict. 

We could use the Lie-Trotter time-splitting scheme \cite{zhang2025low} (see section \ref{pod v.s. podno} for the specific algorithm) to obtain the numerical solutions of the NLS equation.  
\subsection{Kadomtsev-Petviashvili (KP) equation\label{kp}}
We consider the KP equations given by
\begin{equation*}
        \partial_x\left(\partial_t u+u \partial_x u+\epsilon^2 \partial_{x x x} u\right)+\lambda \partial_{y y} u=0, \quad \lambda= \pm 1,
\end{equation*}
where $\epsilon>0$ is a scaling parameter. 
First, we would like to emphasize the oscillations caused by the linear part \cite{klein2007numerical}, i.e.,
\begin{equation}\label{linearKP}
    \partial_x\left(\partial_t u+\epsilon^2 \partial_{x x x} u\right)+\lambda \partial_{y y} u=0,\left.\quad u\right|_{t=0}=u_0(x, y).
\end{equation}
If we denote the Fourier transform by
\begin{equation*}
    \mathscr{F}[u](t, k_x, k_y) \equiv \hat{u}\left(t, k_x, k_y\right):=\int_{\mathbb{R}^2} u(t, x, y) \mathrm{e}^{-\mathrm{i}\left(k_x x+k_y y\right)} \mathrm{d} x \mathrm{d} y,
\end{equation*}
then the Fourier transform of the linear part can be written as:
\begin{equation}\label{fourier_linearKP}
    \mathrm{i} k_x \partial_t \hat{u}+\epsilon^2 k_x^4 \hat{u}-\lambda k_y^2 \hat{u}=0.
\end{equation}
The solution of \eqref{fourier_linearKP} is given by:
\begin{equation*}
\hat{u}\left(t, k_x, k_y\right)=\mathrm{e}^{-\mathrm{i} t\left(\lambda k_y^2 / k_x-\epsilon^2 k_x^3\right)} \hat{u}_0\left(k_x, k_y\right),
\end{equation*}
in the case that $k_x\neq 0$. The solution in the physical space therefore could be expressed as:
\begin{equation*}
    u(t, x, y)=\frac{1}{4 \pi^2}\big[\mathscr{F}^{-1}[\mathrm{e}^{-\mathrm{i} t \lambda k_y^2 / k_x}] * \mathscr{F}^{-1}[\mathrm{e}^{\mathrm{i} t \epsilon^2 k_x^3}]\big] * u_0(x, y).
\end{equation*}
The first inverse Fourier factor on the right-hand side gives the Green's function corresponding to the standard wave function. The second inverse Fourier factor on the right-hand side can be interpreted as a tempered distribution regarding the Airy function and the Dirac function:
\begin{equation*}
\mathscr{F}^{-1}[\mathrm{e}^{\mathrm{i} t \epsilon^2 k_x^3}]=\frac{1}{\left(3 t \epsilon^2\right)^{1 / 3}} \operatorname{Ai}\Big(\frac{x}{(3 t \epsilon^2)^{1 / 3}}\Big) \otimes \delta(y), \quad \forall t>0.
\end{equation*}
% \zw{you should use bracket?}
Thus we could conclude that the $y$-dependence of the solutions to the linear KP equation is completely nonoscillatory and has a rather slow algebraic decay as $|y|\rightarrow \infty$; the oscillatory behavior of the linear KP equation is governed by the Airy function. In our scaling, this yields oscillations with wave length $\mathcal{O}(\epsilon^{\frac{2}{3}})$ in the $x$-direction.

The nonlinear term then also plays a crucial role in shaping the behavior of the solutions \cite{cheng2014interactions,grava2018numerical}. The interplay between nonlinearity and dispersion gives rise to a range of phenomena, including solitons and shocks, leading to highly oscillatory solutions.

We use the ETDRK4 formula \cite{cox2002exponential} to perform the time integration and generate datasets for experiments on PODNO.
\end{appendix}
\end{document}